\documentclass[12pt,twoside]{report}
\usepackage[T1]{fontenc}

\usepackage{fancyhdr}

\usepackage{a4}
\usepackage{amsmath}
\usepackage{amssymb}
\usepackage{amsthm}
\usepackage{verbatim}

\usepackage[all]{xy}
\usepackage[dvips]{graphicx}

\usepackage[utf8x]{inputenc} 

\usepackage[polish,british]{babel}

\def\gotm{\mathfrak{m}}

\def\Inv{\mathcal{INV}}

\def\cDeg{\mathcal{DEG}}
\def\tr{\operatorname{tr}}
\def\Xdeg{X_{\operatorname{deg}}}
\def\Xinv{X_{\operatorname{inv}}}
\def\Ldot{\mathbf{L}^{\bullet}}
\def\thetadot{\theta^{\bullet}}

\newenvironment{prf}
{\medskip\par\noindent{\bf Proof.}}
{\nopagebreak\par\rightline{$\Box$} \medskip}
\def\noprf{\nopagebreak\par\rightline{$\Box$} \medskip}

\numberwithin{equation}{chapter}
\newtheorem{theo}[equation]{Theorem}
\newtheorem{lem}[equation]{Lemma}
\newtheorem{cor}[equation]{Corollary}
\newtheorem{con}[equation]{Conjecture}

\newtheorem*{defin}{Definition}

\newtheorem{ex}[equation]{Example}
\newtheorem{rem}[equation]{Remark}
\newtheorem{prop}[equation]{Proposition}

\def\ideal{\vartriangleleft}

\def\uD{\mathrm{D}}
\def\ud{\mathrm{d}}
\def\sspan{\operatorname{span}}
\def\pd{\partial}
\def\Sp{\mathbf{Sp}}
\def\cSp{\mathbf{cSp}}

\def\Sl{\mathbf{SL}}

\def\GL{\mathbf{GL}}
\def\Gl{\mathbf{GL}}
\def\exp{\operatorname{exp}}

\def\Id{\operatorname{Id}}
\def\I{\mathcal{I}}

\def\Sym{\operatorname{Sym}}
\def\conv{\operatorname{conv}}
\def\pf{\operatorname{Pfaff}}
\def\bS{\mathbb{S}}

\def\mono{\hookrightarrow}
\def\epi{\twoheadrightarrow}

\def\gotg{\mathfrak{g}}

\def\gotsp{\mathfrak{sp}}
\def\gotcsp{\mathfrak{csp}}
\def\gotasp{\mathfrak{wsp}}
\def\gotsl{\mathfrak{sl}}

\def\gotgl{\mathfrak{gl}}
\def\gotm{\mathfrak{m}}

\def\Aut{\operatorname{Aut}}
\def\gotaut{\mathfrak{aut}}
\def\autinf{\gotaut^{\operatorname{inf}}}
\def\autinfF{\autinf_F}
\def\Authol{\Aut^{hol}}
\def\Autalg{\Aut^{alg}}
\def\Autdot{\Aut^{\bullet}}

\def\half{\frac{1}{2}}

\def\C{{\mathbb C}}
\def\H{{\mathbb H}}
\def\N{{\mathbb N}}

\def\Q{{\mathbb Q}}
\def\R{{\mathbb R}}
\def\Z{{\mathbb Z}}
\def\F{{\mathbb F}}

\def\P{{\mathbb P}}
\def\A{{\mathbb A}}

\def\ra{\rightarrow}
\def\lra{\longrightarrow}

\def\ccB{{\cal B}}

\def\ccF{{\cal F}}
\def\ccG{{\cal G}}

\def\ccO{{\cal O}}

\def\ccR{{\cal R}}

\def\ccS{{\cal S}}

\def\cHom{{\mathcal H}om}

\def\Proj{\operatorname{Proj}}
\def\Spec{\operatorname{Spec}}
\def\SheafySpec{\ccS{}pec}

\def\Cl{\operatorname{Cl}}
\def\End{\operatorname{End}}

\def\Pic{\operatorname{Pic}}
\def\Sing{\operatorname{Sing}}

\def\codim{\operatorname{codim}}

\def\diag{\operatorname{diag}}
\def\id{\operatorname{id}}
\def\im{\operatorname{im}}

\def\rk{\operatorname{rk}}

\def\red{\operatorname{red}}
\def\wt{\operatorname{wt}}

\def\set#1{\left\{#1\right\}}

\def\Wedge#1{\textstyle{\bigwedge\nolimits}^{\! #1}}

\renewcommand{\labelenumi}{\theenumi}
\renewcommand{\thechapter}{\Alph{chapter}}

\def\thetitle{%%%%%%% Type the title of the thesis here: 
Algebraic Legendrian varieties}
\def\theauthor{%%%%%% Type the name of the author here:
Jaros\l{}aw Buczy\'n{}ski}

\author{\theauthor}
\title{\thetitle}

\date{20 October 2009}

\begin{document}
\renewcommand{\headrulewidth}{0pt} 
\renewcommand{\headheight}{15pt} 

\lfoot[\fancyplain{\thepage}{\thepage}]{}
\rfoot[\fancyplain{}{}]{\fancyplain{\thepage}{\thepage}}
\cfoot{\fancyplain{}{}}

\rhead{\fancyplain{}{}}
\lhead{\fancyplain{}{}}
\chead[\fancyplain{}{\textsc{\theauthor}}]{\fancyplain{}{\textsc{\thetitle}}}

\maketitle

\begin{abstract}
\thispagestyle{plain}

Real Legendrian subvarieties are classical objects of differential
geometry and classical mechanics and they have been studied since antiquity 
(see \cite{arnold}, \cite{slawianowski} and references therein).
However, complex Legendrian subvarieties 
are much more rigid and have more exceptional properties. 
The most remarkable case is 
the Legendrian subvarieties of projective space 
and prior to the author's research only few smooth examples of these were known
(see \cite{bryant}, \cite{landsbergmanivel04}).
Strong restrictions on the topology of such varieties have been found and studied 
by Landsberg and Manivel (\cite{landsbergmanivel04}).
 
This dissertation reviews the subject of Legendrian varieties and extends some of recent results. 

The first series of results is related to the  automorphism group of any Legendrian subvariety 
in any projective contact manifold. 
The connected component of this group (under suitable minor assumptions) 
is completely determined by 
the sections of the distinguished line bundle on the contact manifold 
vanishing on the Legendrian variety.
Moreover its action 
preserves the contact structure.
The relation between the Lie algebra tangent to automorphisms 
and the sections is given by an explicit formula (see also \cite{lebrun}, \cite{beauville}).
This summarises and extends some earlier results of the author. 

The second series of results is devoted to finding new examples of smooth 
Legendrian subvarieties of projective space.
The examples found by other researchers were
some homogeneous spaces,
many examples of curves
and a family of surfaces birational to some K3 surfaces.
Further the author found a couple of other examples including 
 a smooth toric surface and 
 a smooth quasihomogeneous Fano 8-fold.
Finally, the author proved both of these are special cases of a very general construction:
a general hyperplane section of a smooth Legendrian variety, after a suitable projection, 
is a smooth Legendrian variety of smaller dimension.
We review all of those examples and also add another 
infinitely many new examples in every dimension, with various 
Picard rank,
canonical degree,
Kodaira dimension 
and other invariants.

The original motivation for studying complex Legendrian varieties comes from 
a 50 year old problem of giving compact examples of quaternion-K\"{a}hler manifolds 
(see \cite{berger}, \cite{lebrun_salamon}, \cite{lebrun} and references therein).
Also Legendrian varieties are related to some algebraic structures 
(see \cite{mukai}, \cite{landsbergmanivel01}, \cite{landsbergmanivel02}).
A new potential application to classification of smooth varieties with smooth dual arises 
in this dissertation.

\medskip
\footnotesize{
\noindent\textbf{keywords:} \\
Legendrian variety, 
complex contact manifold, 
automorphism group;

}
\end{abstract}

\tableofcontents

\vfill

\section*{Acknowledgements}

The monograph is built on the author's PhD thesis \cite{jabu_dr}.

The author was supported by 
the research project N20103331\slash{}2715 
funded by Polish financial means for science from 2006 to 2008.

The author would like to thank especially 
his patient advisor Jaros\l{}aw Wi\'s{}\-nie\-wski for his comments, support
and for answering numerous questions.
Parts of the thesis were written while the author enjoyed the hospitality 
of Insong Choe at  KIAS (Korea Institute for Advanced Study), 
Mark Gross at University of California, San Diego and 
Joseph M.~Landsberg at Texas A\&M University.
The author is grateful for their 
invitation, financial support and for creating a stimulating
atmosphere.
Also the author acknowledges many enlightening discussions with
Michel Brion,
Stephen Coughlan,
Zbigniew Jelonek,
Grzegorz Kapustka, 
Micha\l{} Kapustka,
Micha\l{} Krych,
Joseph M.~Landsberg,
Adrian Langer,
Laurent Manivel,
Sung Ho Wang and
Andrzej Weber.
Finally, the author is grateful to the referees of the thesis: 
S\l{}awomir Cynk,
Laurent Manivel,
Henryk Żołądek
 and the anonymous referee of this mongraph 
--- their comments have helped to improve the presentation.

\chapter{Introduction}
\rhead[\fancyplain{}{}]{\fancyplain{}{\footnotesize\textbf{Chapter \thechapter}}}

\section{State of art}\label{section_state_of_art}

In this monograph we study algebraic and geometric properties of complex Legendrian subvarieties.
The main motivation for our research comes from the classification problem of contact 
Fano manifolds\footnote{
A complex manifold $Y^{2n+1}$ is called \textbf{a contact manifold}
if there exists a rank $2n$ vector subbundle $F \subset TY$ of the tangent bundle,
such that the map $F \otimes F \lra TY/F$ determined by the Lie bracket is nowhere degenerate
(see Chapter~\ref{chapter_contact} for more details).
A projective manifold is \textbf{Fano} if its anticanonical line bundle is ample.
}.

\subsection{Contact manifolds and quaternion-K\"{a}hler manifolds} \label{intro_contact_and_qK}

There are two known families of complex projective contact manifolds:
 projectivisations of cotangent bundles to projective manifolds 
 (see Example~\ref{example_projectivised_cotangent})
 and certain homogeneous spaces, namely
 the unique closed orbit of the adjoint action of a simple Lie group $G$ on $\P(\gotg)$,
 where the $\gotg$ is the Lie algebra of $G$ 
 (see \S\ref{section_examples_of_symplectic} and  Example \ref{example_adjoint_orbit_is_contact}
 and for more details). 
 These orbits are also called \textbf{the adjoint varieties}.
 The following theorem summarises works of Demailly \cite{demailly} 
 and Kebekus, Peternell, Sommese and Wi\'s{}niewski \cite{4authors}:
\begin{theo}\label{theorem_KPSW_D}
 If $Y$ is a complex projective contact manifold, then either
 $Y$ is a projectivisation of the cotangent bundle to some projective manifold $M$ or 
 $Y$ is a Fano manifold with second Betti number $b_2=1$.
\end{theo}

The following conjecture would be an important classification result in algebraic geometry 
and it claims that known examples are all existing examples:

\begin{con}[LeBrun, Salamon] \label{con_contact}
If $Y^{2n+1}$ is a Fano complex contact manifold, then $Y$ is 
an adjoint variety.
\end{con}

This conjecture originated with a famous problem in Riemannian geometry.
In 1955 Berger \cite{berger} gave a list of all possible holonomy groups%
\footnote{Given an $m$-dimensional Riemannian manifold $M$,
\textbf{the holonomy group of $M$} is the subgroup of orthogonal group $\mathbf{O}(T_x M)$ 
generated by parallel translations along loops through $x$.}
of simply connected Riemannian manifolds. 
The existence problem for all the cases has been solved locally. 
Compact non-homogeneous examples with  most of the possible holonomy groups were constructed,
for instance the two exceptional cases $G_2$ and $Spin_7$ were constructed by D.~Joyce 
--- see an excellent review on the subject \cite{joyce}. 
Since then all the cases from Berger's list have been illustrated with compact non-homogeneous examples, 
with the unique exception of the quaternion-K\"{a}hler manifolds\footnote{
A Riemannian $4n$-dimensional manifold $M$ is called \textbf{quaternion-K\"{a}hler} if its 
holonomy group is a subgroup of $\Sp(1) \times \Sp(n) / \Z_2$.}.
Although there exist non-compact, non-homogeneous examples, 
it is conjectured that the compact quaternion-K\"{a}hler manifolds must be homogeneous
(see \cite{lebrun} and references therein).

\begin{con}[LeBrun, Salamon]\label{con_qK}
Let $M^{4n}$ be a compact quaternion-K\"{a}hler manifold.
Then $M$ is a homogeneous symmetric space 
(more precisely, it is one of the Wolf spaces  --- see \cite{wolf}).
\end{con}

The relation between the two conjectures is given by the construction of a twistor space $Y$,
an $S^2$-bundle of complex structures on tangent spaces to a quaternion-K\"{a}hler manifold $M$.
If $M$ is compact, it has positive scalar curvature, and then $Y$ has a natural complex structure 
and 
is a contact Fano manifold with a K\"{a}hler-Einstein metric. 
In particular, the twistor space of a Wolf space is an adjoint variety.
Hence Conjecture~\ref{con_contact} implies Conjecture~\ref{con_qK}.
Conversely, LeBrun \cite{lebrun} observed that if $Y$ is a contact Fano manifold with 
K\"{a}hler-Einstein metric, then 
 it is a twistor space of a quaternion-K\"{a}hler manifold. 

A number of attempts have been  undertaken to prove the above conjectures. 
They were proved in  low dimension: for 
$n=1$ by N.~Hitchin \cite{hitchin} and Y.~Ye \cite{ye}, 
$n=2$ by Y.S.~Poon and S.M.~Salamon \cite{poon_salamon} and S.~Druel \cite{druel}
and Conjecture~\ref{con_qK} for $n=3$ by H.\&R.~Herrera \cite{2herreras}. 
Moreover A.~Beauville, J.~Wi\'s{}niewski, S.~Kebekus, T.~Peternell, A.~Sommese, J.P.~Demailly, C.~LeBrun, J-M.~Hwang 
and many other researchers have worked on this problem. 

\subsection{Lines on contact Fano manifold}

If $Y^{2n+1}$ is a contact complex manifold,
then  a subvariety $X\subset Y$ of pure dimension $n$ is \textbf{Legendrian}
if it is maximally $F$-integrable --- see Chapter~\ref{chapter_contact} for the details.

Let $Y^{2n+1}$ be a contact Fano manifold not isomorphic to a projective space.
A rational curve $C\subset Y$ is \textbf{a contact line}
if its intersection with the anticanonical line bundle is minimal possible,
i.e.~equal to $n+1$.
Contact lines on $Y$ are an instance of minimal rational curve on uniruled manifold 
and they are extensively studied by numerous researchers.
Wi\'s{}niewski \cite{wisniewski} and Kebekus \cite{kebekus_lines1}, \cite{kebekus_lines2}
have studied geometric properties of contact lines. 
The following theorem is due to Kebekus \cite{kebekus_lines2}, but some parts of it were known before:

\renewcommand{\theenumi}{(\roman{enumi})}
\begin{theo}\label{theorem_kebekus_lines}
  With $Y$ as above choose a point $y \in Y$. 
  Let 
  \[
    H_y \subset RatCurves^n(Y)
  \]
  be the subscheme parametrising contact lines through $y$,
  let $C_y \subset Y$ be the locus of these lines (i.e.~the subset swept out by contact lines through $y$)
  and let $X_y \subset \P(T_y Y)$ be the projectivised tangent cone to $C_y$ at $y$ 
  (so that $X_y$ is the set of tangent directions at $y$ to lines through $y$). 
  Then:
  \begin{enumerate}
   \item  
       $H_y$ is projective and $C_y$ is closed in $Y$.
   \item
       $C_y \subset Y$ is a Legendrian subvariety 
       and $X_y \subset \P(F_y) \subset \P(T_y Y)$
       and $X_y$ is a non-degenerate Legendrian subvariety in $\P(F_y)$.
   \item 
       If $y$ is a general point, then 
       $H_y$ is isomorphic to $X_y$ 
       (in other words every line through $y$ is determined by its tangent direction)
       and $H_y$ is smooth.
   \item 
       If $y$ is a general point, then 
       $C_y$ is isomorphic to the projective cone
       $\widetilde{X}_y \subset \P(F_y \oplus \C)$
       over $X_y \subset \P(F_y)$. 
       This isomorphism maps the contact line through $y$ tangent to $x \in X_y$
       onto the line $\P^1 = \P(x \oplus \C) \subset \P(F_y \oplus \C)$. 
       In particular two distinct lines through $y$ do not intersect except at $y$.
       Moreover the distinguished line bundle $L = T Y / F$
       restricted to $C_y$ is identified with $\ccO_{\widetilde{X}_y} (1)$ via this isomorphism.
  \end{enumerate}
\end{theo}

Note that in \cite[Thm~1.1]{kebekus_lines2} 
Kebekus also claims that a certain irreducibility holds for $H_y$.
However it was observed by Kebekus himself together with the author that 
there is a gap in the proof\footnote{%
This gap is on page 234 in Step~2 of proof of Proposition~3.2 
where Kebekus claims to construct ``a well defined family of cycles''.
This is not necessarily a well defined family of cycles:
condition (3.10.4) in \cite[\S{}I.3.10]{kollar_book_rational_curves}
is not necessarily satisfied.
As a consequence the map $\Phi$ is not necessarily regular at non-normal points of $D^0$.}.

Thus contact lines through a general point 
behave very much like ordinary lines in a projective space. 
Moreover their geometry is described by a nondegenerate smooth Legendrian subvariety $X:=X_y$ in $\P^{2n-1}$.
If $Y$ is  one of the adjoint varieties, then $X$ will be a homogeneous Legendrian subvariety 
called \textbf{a subadjoint variety} (see \cite{landsbergmanivel04}, \cite{mukai}).
Proving that there is an embedding of $Y$ 
into a projective space which maps contact lines to ordinary lines 
would imply Conjecture~\ref{con_contact}.
Moreover it is proved by Hong \cite{hong}, 
that if $X$ is homogeneous, then so is $Y$. 
Therefore contact lines and particularly the Legendrian varieties determined by them 
are important objects,
useful in the study of Conjecture~\ref{con_contact}.

\begin{table}[hbt]
{\small\centering
\begin{tabular}{|p{0.1\textwidth}|p{0.1\textwidth}|p{0.2\textwidth}|p{0.2\textwidth}|p{0.22\textwidth}|}
\hline
  Lie group  & Type& Contact  manifold $Y^{2n+1}$     & Legendrian variety $X^{n-1}$ & Remarks \\
\hline
&&&&\\
$\textrm{SL}_{n+2}$  & $ A_{n+1}  $              & $\P(T\P^{n+1})$      & $\P^{n-1} \sqcup \P^{n-1}$ \mbox{$\subset \P^{2n-1}$}
& $b_2(Y) =2 $    \\
\hline
$\textrm{Sp}_{2n+2}$ & $ C_{n+1}$              & $\P^{2n+1}$          &
 $\emptyset\subset \P^{2n-1}$                & $Y$ does not have any contact lines\\
\hline
$\textrm{SO}_{n+4}$  & $ B_{\frac{n+3}{2}}$  or $D_{\frac{n+4}{2}}$  & $Gr_O(2,n+4)$        & 
$\P^1 \times Q^{n-2}$ \mbox{$ \subset \P^{2n-1}$}   & $Y$ is the Grassmannian  of projective lines on a quadric $Q^{n+2}$\\
\hline
                           & $ G_2$                  &Grassmannian of special lines on $Q^5$ & $\P^1 \subset \P^3$  & $X$ is the twisted cubic curve\\
\hline
                           & $ F_4$                  & an $F_4$ variety   & $Gr_L(3,6)$  \mbox{$ \subset \P^{13}$}    &   \\
\hline
                           & $ E_6$                  & an $E_6$ variety   & $Gr(3,6) \subset \P^{19}$        &   \\
\hline
                           & $ E_7$                  & an $E_7$ variety   & $\mathbb{S}_6 \subset \P^{31}$  &$X$ is the spinor variety
\\
\hline
                           & $ E_8$                  & an $E_8$ variety   & the $E_7$ variety $\subset \P^{55}$  \includegraphics[width=0.17\textwidth]{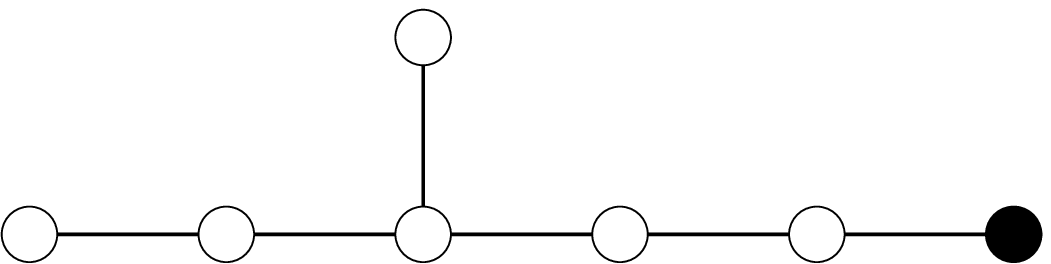}
&   \\
\hline
\end{tabular}
\caption{
Simple Lie groups together with the corresponding adjoint variety $Y$ and
its variety of tangent directions to contact lines: the subadjoint variety $X$
(listed in details also in \S\ref{section_intro_legendrian}).}
\label{table_contacts}}
\end{table}

\subsection{Legendrian subvarieties of projective space}
\label{section_intro_legendrian}

\renewcommand{\theenumi}{\arabic{enumi})}
Prior to the author's research  the following were the only known examples of   
smooth Legendrian subvarieties of projective space
(see \cite{bryant}, \cite{landsbergmanivel04}):
\begin{enumerate}
  \item 
    linear subspaces;
  \item
    some homogeneous spaces called \textbf{subadjoint varieties} 
    (see Table~\ref{table_contacts}):
    the product of a line and  a quadric $\P^1 \times Q^{n-2}$
    and five exceptional cases:
    \begin{itemize}
      \item
        twisted cubic curve  $\P^1 \subset \P^3$,
      \item
        Grassmannian $Gr_L(3,6) \subset \P^{13}$ of Lagrangian subspaces in $\C^6$,
      \item
        full Grassmannian  $Gr(3,6)\subset \P^{19}$,
      \item
        spinor variety $\mathbb{S}_6\subset \P^{31}$
         (i.e.~the homogeneous $\mathbf{SO}(12)$-space
         parametrising the vector subspaces of
         dimension $6$ contained
         in a non-degenerate quadratic cone in $\C^{12}$) and
      \item
        the 27-dimensional $E_7$-variety in $\P^{55}$ corresponding to the marked root:
        \includegraphics[width=0.2\textwidth]{subadjoint_E7};
    \end{itemize}
  \item
    every smooth projective curve admits a Legendrian embedding in $\P^3$ (see \cite{bryant});
  \nopagebreak \item \nopagebreak
    a family of smooth surfaces birational to the Kummer $K3$-surfaces (see \cite{landsbergmanivel04}).
\end{enumerate}

 The subadjoint varieties are expected 
 to be the only homogeneous Legendrian subvarieties in $\P^{2n-1}$ (for $n > 2$).
 Many partial results are known. 
 For instance the following Theorem is traditionally known:
 \begin{theo}
   If $X \subset \P(V)$ is homogeneous and Legendrian and the Legendrian embedding is equivariant 
   (i.e.~the automorphisms of $X$ extend to linear automorphisms of $\P(V)$),
   then $X$ is one of the subadjoint varieties.  
 \end{theo}
 To the author's best knowledge the statement of the theorem has not been explicitly written down before,
 but according to J.M.~Landsberg, 
 already E.~Cartan have known which homogeneous spaces are Legendrian. 
 Also the theorem immediatelly follows from various resources.
 As one instance, the subadjoint varieties are the only smooth irreducible Legendrian varieties
 whose ideal is generated by quadratic polynomials 
 (see \cite{jabu06}).
 But it is a well known theorem of Kostant,
 that homogeneous space in its equivariant embedding has the ideal generated by quadrics
 (see for instance \cite[\S10.6.6]{procesi_book}).

 Also Landsberg and Manivel proved \cite[Thm.~11]{landsbergmanivel04}:
 \begin{theo}\label{theorem_LM_legendrian_homogeneous_b2_one_then_subadjoint}
  Let $X = G/P \ne \P^1$ be a homogeneous space with Picard number one.
  Suppose that $X$ admits a Legendrian embedding into $\P^{2n-1}$,
   not necessarily equivariant a priori. 
  Then $X$ is subadjoint.
  In particular, the Legendrian embedding is the equivariant subadjoint embedding.
 \end{theo}

The subadjoint varieties are strongly related to the group they arise from. 
In Landsberg and Manivel \cite{landsbergmanivel02} 
use the subadjoint varieties to reconstruct the adjoint variety they arise from.
And thus the subadjoint varieties are one of essential ingridients 
of their proof of the classification of simple Lie groups by means of projective geometry only.
Also Mukai \cite{mukai} relates the subadjoint Legendrian varieties 
with another algebraic structure: simple Jordan algebras. 
In \cite{landsbergmanivel01} the authors give a uniform description 
of the exceptional cases (arising from $F_4$, $E_6$, $E_7$ and $E_8$).

For an arbitrary Legendrian subvariety in projective space,
 Landsberg and Manivel proved certain restrictions on its topology,
 the simpliest of which is:
\begin{theo}
   Let $X \subset \P^{2n-1}$ be a smooth Legendrian variety, 
   let $h \in H^2 (X, \Q)$ be the hyperplane class 
    and $c_i \in H^2 (X, \Q)$  be the $i$-th Chern class of $X$.
   Then: 
   \[
     {c_1}^2 − 2c_2 = 2 h c_1 − n h^2.
   \]
   This gives in particular the following restrictions on $X$:
   \begin{itemize}
     \item $X$ cannot be an abelian variety or parallelisable variety.
     \item If $X \simeq \P^{n-1}$ for $n > 2$, then $X \subset \P^{2n-1}$ 
           is a linear subspace.
     \item Suppose that $X \simeq Y \times Z$.
           Then $X = \P^1 × Q^{n−2}$ and the Legendrian embedding is the Segre embedding.
     \item If $X \simeq \P(E)$, the total space of a $\P^p$ bundle over some smooth variety $Y$,
           then $p =1$. 
   \end{itemize}
\end{theo}
Note that the statement in \cite[Prop.~3]{landsbergmanivel04} suggests that the Chern class identities 
hold in $H^{\bullet} (X, \Z)$, yet their proof uses cohomologies with rational coefficients 
(they use $e^h$ and Chern character $ch(TX)$). 

Further restrictions on topology of $X$ where derived if $\dim X =2$, 
see \cite{landsbergmanivel04} for details.

Also the authors there proved \cite[Prop.~17(2)]{landsbergmanivel04}:

\begin{prop} 
  Let $X\subset\P(V)$ be a non-degenerate Legendrian variety. Then:
  \begin{itemize}
   \item the secant variety $\sigma(X)$ fills out $\P(V)$ 
         (i.e.~it is of expected dimension).
   \item if $X$ is in addition smooth,
         then the tangent variety $\tau(X)\subset\P(V)$ 
         and the dual variety
         $X^* \subset\P (V^*)$ are hypersurfaces 
         isomorphic via $\tilde{\omega} : V \ra V^*$ 
  \end{itemize}
\end{prop}
  See \S\ref{section_definition_of_dual_variety} for definitions of secant, tangent and dual varieties.
  See \S\ref{section_bilinear_forms_and_matrices} for definition of $\tilde{\omega}$.
  Here the symplectic form $\omega$ is the form associated with the contact structure on $\P(V)$
  (see Example~\ref{example_projective_space_as_contact}).
  The statement about secant variety is contained in the proof \cite[Prop.~17(2)]{landsbergmanivel04}.

For tangent and dual varieties, 
we note that original formulation in \cite{landsbergmanivel04} 
omits the smoothness assumption. 
If $X$ is not smooth however, their proof does not work
 and the decomposable Legendrian varieties (see \S\ref{definition_decomposable}) are counterexamples. 
In the proof the authors freely interchange the tangent variety 
$\tau(X)$ (which by definition is the union of the limits of secants 
through two points approaching a third fixed point) 
and the closure of the union of embedded tangent spaces at smooth points. 
These are the same for $X$ smooth.
The tangent variety $\tau(X)$ is indeed a hypersurface in the secant variety $\sigma(X)$.
The closure of the embedded tangent spaces at smooth points is indeed isomorphic
to $X^*$.

\section{Topics of the dissertation}

This dissertation addresses two complementary problems regarding Legendrian varieties:
\begin{itemize}
  \item
     write explicit restrictions on the properties of Legendrian varieties;
  \item 
     give examples of smooth Legendrian varieties.
\end{itemize}
The first problem we address by giving a very precise understanding of
the embedded automorphism group of a Legendrian variety.
The second problem is solved by proving that a general hyperplane section 
of a smooth Legendrian variety admits a Legendrian embedding.

In  \cite{jabu_mgr}, we prove that the quadratic part 
of the ideal of a Legendrian subvariety $X$ of projective space $\P^{2n-1}$ produces 
a connected subgroup of projective automorphisms of $X$.
In \cite{jabu06} we improve this result by observing that this group is actually 
the maximal connected subgroup of automorphisms of the contact structure on $\P^{2n-1}$ 
preserving the Legendrian subvariety. 
This relation is explicitely defined:
\begin{theo} \label{theorem_ideal_and_group}
  Let $X \subset \P(V)$ be a Legendrian subvariety.
  Consider the following map $\rho$:
  \[
    H^0(\ccO_{\P(V)}(2)) \simeq \Sym^2 V^* \ni  \ q = \bigl(x \mapsto x^T M(q) x\bigl) 
    \stackrel{\rho}{\longmapsto} 2J \cdot M(q) \ \in \gotsp(V).
  \]
  where $M(q)$ is the matrix of $q$ and 
  $J=M(\omega)$ is the matrix of the symplectic form on $V$ associated with the contact structure on $\P(V)$.
  Let $\widetilde{\I}(X)_2 \subset \Sym^2 V^*$ be the vector space of quadrics containing $X$. Then:
  \begin{itemize} 
    \item
      $\rho(\widetilde{\I}(X)_2)$ is a Lie subalgebra of $\gotsp(V)$ tangent to a closed subgroup 
      $$
        \overline{\exp\Big(\rho\big(\widetilde{\I}(X)_2\big)\Big)} < \Sp(V).
      $$
    \item
      We have the natural action of $\Sp(V)$ on $\P(V)$. 
      The group 
      $\overline{\exp\Big(\rho\big(\widetilde{\I}(X)_2\big)\Big)}$
      is the maximal connected subgroup in $\Sp(V)$ 
      which under this action preserves $X \subset \P(V)$.
  \end{itemize}
\end{theo}

For proof see \cite[Cor.~4.4, Cor.~5.5, Lem.~5.6]{jabu06}.

In the present dissertation we extend this result further. 
Firstly, we replace projective space $\P(V)$ with an arbitrary contact manifold $Y$:
\begin{theo}\label{theorem_introduction_contact_automorphisms_preserving_legendrian}
  Let $Y$ be a compact contact manifold and let $X \subset Y$ be a Legendrian subvariety.  
  Then the connected component of the subgroup of $\Aut(Y)$ that preserves both 
  the contact structure and $X \subset Y$, is completely
  determined by those sections of a distinguished line bundle $L$ on $Y$ that vanish on $X$.
\end{theo}
See Corollary~\ref{corollary_Y_projective_autF_eq_I} for proof. 
In Theorem~\ref{theorem_contact_automorphisms} (quoted from \cite{beauvillefano})
we also make explicit how to obtain an infinitesimal automorphism of $Y$ from a given section of $L$ 
(the analogue of the map $\rho$ in Theorem~\ref{theorem_ideal_and_group}).

Secondly, we try to remove the assumption that the automorphisms preserve the contact structure.
By applying the results of \cite{lebrun} and \cite{kebekus_lines1}
on the uniqueness of contact structures we can deal with this problem for
most projective contact Fano manifolds (see Corollary~\ref{corollary_Y_projective_autF_eq_I}).
The remaining cases are the projectivised cotangent bundles and the projective space.
The first case is not very interesting, as all the Legendrian subvarieties are 
classified for these contact manifolds (see Corollary~\ref{cor_Legendrian_in_cotangent}).
On the other hand the case of projective space is the most important and interesting. 
It is described precisely in Chapter~\ref{chapter_automorphisms}.
We present there the following theorem originally proved in \cite{jabu_toric}:

\begin{theo}\label{corollary_that_conjecture_is_true_for_smooth}
If $X\subset \P^{2n-1}$ is a smooth irreducible Legendrian subvariety which is not a linear subspace
and $G< \P\Gl_{2n}$ is a connected subgroup preserving $X$,
then the action of $G$ on $\P^{2n-1}$ necessarily preserves the contact structure.
Thus in this case the group $\overline{\exp\Big(\rho\big(\widetilde{\I}(X)_2\big)\Big)}$
from Theorem~\ref{theorem_ideal_and_group} is also the maximal connected subgroup in $\Sl(V)$ 
which  preserves $X \subset \P(V)$.
\end{theo}

Our methodology for finding new examples of smooth Legendrian subvarieties is the following.
We pose questions of classification of smooth Legendrian varieties satisfying certain 
additional conditions. For instance, we assume that the variety is toric 
(see Chapter~\ref{chapter_toric}, which follows \cite{jabu_toric}):
\begin{theo}\label{theorem_classification_of_smooth_toric_legendrian}
Every smooth toric Legendrian variety of dimension at least~$2$ in a projective space,
whose embedding is torus equivariant is isomorphic to one of the following:
\begin{itemize}
\item a linear subspace,
\item  $\P^1 \times Q_1 \subset \P^5$,
\item  $\P^1 \times Q_2 \simeq \P^1 \times \P^1 \times \P^1 \subset \P^7$
\item  or $\P^2$ blown up in three non-colinear points.
\end{itemize}
\end{theo}
For proofs see Corollaries~\ref{smooth_toric_surfaces} and \ref{smooth_toric_varieties}.
The linear subspace is not really interesting, 
the products $\P^1 \times Q_1$ and $\P^1 \times Q_2$  are well known 
(see \S\ref{section_intro_legendrian}).
The last case of blow up was an original example of \cite{jabu_toric}.

We also classify those smooth Legendrian varieties,
which are contained in a specific $F$-cointegrable variety
(see Chapter~\ref{chapter_sl}, which follows \cite{jabu_sl}).
In this way we produce a few new smooth examples including a quasihomogeneous Fano 8-fold:
\begin{theo}\label{theorem_hyperplane_section_of_G_3_6}
  A general hyperplane section of $G(3,6)$ admits a Legendrian embedding into $\P^{17}$.
\end{theo}

Finally following \cite{jabu_hyperplane} we generalise this 
and prove that a general hyperplane section of a smooth Legendrian variety 
admits a Legendrian embedding into a smaller projective space:

\begin{theo}
 \label{theorem_hyperplane}
Let $X\subset \P(V)$ be an irreducible Legendrian subvariety, which is
 smooth or has only isolated singularities. 
Then a general hyperplane section of $X$ admits a Legendrian embedding
 into a projective space of appropriate dimension via a specific
 subsystem of the linear system $\ccO(1)$. 

More generally, assume $X\subset \P(V)$ is an irreducible Legendrian
 subvariety with singular locus of dimension $k$ and $H\subset \P(V)$ is
a general hyperplane. Then there exists a variety $\widetilde{X}_H$ whose
 singular locus has dimension at most $k-1$ and which has an open subset
  isomorphic to the smooth locus of $X\cap H$ such that $\widetilde{X}_H$ admits a Legendrian
 embedding.
\end{theo}

The specific linear system and construction of $\widetilde{X}_H$
is described in \S\ref{section_construction} and there we
prove that the resulting variety is Legendrian.
The proof that for a general 
section the result has the required smoothness property is
presented in \S\ref{section_smooth}.

This simple observation has quite strong consequences.
Many researchers, including Landsberg, Manivel, Wi\'s{}niewski, Hwang
and the author of this dissertation,
believed that the structure of smooth Legendrian subvarieties in
projective space had to  be somehow rigid at least in higher dimensions.
So far the only non-rational examples known were in dimensions 1 and 2 
(see \S\ref{section_intro_legendrian}) and
these were also the only known to come in families. 
Already by a  naive application of the theorem to the subadjoint varieties we get many more
examples with various properties:

\begin{ex}\label{many_examples}
The following smooth varieties  and families of smooth varieties admit Legendrian embedding:
\begin{itemize}
\item[(a)]
a family of $K3$ surfaces of genus 9;
\item[(b)]
three different types of surfaces of general type;
\item[(c)]
some Calabi-Yau 3-folds, some Calabi-Yau 5-folds and some Calabi-Yau 9-folds;
\item[(d)]
some varieties of general type in  dimensions 3, 4 (two families for every dimension), 
5,6,7 and 8 (one family per dimension);
\item[(e)]
some Fano varieties, like the blow up of a quadric $Q^n$ in a codimension 2
	  hyperplane section $Q^{n-2}$, 
          a family of Del Pezzo surfaces of degree 4 and others;
\item[(f)]
infinitely many non-isomorphic, non-homogeneous Legendrian varieties in every dimension
	  arising as a codimension $k$ linear section of $\P^1\times Q^{n+k}$.
\end{itemize}
\end{ex}

Example (a) agrees with the prediction of \cite[\S 2.3]{landsbergmanivel04}.
Examples (b) and (d) give a partial answer to the  question of  a possible Kodaira
dimension of a Legendrian variety (also see \cite[\S 2.3]{landsbergmanivel04}).
Example (f) is a counterexample to the naive expectation that Legendrian
variety in a sufficiently high dimension must be homogeneous.

We also note that the previous non-homogeneous examples also arise in this way.
Example (e) for $n=2$ is the toric example described in Example~\ref{example_plane_blown_up}.
Hyperplane sections of $Gr(3,6)$, $Gr_L(3,6)$, $\bS_6$ are studied in
more details in Chapter~\ref{chapter_sl}. 
Also non-homogeneous examples of the other authors, Bryant \cite{bryant},
Landsberg and Manivel \cite{landsbergmanivel04} can be reconstructed by
Theorem~\ref{theorem_hyperplane} from some varieties with only
isolated singularities (see \S\ref{section_extend}).

In this monograph we also demonstrate
a more refined construction, using the decomposable Legendrian varieties, 
which also uses Theorem~\ref{theorem_hyperplane} 
(see \S\ref{definition_decomposable}).
This makes even bigger list of examples,
 in particular we obtained varieties of general type in every dimension:

\begin{theo}\label{theorem_exist_general_type_examples}
  Asume $n$ and $\rho$ are two integers such that either:
  \begin{itemize}
     \item $n \ge 3$ and $\rho \ge n+3$ or 
     \item $ 3 \le n \le 27$ and $\rho \ge 1$.
  \end{itemize}
   Then there exist a smooth irreducible projective variety $X$ of dimension $n$,
   Kodaira dimension $n$    
   and rank of Picard group $\rho$ which 
   admits Legendrian embedding into projective space $\P^{2n+1}$.
\end{theo}

With this method we are also able to give new examples from old ones without dropping the dimension.
This is described in detail in \S\ref{linear_section_decomposable}.

All the varieties arising from Theorem~\ref{theorem_hyperplane}
and our construction in \S\ref{section_construction} are embedded by a non-complete linear system.
Therefore a natural question arises: can the construction  be inverted?
So for a given Legendrian but not linearly normal embedding of some variety $\widetilde{X}$,
can we find a bigger Legendrian variety $X$, such that $\widetilde{X}$
is a projection of a hyperplane section of $X$?

Following \cite{jabu_hyperplane}
and building upon ideas of Bryant, Landsberg and Manivel we suggest a
construction that provides some (but far from perfect)
answer for this question in \S\ref{section_extend}:

\begin{theo}\label{theorem_extending_intro}
  Let $\widetilde{X} \subset \P^{2n-1}$ be an irreducible Legendrian subvariety.
  Then there exist a Legendrian subvariety 
  $X \subset \P^{2n+1}$, a hyperplane $H \subset \P^{2n+1}$
  and a irreducible component of $X \cap H$, 
  such that $\widetilde{X}$ is the image of this component 
  under the projection described in \S\ref{section_construction}.
\end{theo}

This theorem is a simplified version of Theorem~\ref{theorem_extending}.

In particular, we represent the example of Landsberg and Manivel as a
hyperplane section of a 3-fold with only isolated singularities and the
examples of Bryant as hyperplane sections of surfaces with at most
isolated singularities.

\medskip

Chapter~\ref{section_basics} is devoted to introducing our notation and presenting 
some elementary  algebro-geometric facts.

Chapter~\ref{chapter_differential} revises the 
differential geometric properties of infinitesimal automorphisms
that are necessary for Chapter~\ref{chapter_contact}, 
but can be expressed without any explicit reference to the contact structure.

Chapter~\ref{chapter_symplectic} is a brief revision of symplectic geometry 
that will be used in our discussion of contact manifolds. 
Also some statements from \cite{jabu06} are generalised to this context.

Chapter~\ref{chapter_contact} contains an independent review of local geometry 
of contact manifolds, with emphasis on their infinitesimal automorphisms. 
There we compare (after \cite{lebrun} and \cite{beauville}) 
two natural Lie algebra structures related to a contact manifold $Y$: 
the Lie bracket of vector fields and the Poisson bracket on the structure sheaf 
of the symplectisation of $Y$. 
We use this comparison to prove the first theorem on embedded automorphisms of 
Legendrian subvarieties (see Theorem~\ref{theorem_introduction_contact_automorphisms_preserving_legendrian}).

In Chapters~\ref{chapter_automorphisms}--\ref{chapter_hyperplane} we turn our attention
to Legendrian subvarieties of a projective space.

In Chapter~\ref{chapter_automorphisms} we continue the topic of automorphisms of 
Legendrian varieties. We prove the second theorem on embedded automorphisms of 
Legendrian subvarieties in a projective space (see Theorem~\ref{corollary_that_conjecture_is_true_for_smooth}).
The results of this chapter are published in \cite{jabu_toric}.

In Chapter~\ref{chapter_toric} we illustrate, in the case of subvarieties of projective space,
how to classify toric Legendrian subvarieties and  
we give the list of all smooth cases (see Theorem~\ref{theorem_classification_of_smooth_toric_legendrian})
Also the results of that chapter are published in \cite{jabu_toric}.

Chapter~\ref{chapter_sl} contains the classification of Legendrian varieties, 
which are contained in a specific $F$-cointegrable variety. 
Another new example arises in this way: the smooth quasihomogeneous 8-fold
(see Theorem~\ref{theorem_hyperplane_section_of_G_3_6}). 
Also we present two other variants of the construction, 
producing a smooth 5-fold and a smooth 14-fold.
The contents of that chapter are published in \cite{jabu_sl}.

Finally Chapter~\ref{chapter_hyperplane} describes 
a Legendrian embedding of a hyperplane section of a Legendrian variety
(see Theorem~\ref{theorem_hyperplane}).
Also a variant of an inverse construction 
(i.e.~to describe a bigger Legendrian variety from a given one, 
such that a hyperplane section of the big one is the original one)
is presented and is applied to Bryant's, Landsberg's and Manivel's 
examples of smooth Legendrian varieties (see Theorem~\ref{theorem_extending_intro}).
Parts of that chapter are published in \cite{jabu_hyperplane}.
The new parts of this chapter are \S\ref{linear_section_decomposable},
 where we construct a new series of examples 
(see Theorem~\ref{theorem_exist_general_type_examples})
and \S\ref{section_self_dual},
where we explain the relation of certain Legendrian varieties 
with smooth varieties with smooth dual.

\section{Open problems}

Keeping in mind the elegant results sketched in \S\ref{section_state_of_art}
and having many new examples of smooth Legendrian varieties 
(as well as families of such), several natural questions remain unanswered.

\subsubsection{New contact manifolds?}

Can we construct a new example of a contact manifold, 
whose variety of tangent directions to contact lines 
is one of the new Legendrian varieties 
(or is in the given family)?
If Conjecture~\ref{con_contact} is true, then the answer is negative.
If the answer is negative,
then what are the obstructions, i.e., what conditions should we require on the
Legendrian variety to make the reconstruction of contact manifold
possible?

\subsubsection{Further applications to algebra?}

Can the new Legendrian varieties be used in a similar manner 
as the subadjoint cases and will they prove themselves to be equally extraordinary varieties?
The first tiny piece of evidence for this is explained in \S\ref{section_other_examples}.
On the other hand, it is unlikely that such a big variety of examples can have analogous
special properties.

\subsubsection{Self-dual varieties?}

Another problem we want to mention here is a classical question in projective geometry: 
what are the smooth subvarieties of projective space, whose dual variety%
\footnote{%
Given a subvariety $Z\subset \P(W)$, 
\textbf{the dual variety} $Z^*\subset \P(W^*)$ 
is the closure of the set of hyperplanes tangent to $Z$, 
see \S\ref{section_definition_of_dual_variety} for details.} is also smooth?
So far the only examples of these  are the self-dual varieties.
Thanks to L.~Ein \cite{ein}, the classification of smooth self-dual varieties $Z \subset \P^m$ 
is known when $3 \codim Z \ge \dim Z$. 
In Corollary~\ref{corollary_self_dual} we prove that the problem 
of classifying smooth varieties with smooth dual can be expressed in terms of 
Legendrian varieties and possibly we can apply the techniques of Legendrian varieties 
to finish the classification.

\subsubsection{Projectively and linearly normal Legendrian varieties?}

We dare to conjecture:

\begin{con}\label{con_Legendrian_linearly_normal}
  Let $X\subset \P(V)$ be a smooth linearly normal\footnote{
A subvariety $X\subset\P^{m}$ is \textbf{linearly normal} if it is embedded by a complete linear system.}
Legendrian variety. 
Then $X$ is one of the subadjoint varieties.
\end{con}

In the view of Theorems~\ref{theorem_hyperplane} and \ref{theorem_extending_intro},
the classification of linearly normal Legendrian varieties 
might be a necessary step towards a classification of Legendrian varieties.

Furthermore, the conjecture might also contribute to the proof of Conjecture~\ref{con_contact}.
For instance assume Conjecture~\ref{con_Legendrian_linearly_normal} holds and
$Y$ is a contact Fano manifold, for which the variety cut out by contact lines 
through a general point is normal. 
Then by applying Theorem~\ref{theorem_kebekus_lines}
we get that the associated Legendrian variety $X\subset \P^{2n-1}$ 
is projectively normal\footnote{
A subvariety $X\subset\P^{m}$ is \textbf{projectively normal} 
if its affine cone is normal. 
If $X$ is projectively normal, then it is also linearly normal by 
\cite[Ex.~II.5.14(d)]{hartshorne}}  and by the conjecture and results of \cite{hong} 
the manifold $Y$ is an adjoint variety.

The author is able to prove Conjecture~\ref{con_Legendrian_linearly_normal} 
if $\dim X=1$, but this is not an elegant argument nor 
does it have important applications. 
We omit the proof here until we manage to improve the argument or to generalise it to higher dimensions.

\chapter{Notation and elementary properties}\label{section_basics}

In the present monograph we always work over the field of complex numbers $\C$.

\section{Vector spaces and projectivisation}

Let $V$ be a vector space over $\C$.
By $\P(V)$ we mean the naive projectivisation of $V$, i.e.~the
quotient  $(V \backslash \{0\})/ \C^*$.

If $v\in V \setminus \set{0}$, then 
by $[v]\in \P(V)$ we denote the line spanned by $v$.

\medskip

Analogously, if $E$ is a vector bundle, by $\P(E)$ we denote the naive 
projectivisation of $E$.
Let $s_0\subset E$ be the zero section of $E$. 
If $v \in E \setminus s_0$, then by $[v]\in \P(E)$
we denote the line spanned by $v$ 
in the appropriate fibre of $E$.

\section{Bilinear forms and their matrices}
\label{section_bilinear_forms_and_matrices}

Let $V$ be a complex vector space of dimension $m$ and $f$ a bilinear form on $V$. 
Fix a basis $\ccB$ of $V$ and let $M(f)= M(f, \ccB)$ 
be the $m\times m$-matrix such that:
$$
f(v,w) = v^T M(f) w,
$$
where $v$ and $w$ are arbitrary column vectors of $V$. 
We say that $M(f)$ is \textbf{the matrix of $f$ in the basis $\ccB$}.

In particular if $\omega$ is a symplectic form (see \S\ref{section_symplectic_vector_space}),
$\dim V =2n$ 
and $\ccB$ is a symplectic basis, then 
$$
J:=M(\omega, \ccB) = 
\left[ 
\begin{array}{cc}
0          & \Id_{n} \\
-\Id_{n} & 0
\end{array}
\right].
$$
Moreover in such a case $-J$ is also the matrix of the linear map $\tilde{\omega}$:
\begin{align*}
\tilde{\omega}: \  V & \lra  V^* \\
v & \longmapsto  \omega(v, \cdot)  
\end{align*}
in the basis $\ccB$ on $V$ and the dual basis on $V^*$.

Similarly, if $q$ is a quadratic form on $V$, then we denote by $M(q)= M(q, \ccB)$
\textbf{the matrix of $q$ in the basis $\ccB$}:
$$
q(v) = v^T M(q) v.
$$

\section{Complex and algebraic manifolds}

Our main concern is with complex projective manifolds and varieties.
This is where two categories meet: complex algebraic varieties and analytic spaces (see \cite{griffiths_adams}).
Since the author's origins lie in algebraic geometry, this monograph's intention is
to study algebraic Legendrian varieties. 
However, for some statements there is no reason to limit to the algebraic case, 
so we state them also for the analytic situation.

So $Y$ will be usually the ambient manifold (for example contact or symplectic manifold),
either a complex manifold or smooth algebraic variety. 
Some statements are local for $Y$ (in the analytic topology), 
hence it is enough to prove them for $Y\simeq D^{2n}$, where $D^{2n}\subset \C^{n}$ 
is a  complex disc.

Our main interest is in $X\subset Y$, 
which will be either an analytic subspace (if $Y$ is a complex manifold),
or an algebraic subvariety (if $Y$ is algebraic). 
For short, will always say $X\subset Y$ is a subvariety.

\section{Vector bundles, sheaves and sections}

Given an analytic space or algebraic variety $Y$, 
we denote by $\ccO_Y$ both the structure sheaf 
(consisting of either holomorphic or algebraic functions on $Y$
in the appropriate analytic or Zariski topology)
and the trivial line bundle.
If $X\subset Y$ is a subvariety, 
then by $\I(X)$ we mean the sheaf of ideals in $\ccO_Y$ defining $X$.

Given a vector bundle $E$ on $Y$ we will use the same letter $E$ 
for the sheaf of sections of $E$. 
To avoid confusion and too many brackets (for example $\I(X)(U)$) 
given an open subset $U\subset Y$ and a sheaf (or vector bundle) $\ccF$, 
we will write $H^0(U,\ccF)$ rather than $\ccF(U)$ 
to mean the value of the sheaf at the open subset $U$ (or sections of vector bundle).
By $\ccF|_U$ we mean the sheaf (or vector bundle) restriction of $\ccF$ to the open subset $U$.

Where there can be  no confusion, 
given a sheaf $\ccF$
which does not have any natural vector bundle structure 
we will write $s\in \ccF$ to mean:
\[
\exists \text{ an open } U \subset Y \ \text{ with } s \in H^0(U,\ccF).
\]
On the other hand, if $E$ is a vector bundle,
then by $v\in E$, we mean that $v$ is a vector in the bundle.

Given a vector bundle $E$, we denote by $E^*$ the dual vector bundle:
\[
E^*:= \cHom(E,\ccO).
\]

If $\theta \colon \ccF \lra \ccG$ is a map of sheaves or vector bundles
and $s \in H^0(U,\ccF)$, then by $\theta(s)$ we mean the image section of $\ccG$.

\section{Derivatives}

Given a complex manifold or smooth algebraic variety $Y$ 
and a holomorphic or algebraic $k$-form  \linebreak[3] 
\mbox{$\theta\in H^0(U,\Omega^k Y)$} by $\ud \theta$ 
we denote the exterior derivative of $\theta$. 
This convention is also valid for $0$-forms:
\[
 f\in H^0(U,\ccO_Y)= H^0(U,\Omega^0 Y).
\]
By $TY$ we mean the tangent vector bundle. 
Nevertheless we keep in mind, that a vector field $\mu \in H^0(U, TY)$ can also be interpreted as
a derivation 
\[
 \mu: \ccO_Y \to \ccO_Y. 
\]
In particular, we can define the Lie bracket of two vector fields $\mu, \nu \in H^0(U, TY)$ as:
$$
[\mu, \nu] = \nu \mu - \mu\nu.
$$
This convention is in agreement with \cite{arnold}.

Given a holomorphic or algebraic map $\phi:Y\lra Y'$, 
by $\uD \phi$ we mean the derivative map:
\[
\uD\phi: T Y \lra \phi^*TY'.
\]

If $\theta \in H^0(U,\Omega^k Y)$ and $\mu \in H^0(U, TY)$, 
then by $\theta(\mu)$ we mean the contracted $(k-1)$-form. 
For example, if $\theta = \theta_1 \wedge \theta_2$ for 1-forms $\theta_i$, 
then 
\[
\theta(\mu) = \theta_1(\mu) \theta_2 - \theta_2(\mu) \theta_1.
\]

\section{Homogeneous differential forms and vector fields}
\label{section_homogeneous_forms}
Let $Y$, $Y'$ be two complex manifolds and let $\phi: Y' \lra Y$ be a holomorphic map.
For a $k$-form $\omega \in H^0(Y,\Omega^k Y)$, 
by $\phi^* \omega \in H^0(Y',\Omega^k Y')$ we denote the pull-back of $\omega$:
$$
(\phi^*\omega)_y (v_1, \ldots, v_k) := 
\omega_{\phi(y)} \big(\uD_{y}\phi (v_1), \ldots, \uD_{y}\phi (v_k) \big).
$$ 

Now assume we have a $\C^*$-action on $Y$:
$$
(t, y) \longmapsto \lambda_t (y).
$$
We say that $\omega \in H^0(Y,\Omega^k Y)$ is \textbf{homogeneous of weight $\wt (\omega)$} if 
$$
\forall t \in \C^* \quad \lambda_t^* \omega = t^{\wt(\omega)} \omega.
$$

For example, assume $Y = \A^n = \Spec(\C[y_1,\ldots, y_n])$ and $\C^*$ acts via homotheties. 
We say $\omega \in \Omega^k {\A^n}$ is \textbf{constant}, 
if it is a $\C$-linear combination of $\ud y_{i_1} \wedge\ldots  \wedge \ud y_{i_k}$.
Constant $k$ forms are homogeneous of weight $k$ (not of weight 0 as one could possibly expect).
Conversely, if $\omega \in H^0(\A^n, \Omega^k \A^n)$ is homogeneous of weight $k$, 
then it is constant, because every global form can be written as 
$\sum f_{i_1,\ldots, i_k} \ud y_{i_1}\wedge\ldots \wedge \ud y_{i_k}$. 
Since $\ud y_{i_1}\wedge \ldots \wedge \ud y_{i_k}$ are already of weight $k$, it follows that 
$f_{i_1,\ldots, i_k}$ are constant functions.

Let $\mu \in H^0(Y, TY)$ be a vector field. 
We say $\mu$ is \textbf{homogeneous of weight $\wt(\mu)$} if
$$
\uD \lambda_{t^{-1}} \mu = t^{\wt(\mu)} \mu.
$$

\renewcommand{\theenumi}{(\roman{enumi})}
\begin{lem}\label{lemma_on_homogeneous_forms_and_vector_fields}
  Let $Y$, $Y'$ be complex manifolds, both with a $\C^*$-action.
  Moreover assume  $\phi: Y' \lra Y$ is a $\C^*$-equivariant map, 
  $\omega \in H^0(Y,\Omega^k Y)$ is a homogeneous $k$-form
  for some $k \in \set{0, 1, \ldots, \dim Y}$
  and $\mu \in H^0(Y, TY)$, $\nu \in  H^0(Y', TY')$
  are two homogeneous vector fields.
  \begin{enumerate}
    \item 
        \label{item_omega(mu)_is_homogeneous}
      $\omega(\mu)$ is homogeneous and 
      $\wt \big(\omega(\mu)\big)  =  \wt(\omega) + \wt(\mu)$;
    \item
        \label{item_phi_upper_star_preserves_weights}
      $\phi^*\omega$ is homogeneous of weight $\wt(\omega)$ and 
      $\uD \phi (\nu)$  is homogeneous of weight $\wt(\nu)$;
    \item\label{item_derivative_preserves_weights}
      $\ud \omega$ is homogeneous of weight $\wt(\omega)$.
  \end{enumerate}
\end{lem}

\begin{prf}
This is an immediate calculation. For instance \ref{item_omega(mu)_is_homogeneous}:
\begin{multline*}
  \lambda^*_t (\omega(\mu))_x (v_1,\ldots, v_{k-1}) 
  = 
  \omega_{\lambda_t(x)}(\mu, \uD\lambda_t (v_1), \ldots, \uD\lambda_t (v_{k-1}) )  
  =
\\
  (\lambda^*_t\omega)_x (\uD\lambda_{t^{-1}} (\mu) , v_1, \ldots,  v_{k-1} ) 
  =  
  t^{\wt(\omega)} t^{\wt(\mu)}(\omega(\mu))_x (v_1,\ldots,v_{k-1}).
\end{multline*}

\end{prf}

\section{Submersion onto image}

We recall the standard fact, that every algebraic map is generically 
a submersion on the closure of the image.
\begin{lem} \label{lemma_algebraic_map_submersion}
  Let $Y$ and $Y'$ be two algebraic varieties over an algebraically closed field of 
characteristic 0 and let $\pi: Y \lra Y'$ be a map such that 
$Y'=\overline{\pi(Y)}$.
 Then for a  general $y\in Y$, 
the derivative $D_y \pi: T_y Y\lra T_{\pi(y)} Y' $ is surjective.
\end{lem}

\begin{prf}
  See \cite[Cor.~III.10.7]{hartshorne}.
\end{prf}

As a corollary, we prove an easy proposition about subvarieties of product manifolds.

\begin{prop}\label{proposition_decomposable}
  Let $Y_1$ and $Y_2$ be two smooth algebraic varieties and suppose $X \subset Y_1 \times Y_2$ 
  is a closed irreducible subvariety.
  Let $X_i \subset Y_i$ be the closure of the image of $X$ under the projection $\pi_i$ onto $Y_i$. 
  Assume that for a Zariski open dense subset of smooth points $U\subset X$
  the tangent bundle to $X$ decomposes as
  \[
    T X|_U = (T X \cap \pi_1^* T {Y_1})|_U \oplus (T X \cap \pi_2^* T {Y_2})|_U,
  \]
  a sum of two vector bundles.
  Then $X= X_1 \times X_2$.
\end{prop}

\begin{prf}
  Since $X$ is irreducible, so is $X_1$ and $X_2$ and clearly $X \subset X_1 \times X_2$.
  So it is enough to prove that  
  \[
    \dim X_1 + \dim X_2 = \dim X = \dim U.
  \]
  However, the maps $\uD({\pi_i}|_U)$ are surjective onto $T X \cap {\pi_i}^*T {Y_i}$ and hence applying 
  Lemma~\ref{lemma_algebraic_map_submersion}:
  \[
    \dim X_1 + \dim X_2 = \rk (T X \cap {\pi_1}^* T Y_1)|_U + \rk (T X \cap {\pi_2}^* T {Y_2})|_U = \rk T X|_U = \dim X.
  \]
\end{prf}

\section{Tangent cone}\label{properties_of_tangent_cone}

We recall the notion of the tangent cone and a few among many of its
properties. For more details and the proofs we refer to 
\cite[Lecture~20]{harris} and \cite[III.\S3,\S4]{mumford}.

For an irreducible Noetherian scheme $X$ over $\C$ and a closed point $x\in X$ 
we consider the local ring $\ccO_{X,x}$ and we 
let $\gotm_x$ be the maximal ideal in $\ccO_{X,x}$.
Let 
$$
R:=\bigoplus_{i=0}^{\infty} \left(\gotm_x^i \slash \gotm_x^{i+1} \right),
$$
where $\gotm_x^0$ is just the whole of $\ccO_{X,x}$. 
Now we define \emph{the tangent cone $TC_x X$ at $x$ to $X$} to be $\Spec R$.

If $X$ is a subscheme of an affine space $\A^m$ (which we will usually
assume to be an affine piece of a projective space),
the tangent cone at $x$ to $X$ can be understood as a subscheme of $\A^m$. 
Its equations can be derived from the ideal of $X$.
For simplicity assume $x=0 \in \A^m$ and then the polynomials defining $TC_0 X$ 
are the lowest degree homogeneous parts of the polynomials in the ideal of $X$.

Another interesting point-wise definition is that $v \in TC_0 X$ is a closed point if and only if 
there exists a holomorphic map $\varphi_v$ from
the disc $D_t:=\{ t\in \C  \ : \ \lvert t \rvert\ < \delta \}$ to $X$,
such that $\varphi_v(0) = 0$
 and the first non-zero coefficient in the Taylor expansion in $t$ of $\varphi_v(t) $
is $v$, i.e.:
$$ \begin{array}{rccl}
\varphi_v: & D_t & \lra & X\\
           &  t  & \mapsto & t^k v + t^{k+1}v_{k+1} +\ldots
\end {array}$$

\smallskip

We list some of the properties of the tangent cone, which will be used
freely in the proofs:

\begin{itemize}
\item[(1)]
The dimension of every component of $TC_x X$ is equal to the dimension of $X$. 
\item[(2)] $TC_x X$ is naturally embedded in the Zariski tangent space to $X$ at $x$ and 
$TC_x X$ spans (as a scheme) the tangent space.
\item[(3)] 
$X$ is regular at $x$ if and only if $TC_x X$ is equal (as a scheme) to
	  the tangent space.
%\item[(4)]
%If $TC_x X$ is reduced, then $X$ is reduced at $x$.
\end{itemize}

%\begin{prf}
%The property (4) follows  from the definition and Krull Intersection Theorem 
%(see for example \cite[Cor.~5.4]{eisenbud}).
%\end{prf}

\section{Secant, tangent and dual varieties}\label{section_definition_of_dual_variety}

Let $W$ be a vector space of dimension $n+1$.
Let $Z \subset \P^n=\P(W)$ be any subvariety.

We denote by $\sigma(Z) \subset \P^n$ its 
\textbf{secant variety}, 
i.e.,  closure of the union of all projective lines through $z_1$ and $z_2$,
where $(z_1,z_2)$ vary through all pairs of different points of $Z$.

By $\tau(Z)$ we denote \textbf{the tangent variety of $Z$},
i.e.~the closure of the union of limits of secant lines as $z_1$ converges to $z_2$.

Also:
\begin{defin}
  We let $Z^* \subset \check{\P}^n:=\P(W^*)$ 
  be the closure of the set of hyperplanes tangent to $Z$ at some point:
  \[
    Z^*:= \overline{\set{H \in \check{\P}^n \mid \exists z \in Z \text{ s.t.~} T_z Z \subset H}}.
  \]
  We say $Z^*$ is \textbf{the dual variety} to $Z$.
\end{defin}

\chapter{Vector fields, forms and automorphisms}
\label{chapter_differential}

In the course of  this dissertation,
particularly in Chapter~\ref{chapter_contact} 
we use some differential geometric facts, which we summarise in this chapter.
Although all these facts are standard or follow  easily from the standard material,
we reproduce or at least sketch most of the proofs.
We do this for the sake of completeness of the material presented in the monograph 
and also because various authors of textbooks use various notations 
and combining them one can get very confused 
(at least this has happened to the author of this monograph).

\section{Vector fields, Lie bracket and distributions}

Let $Y$ be a complex manifold or a smooth algebraic variety, 
 let $F\subset TY$ be a corank $1$ subbundle%
\footnote{
One could also consider $F$ to be a corank $r$ subbundle for any $r \in \set{1\ldots \dim Y}$.
Some of the statements below can be generalised to any $r$ (not necessary $r=1$),
but the proofs get more complicated, especially in notation. 
We restrict our considerations to the $r=1$ case, 
as this is the only one used in the dissertation.
}
and let $\theta: T Y \to TY/F =:L$ be the quotient map, so that the following sequence is exact:
$$
0\lra F \lra T Y \stackrel{\theta}{\lra} L \lra 0.
$$
Also assume $U$ is an open subset.
We say that a (possibly singular) subvariety $X\subset U$ with its smooth locus $X_0$
is \textbf{$F$-integrable} if $T X_0$ is contained in $F$.

\renewcommand{\theenumi}{(\roman{enumi})}
\begin{prop}\label{properties_of_distribution}
  With the assumptions as above:
  \begin{enumerate}
    \item 
      \label{item_dtheta_well_defined}
      $\ud \theta$ gives a well defined map of $\ccO_Y$-modules: 
      $$
        \ud \theta: \Wedge{2} F \lra L.
      $$
      We refer to this map as \textbf{the twisted 2-form $\ud \theta$}.
    \item 
      \label{item_bracket_depends_only_on_vectors}
      Assume $\mu$ and $\nu$ are two vector fields on $U$, both contained in $F$.
      Then $\theta([\mu ,\nu])(y) =  \ud \theta_y(\mu(y), \nu(y))$.
      In particular $\theta([\mu ,\nu])(y)$ does not depend on the vector fields, 
      but only on their values at $y$.
    \item 
      \label{item_bracket_depends_only_on_vectors2}
      Again assume  $\mu$ and $\nu$ are two vector fields on $U$,
      but now only $\nu$ is contained in $F$. Then again $\theta([\mu ,\nu])(y)$
      depends only on the value of $\nu$ at $y$, but not on the whole vector field.
      In other words the map of sheaves $F \lra L$ 
      given by $\theta\big([\mu , \cdot]\big)$ is $\ccO_Y$-linear 
      and hence it determines a map of vector bundles $F \lra L$.
    \item
      \label{item_Y_tangent_to_F}
      If $X$ is $F$-integrable, then $\ud \theta|_{X_0} \equiv 0$. 
      In particular: 
      $$
        \dim X \le \rk F - \half \min_{x \in X} \left( \rk \ud\theta_x \right).
      $$
  \end{enumerate}
\end{prop}

\begin{prf}
All the statements are analytically local, so it is enough to assume that $Y$ 
is a disc $D^{2n} \subset \C^n$ with coordinates $ y_1, \ldots, y_m$, $U=Y$, $y=0$ 
and that $\theta$ is a nowhere 
vanishing section of $\Omega^1 Y\otimes L \simeq \Omega^1 Y$ 
(the choice of the trivialisation of $L$ is of course not unique):
$$
\theta = \sum_{i} A_i \ud y_i = \boldsymbol{A} \cdot \ud \boldsymbol{y},
$$
where the collection $\left( A_1, \ldots, A_m \right)$ 
(respectively $\left( \ud y_1, \ldots, \ud y_m \right)^T$)
we denote by $\boldsymbol{A}$ 
(respectively $\ud \boldsymbol{y}$). 
Then:
$$
F:= \set{ v \in T D^{2n} \mid \sum_{i} A_i \ud y_i (v) = 0 }.
$$

To prove \ref{item_dtheta_well_defined} note that: 
$$
\ud \theta = 
\sum_{i} \ud A_i \wedge  \ud y_i 
= \ud \boldsymbol{A} \wedge \ud \boldsymbol{y}.
$$ 
We must check that this does not depend on the choice of 
the trivialisation $\boldsymbol{A}$ of $L$.
So assume $\boldsymbol{B}$ is a  different trivialisation, so
there exists $g: Y \lra \GL(1)\simeq \C^*$ such that:
$$
\boldsymbol{B}= g \cdot \boldsymbol{A}. 
$$
We must prove that $\ud \boldsymbol{B}\wedge \ud \boldsymbol{y}$ 
restricted to $F$ transforms in the same manner:
$$
\ud \boldsymbol{B} \wedge \ud \boldsymbol{y}  
= \ud(g \cdot \boldsymbol{A}) \wedge \ud \boldsymbol{y} =
(\ud g \cdot \boldsymbol{A} + g \cdot  \ud \boldsymbol{A}) \wedge \ud \boldsymbol{y} = 
$$
$$ 
\stackrel{\textrm{since $\boldsymbol{A}$ vanish on $F$}}{=}
(g \cdot  \ud \boldsymbol{A}) \wedge \ud \boldsymbol{y}. 
$$

\smallskip

To prove \ref{item_bracket_depends_only_on_vectors} let
\begin{align*}
\mu = &\sum_k \mu_k \frac{\pd}{\pd y_k},\\
\nu = &\sum_k \nu_k \frac{\pd}{\pd y_k}
\end{align*}
for some holomorphic functions $\mu_k$ and $\nu_k$.
Since $\mu$ and $\nu$ are contained in $F$  we have:
$$
\sum_k A_k \mu_k = 0 \quad \textrm{ and } \quad \sum_l A_l \nu_l = 0.
$$
Therefore for every $k$ or $l$ we have:
\begin{subequations}
    \label{equation_derivative_for_lie_bracket}
  \begin{align}
    \sum_k \frac{\pd A_k}{\pd y_l} \mu_k & = - \sum_k A_k \frac{\pd \mu_k}{\pd y_l};\\
    \sum_l \frac{\pd A_l}{\pd y_k} \nu_l & = - \sum_l A_l \frac{\pd \nu_l}{\pd y_k}.
  \end{align}
\end{subequations}

Since
\[
[\mu, \nu] = \sum_{k,l} 
\left(
\nu_k \frac{\pd \mu_l}{\pd y_k} \frac{\pd}{\pd y_l} - 
\mu_l \frac{\pd \nu_k}{\pd y_l} \frac{\pd}{\pd y_k} 
\right),
\]
hence:
\[
\begin{split}
\theta([\mu,\nu])  &
= \sum_{k,l} \left(
A_l \nu_k \frac{\pd \mu_l}{\pd y_k}  - 
A_k \mu_l \frac{\pd \nu_k}{\pd y_l} 
\right)=
\\
&\stackrel{
\textrm{by \eqref{equation_derivative_for_lie_bracket}}
}{=}
 \sum_{k,l} \left(
 - \frac{\pd A_l}{\pd y_k} \mu_l \nu_k + 
   \frac{\pd A_k}{\pd y_l}   \mu_l \nu_k
\right) = 
\\
&=
 \sum_{k,l} \left(
 \frac{\pd A_l}{\pd y_k}\left(\mu_k \nu_l- \mu_l \nu_k \right)
\right) = 
\\
&= \sum_{k,l} \left(
 \frac{\pd A_l}{\pd y_k}\left(\ud y_k \wedge \ud y_l\right) \left(\mu , \nu\right)
\right) = 
\\
&= \ud \theta \left(\mu, \nu\right).
\end{split}
\]
We note that the above calculation is a special case of \cite[Prop.~I.3.11]{kobayashi_nomizu},
though the reader should be careful, 
as the notation in \cite{kobayashi_nomizu} 
is different than ours and as a consequence a constant factor $-2$ 
is ``missing'' in our formula. 

\smallskip

The proof of \ref{item_bracket_depends_only_on_vectors2} 
is identical as the beginning of the proof of \ref{item_bracket_depends_only_on_vectors}.

\smallskip

Finally to prove \ref{item_Y_tangent_to_F} just use \ref{item_bracket_depends_only_on_vectors} 
and the fact that the Lie bracket of two vector fields tangent to $X$ must be tangent to $X$.
\end{prf}

%\begin{lem}
%Let $\omega_1$, $\omega_2$ be two holomorphic $k$-forms on a disc $D^{2n} \subset \C^n$ with coordinates
%$x_1, \ldots, x_n$. 
%Assume there exists $k\in $
%\end{lem}

\section{Automorphisms}\label{section_inf_automorphisms}

Here we introduce the notation about several types of automorphisms of a manifold $Y$ 
and its subvariety $X$. Also we recall some standard properties and relations between them.

\medskip

Let $Y$ be a complex manifold (or respectively, smooth algebraic variety) 
and let $U\subset Y$ be an open subset in analytic (or respectively, Zariski) topology. 
By $\Authol(U)$ (respectively, $\Autalg(U)$) 
we denote the group of holomorphic (respectively, algebraic) automorphisms of $U$.
By $\Autdot(U)$ we mean either $\Authol(U)$ or $\Autalg(U)$,
whenever specifying is not necessary.

Assume that a complex Lie  group (respectively, an algebraic group) $G$ 
acts on $U$, i.e.~we have a group homomorphism $G \lra \Autdot(U)$.
Also let $\gotg$ be the Lie algebra of $G$. 
By $G^0$ we denote the the connected component of identity in $G$.

\textbf{An infinitesimal automorphism of $U$} 
is a vector field 
$\mu\in H^0(U, TY)$. 
Differentiating the action map $G\times U \lra U$ by the first coordinate we get 
the induced map $\gotg\times Y \lra TY$ or more precisely $\gotg \lra H^0(U, TY)$.
This map 
preserves the Lie bracket 
(see \cite[Thm in \S1.7]{akhiezer})  
and if the action is faithful, then it is injective 
(see \cite[Thm in \S1.5]{akhiezer}). 

The particular case is when $G =\C^*$. 
Then we get a map $\C\lra H^0(U,TY)$ and we set $\mu_{\C^*}$ 
to be the image of $1\in \C$ under this map.
 We say $\mu_{\C^*}$ is \textbf{the vector field related to the $\C^*$-action.}
Note that $\mu_{\C^*}$ is homogeneous of weight $0$.

The infinitesimal automorphisms make a sheaf $TY$ of Lie algebras, 
which at the same time is an $\ccO_Y$-module. 
The two multiplication structures are related by the following Leibniz rule:
\begin{equation}\label{equation_leibniz_for_lie_bracket_of_vector_fields}
\forall f \in H^0(U, \ccO_Y), \
\forall \mu, \nu \in H^0(U, TY) \quad
[f\mu,\nu] = f[\mu,\nu] + \ud f(\nu) \mu. 
\end{equation}

The following theorem comparing infinitesimal, algebraic and holomorphic automorphisms 
for a projective variety is well known and standard:

\begin{theo}\label{theorem_automorphisms_of_projective_variety}
Let Y be a projective variety. Then:
\begin{enumerate}
\item \label{item_Authol_is_Lie}
 $\Authol(Y)$ is a complex Lie group.
\item \label{item_automorphisms_hol_alg}
Every holomorphic automorphism of $Y$ is algebraic and hence
$$
\Aut(Y):=\Authol(Y) = \Autalg(Y).
$$ 
\item \label{item_automorphisms_infinitesimal}
By $\gotaut(Y)$ we denote the tangent Lie algebra to $\Aut(Y)$.
Every infinitesimal automorphism is tangent to some 1-parameter subgroup of $\Authol(Y)$, 
so that $\gotaut(Y) = H^0(Y,TY)$.
\end{enumerate}
\end{theo}

\begin{prf}
Part~\ref{item_Authol_is_Lie} is proved in \cite[\S2.3]{akhiezer}.
Part~\ref{item_automorphisms_hol_alg} is a consequence of \cite[Thm~IV.A]{griffiths_adams}.
Part~\ref{item_automorphisms_infinitesimal} is explained in 
\cite[Prop. in \S1.5 \& Cor.~1 in \S1.8]{akhiezer}.
\end{prf}

%From now on, we will only use the superscripts ${}^{hol}$ and ${}^{alg}$,
%when there is a significant difference in the statements. 

Clearly $H^0(U,\ccO_Y)$ is a representation of $G$ 
and hence also of $\gotg$. 
We also have 
the following Lie algebra action
of the sheaf of infinitesimal automorphisms:
$$
\begin{array}{rcl}
TY \times \ccO_Y & \lra &\ccO_Y \\
(\mu, f)  & \mapsto &  \ud f (\mu),
\end{array}
$$
which is given by the derivation in the direction of the vector field.

The action of $\gotg$ on $ H^0(U,\ccO_Y)$ is the composition
$$
\gotg \lra H^0(U, TY) \lra \gotgl\left(H^0(U, \ccO_Y)\right).
$$

\medskip

Let $X\subset Y$ be a subvariety. 
By $\Autdot(U,X)$ we denote the respective subgroup of $\Autdot(U)$ preserving the intersection $U\cap X$. 
If $Y$ is projective, then by $\gotaut(Y,X)$ we mean the Lie algebra tangent to $\Autdot(Y,X)$.
By $\autinf(U,X)$ we denote the Lie algebra of infinitesimal automorphisms 
of $U$ preserving $X$, i.e.:
$$
\autinf(U,X) := \left\{\mu \in H^0(U, TY) \mid
\forall f \in \I(X)|_U \quad (\ud f)(\mu) \in \I(X)|_U \right\},
$$
where $\I(X) \ideal \ccO_Y$ is the sheaf of ideals of $X$.

Clearly, if $G$ preserves $X$,
then the image of $\gotg \lra H^0(U, TY)$ 
is contained in $\autinf(U,X)$. 
Conversely, if  the image is contained in $\autinf(U,X)$,
then the action of the connected component $G^0$ preserves $X$.

\begin{cor}
  If $Y$ is projective, then $\autinf(Y,X) = \gotaut(Y,X)$.
\end{cor}

\noprf

Moreover $\autinf(\cdot, X)$ makes in $TY$ a subsheaf of Lie algebras 
and $\ccO_Y$-modules.

\section{Line bundles and $\C^*$-bundles}\label{notation_Ldot}

Let $Y$ be complex manifold or a smooth algebraic variety 
and let $L$ be a line bundle on $Y$. 
By $\Ldot$ we denote the principal $\C^*$-bundle over $Y$ 
obtained as the line bundle $L^*$ with the zero section removed. 
Let $\pi$ be the projection $\Ldot \lra Y$. 

Let $\ccR_L$ be the sheaf of graded $\ccO_Y$-algebras $\bigoplus_{m \in \Z} L^{m}$ on $Y$.
Given an open subset $U \subset Y$ the ring $\ccR_L(U)$ 
consists of all the  algebraic functions on $\pi^{-1}(U)$,
i.e. $\ccR_L = \pi_*\ccO_{\Ldot}$. Therefore 
$$
\Ldot = \SheafySpec_Y \ccR_L.
$$
Moreover, $H^0(U,L^m) \subset H^0(\pi^{-1}(U), \ccO_{\Ldot})$ is 
the set of homogeneous functions of  weight $m$ 
(see \S\ref{section_homogeneous_forms}).

\label{section_on_Pic_of_Cstar_bundles}

\begin{lem}\label{lemma_on_Pic_of_Cstar_bundles}
  Let $Y$ be a smooth algebraic variety and let $L$ be a line bundle on $Y$.
  Then $\Pic(\Ldot) \simeq \Pic(Y) / \langle L \rangle$ 
  and the map  $\Pic (Y)  \epi \Pic (\Ldot)$ 
  is  induced by the projection $\pi : \Ldot \lra Y$.
\end{lem}

\begin{prf}
  The Picard group of the total space of $L^*$ is isomorphic to $\Pic Y$ 
and the isomorphisms are given by the projection and the zero section $s_0: Y \lra L^*$.
Further, $s_0(Y)$ is a Cartier divisor linearly equivalent to any other rational section 
$s: Y \dashrightarrow L^*$. 
Therefore $s_0^*(s_0(Y)) = L^*$ and hence by \cite[Prop.~6.5(c)]{hartshorne}  
the following sequence is exact:
$$
\begin{array}{cccccccc}
\Z &\lra &\Pic Y& \stackrel{\pi^*}{\lra}& \Pic &\Ldot & \lra & 0\\
1 & \longmapsto & [L^*]
\end{array}
$$
\end{prf}

The relative tangent bundle, i.e.~$\ker \left(\uD \pi : T \Ldot \lra \pi^*TY\right)$ 
is trivialised by the vector field 
$\mu_{\C^*}$ related to the action of $\C^*$ (see \S\ref{section_inf_automorphisms})
and hence we have the short exact sequence:
$$
0 \lra \ccO_{\Ldot} \lra T\Ldot \lra \pi^*TY \lra 0.
$$
In particular $K_{\Ldot} = \pi^* K_Y$.

\section{Distributions and automorphisms preserving them}
 \label{section_automorphisms_of_distributions}

Let $F\subset TY$ be a corank 1 vector subbundle. 
Our particular interest will be in the case where $F$ is a contact distribution - 
see \S\ref{section_contact} for the definition.
But what follows holds true for any $F$.
By 
$\Autdot_{F}(U)$, $\gotaut_{F}(Y)$, $\autinfF(U)$,
$\Autdot_{F}(U,X)$, $\gotaut_{F}(Y,X)$ and $\autinfF(U,X)$ 
we denote the appropriate automorphisms or infinitesimal automorphisms preserving $F$
and possibly the subvariety $X$.

For instance, 
\begin{equation}
  \label{equ_definition_of_autinfF}
  \autinfF(U) = \left\{ \mu \in H^0(U,TY) \mid [\mu,F] \subset F \right\}.
\end{equation}

Also $\autinfF$ makes a sheaf of Lie algebras, but usually
it is not an $\ccO_Y$-sub\-mo\-dule of $TY$. 
To see that take any $\mu\in \autinfF(U)$ for $U$ small enough. 
Assume for all $f \in \ccO_Y(U)$ we have $f \mu \in \autinfF(U)$.
Then by Leibniz rule (see Equation~\ref{equation_leibniz_for_lie_bracket_of_vector_fields}):
\[
  \forall \nu \in H^0(U,F) \quad \ud f (\nu) \cdot \mu \in H^0(U,F).
\]
This can only happen if either:
\begin{itemize}
\item $\mu \in H^0(U,F)$ or
\item $F=0$, i.e.~$F$ is the rank 0 bundle.
\end{itemize}
We will see that the first case does not happen if $F$ is a contact distribution 
(unless $\mu=0$, 
see Theorem~\ref{theorem_contact_automorphisms}\ref{item_contact_automorphisms}).
In fact, one can prove that it never happens for all $\mu \in H^0(U,F)$ 
(remember that $U$ is small enough), unless $F=0$.

If $G$ acts on $U$ and preserves the distribution $F$,
then the map $\gotg \to H^0(U,TY)$ factors through $\autinfF(U)$.
Conversely, if $G$ is connected, it acts on $U$ 
and the map $\gotg \to H^0(U,TY)$ factors through $\autinfF(U)$,
then the action of $G$ preserves $F$. 
As a consequence we get:

\begin{cor}\label{corollary_Y_projective_infinitesimal_and_tangent_to_global}
  If $Y$ is projective and $X\subset Y$  is a subvariety,
  then:
  \begin{enumerate}
    \item $\gotaut_F(Y) = \autinfF(Y)$
    \item $\gotaut_F(Y,X) = \autinfF(Y,X)$ 
  \end{enumerate}
\end{cor}

\begin{prf}
This follows from the above considerations and from 
Theorem~\ref{theorem_automorphisms_of_projective_variety}.

\end{prf}

Further, let $L$ be the quotient bundle and $\theta$ be the quotient map:
$$
0 \lra F \lra TY \stackrel{\theta}{\lra} L \lra 0.
$$
If the action of $G$ on $U$ extended to $TY|_U$ preserves $F$, 
then in the obvious way we get 
\textbf{the induced action of $G$ on the total spaces of $L|_U$ and $L^*|_U$}.
These actions preserve the zero sections.

\medskip

Let $\Ldot$ and  $\ccR_L$ be as in \S\ref{notation_Ldot}. 

By analogy with above we want to define the action of $\autinfF$ on $\Ldot$. 
In other words, we define a special lifting of the vector fields from 
$\autinfF \subset TY$ to vector fields on $\Ldot$.

First observe that the sheaf of Lie algebras $\autinfF$ acts on the sheaf $L$: 
if $s \in H^0(U,L)$,  
then choose an open subset $V\subset U$ small enough 
and any lifting $s_{TY}\in H^0(V,TY)$, $\theta(s_{TY})=s|_V$ 
and let $\mu \in \autinfF(U)$ act on $H^0(U,L)$ locally by 
\begin{equation}
  \label{equ_definition_of_action_of_autinfF_on_L}
  s|_V \mapsto (\mu.s)|_V := \theta\left(\left[s_{TY}, \mu|_V \right]\right).
\end{equation}
By Equation~\eqref{equ_definition_of_autinfF} defining $\autinfF$, 
this does not depend on the choice of $s_{TY}$ and 
hence, by elementary properties of sheaves, 
it glues uniquely to an element of $H^0(U,L)$.
Hence we get a Lie algebra representation $\autinfF(U)\lra \gotgl\left(H^0(U,L)\right)$.

Secondly, we can extend the action of $\autinfF$ on the locally free sheaf $L$
defined in Equation~\eqref{equ_definition_of_action_of_autinfF_on_L}
to an action on $\ccR_L$, by requesting that the action must satisfy the Leibniz rule:
\begin{equation}\label{equ_leibnitz_for_autinfF}
t, s \in \ccR_L,  \  \mu \in \autinfF  \ \Longrightarrow \ \mu.(ts) = (\mu.t)s + t(\mu.s)
\end{equation}
--- locally every section of $L^m$ can be written as
a sum of products of sections of $L$  (or their inverses, if $m < 0$).

Finally, we can extend this action to $\ccO_{\Ldot}$, again requesting the Leibniz rule.
Eventually, we get the action, 
which we will call \textbf{the induced action of $\autinfF$ on $\Ldot$}.
The following property justifies the name:

\begin{prop}
If the action of $G$ preserves $F$, then the tangent action 
to the induced action of $G$ on $\Ldot|_U:= \pi^{-1}(U)$ is the composition of 
$\gotg \lra \autinfF(U)$  and the induced action of $\autinfF$ on $\Ldot$.
\end{prop}
\noprf

For a fixed $\mu \in \autinfF(U)$, the induced map 
$\ccO_{\Ldot}|_{\pi^{-1}(U)} \lra \ccO_{\Ldot}|_{\pi^{-1}(U)}$
is a derivation, so it corresponds to a vector field $\breve{\mu} \in H^0(\pi^{-1}(U), T\Ldot)$,
such that
\begin{equation}
  \label{equation_define_mu_tilde}
  \forall f \in \ccO_{\Ldot} \quad \mu . f = \ud f (\breve{\mu}).
\end{equation}
By construction we also have $\uD \pi (\breve{\mu}) = \mu$.

\section{Lifting and descending twisted forms}
  \label{section_thetadot}

In this section we explain how to lift the twisted form $\theta:T Y \to L$ 
to a form $\thetadot$ on $\Ldot$. Later we study the derivative $\ud \thetadot$ 
and show its relation to $\ud \theta : \Wedge{2} F \to L$.
Finally we prove the homogeneity of $\thetadot$ and $\ud \thetadot$ 
and thus explicitly set conditions when a $\C^*$-bundle together with 
a $2$-form $\omega$ arises as $(\Ldot, \ud \thetadot)$ from some distribution
$F \subset TY$.

With the notation and assumptions as in the previous sections,
we have a canonical isomorphism of line bundles 
$
\tau: \pi^* L \stackrel{\simeq}{\lra} \ccO_{\Ldot}
$:
if $y\in Y$, $\lambda \in \Ldot_y=\pi^{-1}(y)$, $l \in L_y$, then we set
$$
\tau(y,\lambda, l):= (y,\lambda, \lambda(l)).
$$ 

We let $\thetadot:= \tau \circ \pi^*\theta \circ \uD \pi$:%
\footnote{%
In \cite{beauvillefano}, \cite{lebrun} the authors denote $\thetadot$ 
simply as $\pi^*\theta$, since the other maps are natural.
 This is a bit confusing to some people (including the author of this monograph, 
but see also a comment in \cite{sola_wisniewski} about a small mistake in 
\cite{4authors}) and therefore we underline that $\thetadot$ is the composition of three maps.
}
\[
 T \Ldot 
  \stackrel{\uD \pi}{\lra} 
 \pi^* TY  
  \stackrel{\pi^* \theta}{\lra} 
 \pi^* L  
  \stackrel{\tau}{\lra} 
 \ccO_{\Ldot}.
\]

\begin{lem}\label{lemma_Lie_derivative_of_thetadot_is_0}
  For every $\mu \in \autinfF(U)$ 
  the induced infinitesimal automorphism $\breve{\mu}$
  preserves $\thetadot$, i.e.:
  \[
    L_{\breve{\mu}} (\thetadot) 
    := \lim_{t\to 0} \frac{\gamma_{\breve{\mu}}(t)^* \thetadot - \thetadot}{t}
    = 0,
  \]
  where $L_{\breve{\mu}}$ is the Lie derivative operator 
  and $\gamma_{\breve{\mu}}(t)$ 
  is the local 1-parameter group of transformations of $\Ldot$ determined by $\breve{\mu}$.
\end{lem}

\begin{prf}
  For the simplicity of notation assume $\gamma_{\breve{\mu}}(t)$ is a global transformation.
  The following diagram of vector bundles is commutative:
  \[ 
   \xymatrix{  
      T {\Ldot} \ar@{->}[rr]^{\uD \pi}      \ar@{->}[d]^{\uD\gamma_{\breve{\mu}}(t)} && 
      \pi^* TY  \ar@{->}[rr]^{\pi^* \theta} \ar@{->}[d]^{\uD^{\pi}\gamma_{\mu}(t)}         &&
      \pi^* L   \ar@{->}[rr]^{\tau}         \ar@{->}[d]^{\gamma^{\pi^*L}_{\mu}(t)}         &&
      \ccO_{\Ldot}                         \ar@{->}[d]^{\gamma_{\breve{\mu}}(t) \times \id_{\C} }\\
      T {\Ldot} \ar@{->}[rr]^{\uD \pi} && 
      \pi^* TY  \ar@{->}[rr]^{\pi^* \theta} &&
      \pi^* L   \ar@{->}[rr]^{\tau} &&
      \ccO_{\Ldot}  
   }
   \]
   where by $\uD^{\pi}\gamma_{\mu}(t)$ we mean the automorphism of $\pi^* TY$, 
   which is determined by $\uD\gamma_{\mu}(t): TY \to TY$ and 
   $\gamma_{\breve{\mu}}(t):\Ldot\to \Ldot$; 
   similarly $\gamma^{\pi^*L}_{\mu}(t)$ is determined by 
   $\uD\gamma_{\mu}(t): TY/F \to TY/F$ and  $\gamma_{\breve{\mu}}(t):\Ldot\to \Ldot$.
   The composition of the whole upper row is equal to $\thetadot$. 
   The composition of the left most vertical arrow and the whole lower row is equal to 
   $\gamma_{\breve{\mu}}(t)^* \thetadot$. 
   Since the right arrow is the identity on the second component 
   $\ccO_{\Ldot} = \Ldot \times \C$ and since the diagram is commutative,
   both forms take the same values on every vector $v \in T {\Ldot}$, hence are equal and
   the claim follows.   
\end{prf}

We also give a local description of $\thetadot$ and $\ud \thetadot$.
So now assume $Y \simeq D^{2m}$ and let $y_1, \ldots, y_m$ be some coordinates on $Y$.
Let $z$ be a linear coordinate on the fibre of $\Ldot \simeq Y \times \C^*$.
This means that $z$ determines a section of $L$ which trivialises $L$ over $D^{2m}$.
So we can think of $\theta$ as of a holomorphic $1$-form on $\Ldot$ 
depending only on $y_i$'s and $\ud y_i$'s.
Let $(y,z_0)$ be any point of $\Ldot$ 
and let $\bar{v}$ be any vector tangent to $\Ldot$ at $(y,z_0)$.
We write $\bar{v} = v +w$, where $v$ is the component tangent to $Y$,
while $w$ is tangent to $\C^*$.
Then:
$$
  \thetadot_{(y,z_0)}(\bar{v})  = (\tau \circ \pi^*\theta \circ \uD \pi)_{y,z_0}(\bar{v}) = 
   (\tau \circ \pi^*\theta)_{y,z_0}(v) = z_0 (\theta_y(v)) = z_0 \cdot \theta_y(v)
$$ 
or more concisely (in local coordinates)
\begin{equation}\label{equ_locally_thetadot}
  \thetadot = z \theta,
\end{equation}
and therefore
\begin{equation}\label{equ_locally_dthetadot}
\ud \thetadot = \ud \left( z \theta \right) = 
z\ud \theta + \ud z \wedge \theta.
\end{equation}

Since in this notation $\theta$ is a homogeneous 1-form of weight $0$ and $\wt(z) =1$,
$\thetadot$ and $\ud \thetadot$ are homogeneous forms of weight 1
(see \S\ref{section_homogeneous_forms}).

In the above coordinates, 
the vector field $\mu_{\C^*}$ related to the $\C^*$-action can be expressed as follows:
$$
\mu_{\C^*} = z \frac{\pd}{\pd z}.
$$

\renewcommand{\theenumi}{(\roman{enumi})}
\begin{prop}
    \label{prop_on_symplectisation}
  Let $Y$ be a complex manifold or smooth algebraic variety and let $L$ be a line bundle on $Y$.
  Also let $\Ldot$ be the principal $\C^*$-bundle over $Y$ as in \S\ref{notation_Ldot} and
  let $\mu_{\C^*}$ be the vector field on $\Ldot$ associated to the action of $\C^*$.
  Finally, let $\omega$ be a homogeneous closed 2-form on $\Ldot$ of weight $1$. 
  Then:
  \begin{enumerate}
    \item \label{item_omega_and_omega(mu)}
       $\omega = \ud \big(\omega(\mu_{\C^*})\big)$; 
    \item \label{item_omega_comes_from_theta}
       There exists a unique twisted 1-form $\theta: TY \lra L$,
       such that $\omega(\mu_{\C^*})=\thetadot$,
       where $\thetadot$ is defined from $\theta$ as above;
    \item \label{item_omega_and_theta_non_vanishing}
       Moreover, $\omega(\mu_{\C^*})$ is nowhere vanishing 
       if and only if $\theta$ is nowhere vanishing.
       If this is the case, 
       then $\omega$ is non-degenerate if and only if $\ud \theta|_F$ is non-degenerate.  
  \end{enumerate}
\end{prop}

\begin{prf}
  To prove \ref{item_omega_and_omega(mu)} 
  let $z$ be a local coordinate linear on the fibres of $\pi: \Ldot \to Y$.
  Since $\omega$ is closed, locally it is exact, so 
  \[
  \omega = \ud (z \phi' + g \ud z)
  \]
  for some function $g$ and 1-form $\phi'$, both homogeneous of weight 0. However,
  \[
    \ud (z \phi' + g \ud z) = \ud \big(z (\phi' - \ud g)\big).
  \]
  Set $\phi:= \phi' - \ud g$, so that $\omega = \ud (z \phi)$.
  Note that although $\phi'$ and $g$ are not uniquely determined, 
  $\phi$ is the unique homogeneous 1-form of weight 0 
  such that $\omega = \ud (z \phi)$.
  Then,
\[ 
\omega \left(\mu_{\C^*}\right) 
= (\ud z \wedge \phi)
\left(z \frac{\pd}{\pd z}\right)
+ z \ud \phi  \left(z \frac{\pd}{\pd z}\right)
=
\ud z \left(z \frac{\pd}{\pd z}\right) \cdot \phi = z \phi.  
\] 
Hence $\ud\bigl(\omega \left(\mu_{\C^*}\right)\bigr)= \omega$, as claimed in \ref{item_omega_and_omega(mu)}.

\smallskip
 
To prove \ref{item_omega_comes_from_theta}, 
define $\theta$ to be locally the form $\phi$ from the above argument. 
One must verify that $\phi$ glues uniquely to a twisted 1-form $\theta: TY \lra L$.

\smallskip

Part  \ref{item_omega_and_theta_non_vanishing} follows from the local descriptions 
of $\thetadot$ and $\ud \thetadot$, 
see Equations~\eqref{equ_locally_thetadot} and \eqref{equ_locally_dthetadot}.
For instance, if $n=\half (\dim Y-1) $, then:
$$
  (\ud \thetadot)^{\wedge^{n+1}} = (n+1) \ud z \wedge \theta \wedge (\ud\theta)^{\wedge^n}.
$$
Therefore $\ud \thetadot$ is non-degenerate at a given point 
if and only if $\theta$ does not vanish at that point  
and $\ud \theta$ is non-degenerate on the kernel of $\theta$.
\end{prf}

\begin{lem}
  \label{lemma_integrable_iff_dthetadot_vanish}
Let $X\subset Y$ be any subvariety and $X_0$ its smooth locus. 
Then $X$ is $F$-integrable if and only if $\ud \thetadot$ 
vanishes identically on the tangent space to $\pi^{-1}(X_0)$.
\end{lem}

\begin{prf}
First assume $X$ is $F$-integrable. 
Then $\ud \theta$ vanishes on $T \left(\pi^{-1}(X_0)\right)$ 
by Proposition~\ref{properties_of_distribution}\ref{item_Y_tangent_to_F}
and $\theta$ vanishes by definition. 
Hence from the local description of $\ud \thetadot$ 
(see Equation~\eqref{equ_locally_dthetadot}) we get the result.

\smallskip

On the other hand if $\ud \thetadot|_{T \left(\pi^{-1}(X_0)\right)} \equiv 0$, since
$$
\mu_{\C^*}|_{\pi^{-1}(X_0)} \in H^0\left(\pi^{-1}(X_0), T \left(\pi^{-1}(X_0)\right)\right),
$$
then in particular  
$$
\ud \thetadot \left(\mu_{\C^*}, T\left(\pi^{-1}(X_0)\right)\right) \equiv 0.
$$
But $\ud \thetadot \left(\mu_{\C^*}\right)= \thetadot$ 
(see Proposition~\ref{prop_on_symplectisation}\ref{item_omega_comes_from_theta}),
hence $\pi^{-1}X$ is $(\pi^*F)$-integrable and therefore $X$ is $F$-integrable.
\end{prf}

\medskip

For $s \in \ccR_L = \pi_* \ccO_{\Ldot}$,
by $\tilde{s} \in \ccO_{\Ldot}$
we denote the lifting of $s$, i.e.~$\tilde{s}:=\tau\circ\pi^*s$.

Hence we have two different possibilities of lifting an infinitesimal automorphism 
\mbox{$\mu \in \autinfF$} to an object on $\Ldot$: 
either we lift it to vector field $\breve{\mu}$ (see Equation~\eqref{equation_define_mu_tilde})
or  we lift $\theta(\mu)$ to function $\widetilde{\theta(\mu)}$.
We will compare these two liftings 
and how they behave with respect to the Lie bracket of vector fields in 
the following statements.

\begin{lem}
  \label{lemma_on_mu_nu_theta}
  We have:
  $$
    \forall \nu \in \autinfF(U), \mu \in H^0(U,TY) 
    \qquad \widetilde{\theta([\mu,\nu])} = \ud \left(\widetilde{\theta(\mu)}\right)(\breve{\nu}).
  $$
\end{lem}

\begin{prf}
  By Equation~\eqref{equ_definition_of_action_of_autinfF_on_L}:
  $$
    \theta([\mu,\nu]) = \nu . {\theta(\mu)}
  $$
  and hence $ \widetilde{\theta([\mu,\nu])} = \nu . \widetilde{\theta(\mu)} $.
  By Equation~\eqref{equation_define_mu_tilde}, 
  this is equal to $\ud \left(\widetilde{\theta(\mu)}\right)(\breve{\nu})$.
\end{prf}

\begin{prop}\label{proposition_on_mu_and_thetadot}
  If $\mu \in \autinfF(U)$, then:
  $$
    \ud \left(\widetilde{\theta(\mu)}\right) = - (\ud \thetadot) (\breve{\mu}).
  $$
\end{prop}

\begin{prf}
  The following proof is quoted from \cite[Prop.~1.6]{beauvillefano}.
  Since $L_{\breve{\mu}}(\thetadot) =0$ (see Lemma~\ref{lemma_Lie_derivative_of_thetadot_is_0}),
  by \cite[Prop.~I.3.10(a)]{kobayashi_nomizu}
  we have:
  $$
    (\ud \thetadot)(\breve{\mu}) = - \ud \bigl(\thetadot (\breve{\mu})\bigr).
  $$
  On the other hand:
  $$
  \thetadot (\breve{\mu}) = \tau \circ \pi^*\theta \circ \uD \pi (\breve{\mu})
  = \tau \circ \pi^*\left(\theta(\mu)\right) = \widetilde{\theta(\mu)}.
  $$
  Combining the two equalities, we get the result. 
\end{prf}

\begin{cor}
    \label{corollary_bracket_of_mu_nu}
  If $\mu, \nu \in \autinfF (U)$, then 
  $$
    \widetilde{\theta([\mu, \nu])} =  - (\ud \thetadot)(\breve{\mu}, \breve{\nu}).
  $$
\end{cor}
\begin{prf}
  This combines Lemma~\ref{lemma_on_mu_nu_theta} 
  and Proposition~\ref{proposition_on_mu_and_thetadot}.
\end{prf}

\chapter{Elementary symplectic geometry}
\label{chapter_symplectic}

We introduce some elementary facts from symplectic geometry, 
having in mind the needs of subsequent chapters.
Most of this material is contained in (or can be easily deduced from)
classical textbooks on symplectic geometry, such as \cite{mcduff_salamon},
although we rewrite this over the ground field $\C$ rather than $\R$.

\section{Linear symplectic geometry}

In this section we study linear algebra of vector space, 
which has a symplectic form.
Although it is elementary, 
it is very important for our considerations as it has threefold application:
Firstly, the content of this section describes the local behaviour of symplectic manifolds 
(see \S\ref{section_symplectic_manifolds}), 
particularly the symplectisations of contact manifolds (see \S\ref{section_sandwich}).
Secondly, it describes very much of global geometry of projective space as a contact manifold
(see Example~\ref{example_projective_space_as_contact}, but also look through 
Chapters~\ref{chapter_automorphisms}--\ref{chapter_hyperplane}).
Finally, it explains the fibrewise behaviour of contact distribution (see \S\ref{section_contact}).

\subsection{Symplectic vector space}
\label{section_symplectic_vector_space}

A symplectic form on a vector space $V$ is 
a non-degenerate skew-symmetric bilinear form.
So $\omega \in \Wedge{2} V^*$ is a symplectic form 
if and only if 
$$
\forall v \in V  \ \exists w \in V 
\ \textrm{ such that } \ \omega(v, w) \ne 0
$$
or equivalently the map 
$$
\begin{array}{rcl}
\tilde{\omega}: \  V & \lra & V^* \\
v & \longmapsto & \omega(v, \cdot)
\end{array}
$$
is an isomorphism.

If a vector space $V$ has a symplectic form $\omega$,
we say that $V$ (or $(V,\omega)$ if specifying the form is important)
is \textbf{a symplectic vector space}.
In such a case the dimension of $V$ is even  and there exists 
a basis $v_1,\ldots, v_n, w_1,\ldots, w_n$ (where $n=\half \dim V$)
 of $V$ such that 
$\omega(v_i, w_i)=1$,   
$\omega(v_i, v_j)=0$ and
$\omega(v_i, w_j)=0$ for $i \ne j$.
Such a basis is called \textbf{a symplectic basis}.

By $\omega^{\vee}$ we  
denote the corresponding symplectic form on $V^*$:
$$
\omega^{\vee}:= \left(\tilde{\omega}^{-1} \right)^* \omega.
$$
Note that if $v_1,\ldots, v_n$, $w_1,\ldots, w_n$ is a symplectic basis of $V$
and $x_1,\ldots, x_n$, $y_1,\ldots, y_n$ is the dual basis of $V^*$, 
then $x_1,\ldots, x_n, y_1,\ldots, y_n$ is a symplectic basis of $V^*$. 
In such a case $x_1,\ldots, x_n, y_1,\ldots, y_n$ 
are also called \textbf{symplectic coordinates on $V$}.

\subsection{Subspaces in a symplectic vector space}

Assume $V$ is a vector space of dimension $2n$ and $\omega$ is a symplectic form on $V$.
Now suppose $W\subset V$ is a linear subspace. 
By $W^{\perp_{\omega}}$ we denote the $\omega$ perpendicular complement of $W$:
$$
W^{\perp_{\omega}}:=\left\{v \in V \mid \forall w\in W \quad \omega(v,w) = 0  \right\}.
$$
Denote by $\pi$ the natural projection $V^*\ra W^*$.
We say that the subspace $W$ is:

\bigskip

\begin{centering}
\begin{tabular}{|p{0.155\textwidth}|cp{0.225\textwidth}ccc|p{0.14\textwidth}|}

\hline
\textbf{isotropic}& 
$\Leftrightarrow$ & $\omega|_W \equiv 0$&
$\Leftrightarrow$ & $W \subset W^{\perp_\omega}$&
$\Leftrightarrow$ & $\ker \pi$ is coisotropic;\\
\hline
\textbf{coisotropic} (or sometimes called \textbf{involutive})&
$\Leftrightarrow$ & $\omega^{\vee}|_{\ker \pi} \equiv 0$&
$\Leftrightarrow$ & $W \supset W^{\perp_\omega}$&
$\Leftrightarrow$ & $\ker \pi$ is isotropic;\\
\hline
\textbf{Lagrangian}&
$\Leftrightarrow$ & $W$ is isotropic or involutive and \mbox{$\dim W=n=\half \dim V$}&
$\Leftrightarrow$ & $W = W^{\perp_\omega}$&
$\Leftrightarrow$ & $\ker \pi$ is Lagrangian;\\
\hline
\textbf{symplectic}&
$\Leftrightarrow$ & $\omega|_W$ is a symplectic form on $W$ &
$\Leftrightarrow$ & $W \cap W^{\perp_\omega} = 0$&
$\Leftrightarrow$ & $\ker \pi$ is symplectic.\\
\hline
\end{tabular}
\end{centering}

\subsection{Symplectic reduction of vector space}\label{notation_symplectic_reduction}

With the assumptions as above let $W\subset V$ 
be any linear subspace and let $W':=W \cap W^{\perp_{\omega}}$.
Define $\omega'$ to be the following bilinear form on $V':=W\slash W'$:
$$
\textrm{for $\ w_1,w_2 \in W \ $ let } \  \omega'([w_1], [w_2]):= \omega(w_1,w_2).
$$ 
Then $(V',\omega')$ is a symplectic vector space.

The particular case we are mostly interested in is when $W$ is a hyperplane 
or more generally a coisotropic subspace.

Note the following elementary properties of this construction:
\begin{prop}
\label{proposition_symplectic_reduction}
For a subspace  $L\subset V$ let $L'$ be the image of 
$L\cap W$ in $V'$.
\renewcommand{\theenumi}{(\alph{enumi})}
\begin{enumerate}
\item 
If $L$ is isotropic (resp.~coisotropic or Lagrangian) in $V$,
then $L'$ is isotropic (resp.~coisotropic or Lagrangian) in $V'$.
\item \label{proposition_symplectic_reduction_Lagrangian}
Conversely, if $W$ is coisotropic, $L \subset W$ and 
$L'$ is isotropic (resp.~coisotropic or Lagrangian) in $V'$,
then $L$ is isotropic (resp.~coisotropic or Lagrangian) in $V$.
\end{enumerate}
\end{prop}
\noprf

\subsection{Symplectic automorphisms and weks-symplectic matrices}\label{section_symplectic_automorphisms}

A linear automorphism $\psi$ of a symplectic vector space $(V,\omega)$ is called
a \textbf{symplectomorphism}
if $\psi^* \omega =  \omega$ i.e.:
\[
\forall u,v \in V \quad \omega\bigl(\psi(u), \psi(v) \bigr) = \omega\bigl(u, v \bigr). 
\]
We denote by $\Sp (V)$ the group of all symplectomorphisms of $V$ and
by $\gotsp(V)$ its Lie algebra:
\[
\gotsp(V)= 
\set{  
g \in \End(V) \mid \forall u,v \in V \quad 
\omega\bigl(u, g(v) \bigr) + \omega\bigl(g(u), v \bigr) = 0
}.
\]

A linear automorphism $\psi$ of $V$ is called
a \textbf{conformal symplectomorphism}
if $\psi^* \omega = c \omega$ for some constant $c \in \C^*$. 
We denote by $\cSp (V)$ the group of all conformal symplectomorphisms of $V$ and
by $\gotcsp(V)$ the tangent Lie algebra.

\medskip

Fix a basis $\ccB$ of $V$ and note
 that a matrix $g\in \gotgl(V)$ is in the symplectic algebra $\gotsp(V)$ if and only if 
$$g^T J + Jg =0$$
where $J:=M(\omega, \ccB)$. 
For the sake of Chapter~\ref{chapter_automorphisms}
we also need to define a complementary linear subspace to $\gotsp(V)$:
\begin{defin}
A matrix $g\in \gotgl(V)$ is \textbf{weks-symplectic}\footnote{
A better name would be \emph{skew-symplectic} or \emph{anti-symplectic}, but these are reserved for some different notions.
}
if and only if it satisfies the equation:
$$
g^T J - Jg =0.
$$
The vector space of all weks-symplectic matrices will be denoted by $\gotasp(V)$
(even though it is not a Lie subalgebra of $\gotgl(V)$).
\end{defin}

We immediately see that a matrix is weks-symplectic if and only if it corresponds to 
a linear endomorphism $g$, 
such that for every $u,v\in V$:
\begin{equation}\label{weks_symplectic_equation}
\omega(g u, v) - \omega (u, g v) =0.
\end{equation}
This is a coordinate free way to describe $\gotasp(V)$.

\medskip

Assume that our basis $\ccB$ is symplectic.
In particular $J^2=M(\omega,\ccB)^2  = -\Id_{2n}$.

\begin{rem}\label{weks_corresponds_to_skew}
For a matrix $g\in \gotgl(V)$ we have:
\begin{itemize}
\item[(a)] $g \in \gotsp(V) \  \iff  \  Jg$ is a symmetric matrix;
\item[(b)] $g \in \gotasp(V) \  \iff  \  Jg$ is a skew-symmetric matrix.
\end{itemize}
\end{rem}
\noprf

Note that if $g\in \gotgl(V)$, then we can write:
$$
g = \half( g + J g^T J) + \half (g - J g^T J)
$$
and the first component $g_+:= \half( g + J g^T J)$ is in $\gotsp(V)$, 
while the second $g_- := \half (g - J g^T J)$ is in $\gotasp(V)$. 
Obviously, this decomposition corresponds to expressing the matrix $Jg$ 
as a sum of symmetric and skew-symmetric matrices.

We list some properties of $\gotasp(V)$:
\begin{prop}\label{properties_of_asp}
Let $g, h \in \gotasp(V)$. The following properties are satisfied: 
\begin{itemize}
\item[(i)] Write the additive Jordan decomposition for $g$:
$$
g = g_s + g_n 
$$
where $g_s$ is semisimple and $g_n$ is nilpotent.
Then both $g_s \in \gotasp(V)$ and $g_n \in \gotasp(V)$.

\item[(ii)] For $\lambda \in \C$, 
denote by $V_{\lambda}$ the \mbox{$\lambda$-eigenspace} of $g$.
For any $\lambda_1, \lambda_2 \in \C$ two different eigenvalues  
$V_{\lambda_1}$ is $\omega$-perpen\-di\-cu\-lar to $V_{\lambda_2}$.

\item[(iii)] If $g$ is semisimple, then each space $V_{\lambda}$ is symplectic.

\end{itemize}
\end{prop}
\noprf

\subsection{Standard symplectic structure on $W \oplus W^*$}
  \label{section_symplectic_on_W_plus_W_dual}
Let $W$ be any finite dimensional vector space. Set $V:=W \oplus W^*$ 
and there is a canonical symplectic form on $V$:
$$
\omega\big((v,\alpha), (w, \beta)\big):= \beta(v) - \alpha(w).
$$
If $a_1, \ldots, a_n$ is any basis of $W$ and 
$\lambda_1, \ldots, \lambda_n$ is the dual basis of $W^*$,
then 
$$
  a_1, \ldots, a_n,\lambda_1, \ldots, \lambda_n
$$
is a symplectic basis of $V$.
In particular, we have the natural embedding 
\begin{align*}
\GL(W) &\mono \Sp(V)\\
A &\mapsto A\oplus (A^{-1})^T.
\end{align*}

We note the following elementary lemma:

\begin{lem}\label{lemma_L+ker_is_lagrangian}
  Let $L\subset W$ be any linear subspace. 
  Then $L \oplus \ker(W^* \to L^*) \subset V$ 
  is a Lagrangian subspace.
\end{lem}
\noprf

\section{Symplectic manifolds and their subvarieties}
\label{section_symplectic_manifolds}

Symplectic manifolds will serve us to understand some geometric and algebraic structures 
of the symplectisations of contact manifolds (see \S\ref{section_sandwich}).

\medskip

A complex manifold or a smooth complex algebraic variety $Y$ is \textbf{a symplectic manifold} 
if there exists a global closed holomorphic 2-form $\omega \in H^0(\Omega^2 Y)$ with
\mbox{$\ud \omega =0$} which restricted to every fibre is a symplectic form 
on the tangent space. 
In other words, 
$\omega^{\wedge^n}$ is a nowhere vanishing section of $K_Y=\Omega^{2n} Y$.
The form $\omega$ is called \textbf{a symplectic form on $Y$}.

Similarly as in the case of the vector space, 
the symplectic form determines an isomorphism:
$$
\begin{array}{rcl}
\tilde{\omega}: \  TY & \stackrel{\simeq}{\lra} & T^* Y  \\
v & \longmapsto & \omega(v, \cdot).
\end{array}
$$

The theory of compact (or projective) complex symplectic manifolds is well developed 
and has a lot of beautiful results 
(see for example \cite{lehn}, \cite{huybrechts} and references therein).
Yet here we will only use some non-compact examples 
as a tool for studying contact manifolds 
and we will only need a few of their basic properties.
Also some extensions of the symplectic structure to the singularities of $Y$ are studied, 
but we are interested only in the case where $Y$ is smooth.

\subsection{Subvarieties of a symplectic manifold}
\label{section_lagrangian_subvarieties}

Let $(Y,\omega)$ be a symplectic manifold.
For a subvariety $X \subset Y$ 
we say $X$ is respectively
\begin{enumerate}
  \item 
    \textbf{isotropic},
  \item 
    \textbf{coisotropic},
  \item 
    \textbf{Lagrangian},
\end{enumerate}
if and only if there exists an open dense subset $U$ 
(equivalently, for any open dense subset $U$)
 of smooth points of $X$, 
such that for every $x \in U$ the tangent space $T_x X \subset T_x Y$ is respectively
\begin{enumerate}
  \item 
    isotropic,
  \item 
    coisotropic,
  \item 
    Lagrangian.
\end{enumerate}
Or equivalently, for every  $x \in U$ the conormal space 
$N^*_x X \subset T^*_x Y$ is respectively
\begin{enumerate}
  \item 
    coisotropic,
  \item 
    isotropic,
  \item 
    Lagrangian.
\end{enumerate}

Note that a subvariety is Lagrangian 
if and only if it is isotropic (or coisotropic) and the dimension is equal to $n$.

\subsection{Examples}
\label{section_examples_of_symplectic}

The following examples are important for our considerations, as they will appear 
as symplectisations of projective contact manifolds (see \S\ref{section_sandwich}).

\subsubsection{The affine space}

Our key example is the simplest possible:
 an affine space of even dimension.
So assume $(V,\omega)$ is a symplectic vector space of dimension $2n$.
Then take the affine space $\A^{2n}$ of the same dimension, 
whose tangent space at every point is $V$ and globally
$T \A^{2n} = \A^{2n} \times V$. 
Then $\omega$ trivially extends to the product and it is a symplectic form 
on $\A^{2n}$.

By an abuse of notation, we will denote the affine space by $V$ as well 
(so in particular a $0$ is fixed in the affine space 
and the action of $\C^*$ by homotheties is chosen).
In this setup, the form $\omega$ is homogeneous of weight 2
(see \S\ref{section_homogeneous_forms}).

\subsubsection{Products}\label{example_product_Lagrangian}

Assume $Y_1$ and $Y_2$ are two symplectic manifolds 
with symplectic forms $\omega_1$ and $\omega_2$ respectively. 
Clearly $Y_1 \times Y_2$ is a symplectic manifold with 
the symplectic form $p_1^*\omega_1 +p_2^*\omega_2$,
where the $p_i$'s are the appropriate projections.

Next, let $X_i \subset Y_i$ be two subvarieties. 
Both the $X_i$'s  are respectively 
\begin{enumerate}
  \item 
    isotropic,
  \item 
    coisotropic,
  \item 
    Lagrangian,
\end{enumerate}
 if and only if the product $X_1 \times X_2 \subset Y_1 \times Y_2$ is respectively
\begin{enumerate}
  \item 
    isotropic
  \item 
    coisotropic,
  \item 
    Lagrangian.
\end{enumerate}

\subsubsection{Cotangent Bundle}
Let $M$ be a complex manifold or a smooth algebraic variety of dimension $n$. 
Set $Y$ to be the total space of the cotangent vector bundle $T^* M$
 and let $p:Y\lra M$ be the projection map.
If $x_1, \ldots, x_n$ are local coordinates on $U \subset M$,
then   $x_1, \ldots, x_n, y_1=\ud x_1, \ldots, y_n = \ud x_n$ form 
the local coordinates on $Y|_U$. Then we can set:
$$
\omega|_U:= \ud x_1 \wedge \ud y_1 +\ldots +\ud x_n \wedge \ud y_n \in H^0(U,\Omega^2 Y),
$$
and these glue to a well defined symplectic form $\omega \in H^0(Y,\Omega^2 Y )$.
This symplectic form is homogeneous of weight $1$ 
with respect to the usual $\C^*$-action on the cotangent spaces.

Since for $m \in M$, $x \in T^*_m M$ we have $T_{(m,x)}Y =T_m M \oplus T^*_m M$
this example of symplectic manifold, generalises the standard symplectic structure on 
$W \oplus W^*$ (see \S\ref{section_symplectic_on_W_plus_W_dual}).
The following example generalises Lemma~\ref{lemma_L+ker_is_lagrangian}:

\begin{ex}\label{example_conormal_Lagrangian}
Let $Z\subset M$ be any subvariety. 
Define $\hat{Z}^{\#} \subset Y$ to be 
\textbf{the conormal variety to $Z$},
i.e.~the closure of the union of conormal spaces to smooth points of $Z$:
$$
\hat{Z}^{\#}:= \overline{N^* Z_0/M}.
$$
Then $\hat{Z}^{\#}$ is a Lagrangian subvariety in $Y$.
\end{ex}

\begin{prf}
    Let $z \in Z$ be a smooth point and let $x \in N^*_z Z_0/M$. 
  Then one can choose appropriate local coordinates on $M$ around $z$ 
  and an appropriate local trivialisation of the cotangent bundle $T^* M$,
  such that:
  \[
    T_x \hat{Z}^{\#} = T_z Z \oplus N^*_z Z_0/M \subset T_z M \oplus T^*_z M. 
  \]
  This is a Lagrangian subspace by Lemma~\ref{lemma_L+ker_is_lagrangian}.
\end{prf}

\begin{lem}\label{lemma_Lagrangian_in_cotangent}
  Conversely, assume $M$ is a smooth algebraic variety 
  and $Y$ is the total space of $T^* M$.
  Moreover assume $X \subset Y$ is an irreducible closed Lagrangian subvariety
  invariant under the $\C^*$-action on $Y$. If $Z = \overline{p(X)}$, 
  then $X = \hat{Z}^{\#}$.
\end{lem}

\begin{prf}
  Let $x\in X$ be a general point  and let $z:=p(x)$.
  So $x$ is a point in $T^*_z M$ and 
  \[ 
    T_x Y = T_z M \oplus T^*_z M.
  \]
  Since $X$ is $\C^*$-invariant, under the above identification 
  \[ 
    (0,x) \in  T_x X \subset T_x Y.
  \]
  We want to prove that $(0,x) \in N^*_z Z/M$ and this will follow 
  if we prove 
  \[
    T_x X \cap T^*_z M  =  N^*_z Z/M.
  \]

  By Lemma~\ref{lemma_algebraic_map_submersion} 
  the map $\uD p: T_x X\lra T_z Z$ is surjective, 
  so 
  \[
    T_x X + T^*_z M =  T_z Z \oplus T^*_z M.
  \]   
  Since $X$ is Lagrangian, we also have the dual equality:
  \begin{align*}
    T_x X \cap T^*_z M   
    &  = \left(T_x X\right)^{\perp_{\omega}} \cap  \left(T^*_z M\right)^{\perp_{\omega}}  \\
    &  = \left(T_x X + T^*_z M \right)^{\perp_{\omega}}                                  \\
    &  = \left(T_z Z \oplus T^*_z M  \right)^{\perp_{\omega}}                            \\
    &  =  N^*_z Z/M.
  \end{align*}
  Hence $T_x X \cap T^*_z M = N^*_z Z/M$ as claimed and therefore $x \in N^*_z Z/M$.
  Since $x$ was a general point of $X$ and both $X$ and $Z$ were irreducible, 
  we have $X \subset \hat{Z}^{\#}$ and 
  by dimension counting $X = \hat{Z}^{\#}$.
\end{prf}

\subsubsection{Adjoint and coadjoint orbits}

Let $G$ be a semisimple complex Lie group and consider the coadjoint action 
on the dual of its Lie algebra $\gotg^*$.
Let $Y$ be an orbit of this action. 
The tangent space at $\xi \in Y$ is naturally isomorphic to $\gotg / Z(\xi)$,
where 
\[
  Z(\xi)= \set{v\in \gotg \mid  \forall w\in \gotg  \ \ \xi([v,w])=0 }.
\] 
Here $[v,w]$ denotes the Lie algebra operation in $\gotg$.
For $v, w \in \gotg$ let $[v]$ and $[w]$  be the corresponding vector fields on $Y$ determined by 
$v$ and $w$. We define:
$$
\omega_{\xi}([v],[w]):= \xi([v,w]).
$$
Then $\omega$ is a symplectic form on $Y$, which is called the Kostant-Kirillov form, see for instance 
\cite[(2.1)]{beauvillefano}.

Now assume $G$ is simple and 
$Y$ is invariant under homotheties 
(for instance $Y$ is the unique minimal nonzero orbit 
--- see \cite[Prop.~2.2 and Prop.~2.6]{beauvillefano}).
Then the actions of $G$ and $\C^*$ commute 
(because $G$ acts on $\gotg^*$ by linear automorphisms, 
$\C^*$ via homotheties and every linear map commutes with a homothety). 
Therefore the vector fields of the form $[v]$ for some $v\in \gotg$ are homogeneous of weight 0 
and hence:
$$
(\lambda_t^*\omega)_{\xi}([v],[w])= 
\omega_{\lambda_t(\xi)}([v],[w]) =
t \xi ([v,w])= 
t \omega_{\xi}([v],[w]).
$$
i.e.~$\omega$ is homogeneous of weight $1$.

We can identify $\gotg^*$ and $\gotg$ by Killing form (see \cite{humphreys}), 
so equally well we can consider adjoint orbits. 
Therefore if $Y$ is as above, then it is isomorphic to a \mbox{$\C^*$-bundle}
over an adjoint variety (see \S\ref{intro_contact_and_qK}). 
More precisely $Y$ is a symplectisation (see \S\ref{section_sandwich}) of the adjoint variety.

\subsubsection{Open subsets}

Let $(Y,\omega)$ be a symplectic manifold and let $U$ be an open subset.
Then $(U, \omega|_U)$ is again a symplectic manifold.

\section{The Poisson bracket}\label{section_Poisson_bracket}

The Poisson bracket is an important algebraic structure of a symplectic manifold.
In Corollary~\ref{corollary_Poisson_for_contact} 
we observe that given a contact manifold and its symplectisation,
the Poisson bracket descents from the symplectisation 
to a bracket on a specific sheaf of rings on the contact manifold.
Moreover, this descended structure is strictly related to the automorphisms of 
the contact manifold (see Theorem~\ref{theorem_contact_automorphisms}).

\medskip

Let $(Y, \omega)$ be a symplectic manifold and let $\ccO_Y$ be the sheaf of holomorphic 
(or algebraic) functions on $Y$. 
Given $f, g \in H^0(U,\ccO_Y)$ let $\xi_g \in H^0(U,TY)$ 
be the unique vector field such that $\omega(\xi_g) = \ud g$. Then we set:
\begin{align*}
   \left\{ f , g \right\}:= & \ud f(\xi_g),\\
      \intertext{or equivalently:}
   \left\{ f , g \right\}(x):= & \omega^{\vee}_x (\ud g_x , \ud f_x).
\end{align*}

The bilinear skew-symmetric map 
$\left\{ \cdot , \cdot \right\}: \ccO_Y \times \ccO_Y \ra \ccO_Y$ is called 
\textbf{the Poisson bracket}.

\begin{lem}\label{lemma_poisson_jacobi_and_leibnitz}
The Poisson bracket satisfies the Jacobi identity and therefore makes $\ccO_Y$
into a sheaf of Lie algebras.
The compatibility between the Poisson bracket and the standard ring multiplication on $\ccO_Y(U)$ 
is given by the following Leibniz rule:
$$
\left\{ f g , h \right\}=  f \left\{g , h \right\} + g \left\{ f , h \right\}.
$$
\end{lem}
\begin{prf}
See for example \cite[\S40]{arnold} --- the proof is identical to the real case.
\end{prf}

The Poisson bracket is determined by the symplectic form and 
moreover it is defined locally. 
Hence we have the following property:
\begin{prop}\label{proposition_Poisson_compatibile}
Assume $(Y,\omega)$ and $(Y', \omega')$ are two symplectic manifolds of dimension $2n$.
 Assume moreover, that  we have a finite covering map:
$$
\psi: Y \lra Y'
$$
such that $\psi^*\omega' = \omega$. 
Then the Poisson structures are compatible: for $f,g \in \ccO_{Y'}$ we have:
$$
\psi^*\{f, g \} = \{\psi^*f, \psi^*g\}.
$$
\end{prop}
\noprf

\renewcommand{\theenumi}{(\roman{enumi})}

\begin{theo}\label{Poisson_structure}
Assume $Y$ is a symplectic manifold.
\begin{enumerate}
\item Suppose $X \subset Y$ is a coisotropic subvariety.
Then the sheaf of ideals $\I(X) \subset \ccO_Y$ is a subalgebra with respect to the Poisson bracket.
\item Conversely, suppose $X \subset Y$ is a closed, generically reduced subscheme 
and that $\I(X)$ is preserved by the Poisson bracket.
Then the corresponding variety $X_{\red}$ is coisotropic.
\end{enumerate}
\end{theo}

Versions of the theorem can be found in \cite[Prop.~11.2.4]{coutinho} and in \cite[Thm 4.2]{jabu06}.
We follow more or less  the proof from \cite{jabu06}:

\begin{prf}
Let $X_0$ be the locus of smooth  points of $X$.
We must show that $\omega^{\vee}|_{N^* {X_0\slash Y}} \equiv 0$ if and only if $\I(X)$
 is a Lie subalgebra sheaf  in $\ccO_Y$.

Suppose that $x \in X_0$ is any point, $U \subset Y$ is an open neighbourhood of $x$
 and that $f,g \in H^0\left(U,\I(X)\right)$ 
are some functions vanishing on $X$. 
Then 
\[
     \ud f_x, \ud g_x \in N^* {X_0 \slash Y}.
\]

If $\omega^{\vee}|_{N^* {X_0\slash Y}} \equiv 0$, then
\[
\left\{f,g\right\}(x) = \omega^{\vee}_x \left(\ud g_x, \ud f_x \right) = 0,
\]
i.e. $\left\{f, g \right\}|_{X_0} = 0$,
 so extending the equality to the closure of  $X_0$ we get 
 $$
\left\{f,g\right\} \in H^0\left(U,\I(X)\right).
$$
Hence $\I(X)$ is a Lie subalgebra.

\smallskip

Conversely, if $\I(X)$ is a Lie subalgebra, then
$$
\omega^{\vee}(\ud g_x, \ud f_x )= \left\{f, g \right\}(x) = 0.
$$
Since the map
$$
\begin{array}{ccl}
  H^0\left(U,\I(X)\right) & \lra & N^*_{x} X_0\slash Y \\
  f & \longmapsto & \ud f_x
\end{array}
$$
is an epimorphism of vector spaces for each $x \in X_0$ and for $U$ sufficiently small, 
we have $\omega|_{N^* {X_0\slash Y}} \equiv 0$.
\end{prf}

\subsection{Properties of Poisson bracket}

In our considerations on contact manifolds and their various subvarieties
we will need the three lemmas that are explained in this subsection. 
These lemmas refer to Proposition~\ref{Poisson_structure} 
--- we have seen that there is a relation between coisotropic varieties and 
Lie subalgebras of $\ccO_Y$ that are ideals under the standard ring multiplication.

The first lemma claims that to test if an ideal is a subalgebra 
it is enough to test it on an appropriate open cover of $Y$.

\begin{lem}\label{lemma_poisson_is_affine}
  Let $Y$ be a symplectic manifold and let $\I\ideal \ccO_Y$ 
be a coherent sheaf of ideals. 
In such a case $\I$ is preserved by the Poisson bracket
if and only if there exists an open cover $\{U_i\}$ of $Y$  
such that for each $i$:
\begin{itemize}
\item 
if $V\subset U_i$ is another open subset,
then the functions in $H^0(V,\ccO_Y)$ are determined by the functions in $H^0(U_i,\ccO_Y)$
--- this means that if $Y$ is algebraic variety (respectively, analytic space),
then the elements of $H^0(V,\ccO_Y)$ can all be written as 
quotients (respectively, Taylor series) of elements of  $H^0(U_i,\ccO_Y)$;
such property holds for instance if $U_i$ is affine or if $U_i$ is biholomorphic to a disk 
$D^{4n} \subset \C^{2n}$ or it is biholomorphic to $D^{4n-2}\times \C^*$; 
\item 
and the ideal $H^0\left(U_i,\I\right) \ideal H^0\left(U_i,\ccO_Y\right)$ 
is preserved by the Poisson bracket.
\end{itemize}
\end{lem}

\begin{prf}
 One implication is obvious, while the other follows from the Leibniz rule
(see Lemma~\ref{lemma_poisson_jacobi_and_leibnitz}) %, 
%from \cite[Prop.~II.2.2(b)]{hartshorne}
  and from elementary properties of coherent sheaves.
\end{prf}

The second lemma asserts that for an isotropic subvariety $X$, 
only functions constant on $X$ can preserve $\I(X)$ by Poisson multiplication.

\begin{lem}\label{lemma_isotropic_then_h_is_constant}
Assume $Y$ is a symplectic manifold, $X$ is a closed irreducible isotropic subvariety.
Let $h \in H^0\left(Y,\ccO_Y\right)$ be any function such that 
\[
\Bigl\{h|_U, H^0\bigl(U,\I(X)\bigr) \Bigr\} \subset H^0\bigl(U,\I(X)\bigr)
\text{ for any open subset $U\subset Y$.}
\]
Then $h$ is constant on $X$.
\end{lem}

\begin{prf}
Choose an arbitrary $x \in X_0$, a small enough open neighbourhood $U \subset Y$ of $x$,
and take any $f \in H^0\bigl(U,\I(X)\bigr)$.

Since $\left\{h|_U,f \right\} \in H^0\bigl(U,\I(X)\bigr)$:
$$
0=\left\{h|_U,f\right\}(x)=\omega(\ud f_x, \ud h_x),
$$
and since $U$ can be taken so small that 
$\left\{\ud f_x \mid f \in H^0\bigl(U,\I(X)\bigr)\right\}$ span the conormal space we have:
$$
\ud h_x \in \left(N^*_x X \slash Y \right)^{\perp_{\omega}}
 \ \stackrel{\textrm{since $X$  is isotropic}}{\subset} \ 
  N^*_x X \slash Y. 
$$
So $\ud h$ vanishes on $TX_0$ and hence $h$ is constant on $X$.
\end{prf}

Finally the third lemma states that for isotropic $X$, very few subvarieties $S\subset X$
 can have the property that $\left\{\I(S), \I(X) \right\} \subset \I(S)$.

\begin{lem}
  \label{lemma_isotropic_and_its_singular_subset}
Assume $Y$ is a symplectic manifold, $X$ is a closed irreducible isotropic subvariety 
and $S \subset X$ is a closed subvariety.
If $\left\{\I(S), \I(X) \right\} \subset \I(S)$, then either $S$ is contained in the singular locus of $X$ or
$X$ is Lagrangian and $S=X$.
\end{lem}

\begin{prf}
 The proof goes along the lines of the proof of \cite[Thm~5.8]{jabu06}.
Suppose $S$ is not contained in the singular locus of $X$, 
so that a general point $s \in S$ is a smooth point of both $X$ and $S$. 
Let $U \subset Y$ be an open neighbourhood of $s$.
Then for all $f\in H^0\left(U,\I(S)\right)$ and  $g \in H^0\left(U,\I(X)\right)$
\begin{equation} \label{equation_for_isotropic_X_and_S}
0 \; = \; \left\{f,g\right\}(s) \; = \; \omega(\ud f_s, \ud g_s),
\end{equation}
so 
\begin{alignat*}{2}
N^*_{s} X \slash Y \quad
&= \quad
\sspan\set{(\ud g)_s \mid g \in H^0\bigl(U,\I(X)\bigr)}&& \\
 &{\subseteq} \quad
\left(N^*_{s} S \slash Y\right)^{\perp_{\omega}} 
 && \text{by \eqref{equation_for_isotropic_X_and_S}}
\\
 &{\subseteq} \quad
\left(N^*_s X \slash Y\right)^{\perp_{\omega}} && \\
&{\subseteq}  \quad
N^*_{s} X \slash Y && \text{because $X$ is isotropic.}
\end{alignat*}
Therefore we have all inclusions becoming equalities 
and in particular 
\[
  \codim S =\codim X
\]
and hence $S=X$. 
Moreover $\left(N^*_{s} X \slash Y\right)^{\perp_{\omega}} = N^*_{s} X \slash Y$, where $s$ is a general point of $X$,
so $X$ is Lagrangian.
\end{prf}

\subsection{Homogeneous symplectic form}

\begin{lem}\label{lemma_homogeneous_symplectic_and_Poisson}
  Assume $(Y,\omega)$ is a symplectic manifold with a $\C^*$-action and 
  that $\omega$ is homogeneous. 
  Let $U\subset Y$ be a $\C^*$-invariant open subset 
  and let $f,g \in H^0(U,\ccO_Y)$ be some homogeneous functions. 
  Then $\{f,g\}$ is homogeneous of weight $\wt (f) + \wt( g) -\wt (\omega)$.
\end{lem}

\begin{prf}
  Let $\xi_g \in H^0(U, TY)$ be such a vector field, that $\omega(\xi_g) = \ud g$.
  By Lemma~\ref{lemma_on_homogeneous_forms_and_vector_fields}\ref{item_omega(mu)_is_homogeneous}
  we have $\wt(\xi_g) = \wt(g)-\wt(\omega)$ and since 
  $\{f,g\} = (\ud f)(\xi_g)$, the claim follows from 
  Lemma~\ref{lemma_on_homogeneous_forms_and_vector_fields}\ref{item_omega(mu)_is_homogeneous}
   and \ref{item_derivative_preserves_weights}.
\end{prf}

\subsection{Example: Veronese map of degree 2}\label{section_example_veronese2}

The following example is important for our considerations,
as it proves that for the contact manifold $\P^{2n-1}$, 
we can equally well consider the Poisson structure on 
$\bigoplus_{i\in \N} \Sym^i \C^{2n}$ (as we do in \cite{jabu_mgr} and \cite{jabu06})
and the Poisson structure on $\bigoplus_{i\in 2\N} \Sym^i \C^{2n}$
(as naturally will arise from the point of view of contact manifolds 
--- see \S\ref{section_sandwich}).
Also this example will be used to illustrate that every contact structure on $\P^{2n-1}$
comes from a symplectic structure on $\C^{2n}$.
Moreover $Y'$ defined below is the minimal adjoint orbit (see \S\ref{section_examples_of_symplectic})
for the simple group $\Sp_{2n}$.
This simple Lie group and its minimal adjoint orbit have quite exceptional behaviour 
(see Table~\ref{table_contacts})
and thus we find it worth to explain this example in more detail.

Let $(V,\omega)$ be a symplectic vector space. 
We let 
$$
\C[V] = \C[x_1, \ldots x_{2n}] = \bigoplus_{i\in \N} \Sym^i V^*
$$ 
be the coordinate ring of $V$.
Also consider 
$$
\ccS := \C[V]^{even} = \bigoplus_{i\in 2\N} \Sym^i V^*
$$
and let $Y':= \Spec \ccS \setminus \{0\}$.
Then we have the following $\Z_2$ covering map:
$$
\psi:V\setminus \set{0} \lra Y',
$$
which is the restriction of the map induced by $\ccS\mono \C[V]$.
This is the underlying map of the second Veronese embedding of $\P(V)$.
In the language of \S\ref{notation_Ldot}, 
we have $Y'= \left(\ccO_{\P(V)}(2)\right)^{\bullet}$ and 
$V\setminus \set{0}= \left(\ccO_{\P(V)}(1)\right)^{\bullet}$.
 
The symplectic form $\omega$ is $\Z_2$ invariant:
$$
\omega(-v,-w) = \omega(v,w),
$$
hence it descents 
to a symplectic form $\omega'$ on $Y'$, making $Y'$ a symplectic manifold, such that:
$$
\psi^*\omega' = \omega.
$$
The natural gradings on $\C[V]$ and on $\ccS$ induce the actions of $\C^*$ on 
$V\setminus \set{0}$ and on $Y'$ 
(note that the action on $Y'$ is not faithful, its kernel is $\Z_2$) 
and $\psi$ is equivariant with respect to  these actions.

\begin{cor}\label{corolary_forms_on_Veronese_2}
With the setup as above,
the form $\omega'$ is homogeneous of weight $2$ 
with respect to the $\C^*$-action described above, so it is of weight $1$ 
with respect to the faithful action of $\C^*/\Z_2 \simeq \C^*$. 
Conversely, if $\omega'$ is a homogeneous symplectic form on $Y'$ of weight $2$, 
then $\psi^*\omega'$ is a constant symplectic form on $V \setminus \set{0}$.
\end{cor}

\begin{prf}
  This follows from Lemma 
  \ref{lemma_on_homogeneous_forms_and_vector_fields}\ref{item_phi_upper_star_preserves_weights}
  and the characterisation of constant forms on an affine space in 
  \S\ref{section_homogeneous_forms}.
\end{prf}

\begin{cor}\label{corollary_Poisson_equal_on_V_and_V_mod_Z2}
The Poisson bracket on $\ccS$ induced by $\omega'$ is the restriction of the Poisson bracket 
on $\C[V]$ induced by $\omega$.
\end{cor}
\begin{prf}
  This follows immediately from Proposition~\ref{proposition_Poisson_compatibile}.
\end{prf}

\chapter{Contact geometry}
\label{chapter_contact}

A projective space seems to be the most standard example of a projective variety. 
Yet, as a contact manifold, 
the projective space of odd dimension is the most exceptional among exceptional examples.
As a consequence, the study of its Legendrian subvarieties is quite complicated 
and very interesting. 
We start our considerations by introducing this case.
Further we generalise to the other contact manifolds.

\section{Projective space as a contact manifold}
\label{section_projective_space}

Let $(V,\omega)$ be a symplectic vector space and let $\P(V)$ be its naive projectivisation.
Then for every $[v] \in \P(V)$ the tangent space $T_{[v]} \P(V)$ 
is naturally isomorphic to the quotient $V \slash [v]$. 
Let $F=F_{\P(V)} \subset T \P(V)$ be a corank 1 vector subbundle defined fibrewise:
$$
F_{[v]}:= \left([v]^{\perp_{\omega}}\right) \slash [v].
$$
Also let $L$ be the quotient line bundle, so that we have the following short exact sequence:
$$
0\lra F \lra T \P(V) \stackrel{\theta}{\lra} L \lra 0.
$$
We say that $F$ (respectively $\theta$) is \textbf{the contact distribution} 
(respectively \textbf{the contact form}) 
\textbf{associated with the symplectic form $\omega$}.

By \S\ref{notation_symplectic_reduction} the vector space $F_p$ 
carries a natural symplectic structure $\omega_{F_p}$.
By Proposition~\ref{properties_of_distribution}\ref{item_dtheta_well_defined}
 $\ud \theta$ gives a well defined twisted 2-form on $F$:
$$
\ud \theta := \Wedge{2} F  \lra L.
$$

\begin{prop}\label{proposition_dtheta_equals_omega}
With an appropriate choice of local trivialisation of $L$, for every $p\in \P(V)$ one has
$\omega_{F_p} = (\ud \theta)_p$. 
In particular $\ud \theta$ is nowhere degenerate and it determines an isomorphism:
$$
F \simeq F^* \otimes L.
$$
Moreover $L \simeq \ccO_{\P(V)}(2)$.
\end{prop}

\begin{prf}
See also \cite[Ex.~2.1]{lebrun}.

Let $x_1\ldots, x_n, y_1, \ldots, y_n$ be some symplectic coordinates on $V$.
Then the \linebreak[3] \mbox{$\omega$-perpendicular} space to $(a_1,\ldots,a_n, b_1, \ldots,b_n)$ is given by the equation 
\[
b_1 x_1+ \ldots+b_n x_n - a_1 y_1 -\ldots - a_n y_n =0.
\] 
We look for a twisted $1$-form $\theta$ on $\P(V)$
whose kernel at each point is exactly as above. 
This is for instance satisfied by 
\[
\theta = 
\half(-y_1\ud x_1 -\ldots - y_n \ud x_n + x_1\ud y_1 +\ldots + x_n \ud y_n ).
\]
The ambiguity is only in the choice of the scalar coefficient 
--- we choose $\half$ in order to acquire the right formula for $\ud \theta$.
Choose an affine piece $U\subset \P(V)$, say where $x_1=1$. 
On $U$ we have 
\[
\theta|_U = 
\half(-y_2\ud x_2 -\ldots - y_n \ud x_n + \ud y_1 + x_2\ud y_2+\ldots + x_n \ud y_n )
\]
and then:
\[
\ud \theta|_U = \ud x_2 \wedge \ud y_2  + \ldots + \ud x_n \wedge \ud y_n.
\]
On the other hand, fixing $p \in U$, $p = [1, a_2, \ldots a_n, b_1, \ldots b_n]$: 
\begin{multline*}
F_p = \Bigl\{(x_1,\ldots,x_n, y_1,\ldots, y_n) \in V \mid 
b_1x_1+ b_2 x_2+ \ldots b_n x_n -\\ 
- y_1- a_2y_2 - \ldots - a_n y_n =0 
\Bigr\} \Big\slash [p].
\end{multline*}
Therefore $F$ is the image under the projection $V \to V/[p]$ of:
\[
  \hat{F}_p := \Bigl\{(0,x_2,\ldots,x_n, \
  a_2 y_2 +\ldots+ a_n y_n - b_2x_2 - \ldots - b_n x_n, y_2, \ldots y_n) \in V \Bigr\}
\]
and
\[
\omega|_{\hat{F}_p} = \ud x_2 \wedge \ud y_2  + \ldots + \ud x_n \wedge \ud y_n.
\]

\smallskip

To see that $L \simeq \ccO_{\P(V)}(2)$ take a section of $T \P(V)$, for instance  $x_1 \frac{\pd}{\pd x_1}$. 
Then 
\[
\theta\left( x_1 \frac{\pd}{\pd x_1} \right)= - \half x_1 y_1
\]
is a section of $L$ and hence $L \simeq \ccO_{\P(V)}(2)$.
\end{prf}

\section{Legendrian subvarieties of projective space}
\label{notation_Legendrian}

Assume $(V, \omega)$ is a symplectic vector space of dimension $2n$.

The author in his previous  articles
\cite{jabu_mgr}, \cite{jabu06}, \cite{jabu_toric}, \cite{jabu_sl} and 
\cite{jabu_hyperplane} found convenient to use  the following definition:

\begin{defin}
  We say that a subvariety $X\subset \P(V)$ is \textbf{Legendrian} 
if the affine cone $\hat{X} \subset V$ is a Lagrangian subvariety
(see \S\ref{section_lagrangian_subvarieties}). 
\end{defin}

Yet the original definition is formulated in a slightly different, 
but equivalent manner:

\begin{prop}\label{prop_equivalent_for_Legendrian_in_projective_space}
Let $X\subset \P(V)$ be a subvariety. 
The following conditions are equivalent:
\begin{itemize}
\item $X$ is Legendrian;
\item $X$ is $F_{\P(V)}$-integrable and it is of pure dimension $n-1$.;
%\item for every smooth point $x \in \hat{X}$  and every $v\in T_x \hat{X}$ 
%we have 
%$$
%\omega(x,v)=0.
%$$
\end{itemize}
\end{prop}

\begin{prf}
If $X$ is $F_{\P(V)}$-integrable, then $X$ is Legendrian by 
Propositions~\ref{proposition_dtheta_equals_omega} and 
\ref{properties_of_distribution}\ref{item_Y_tangent_to_F}.
The other implication is follows from definition of $F_{\P(V)}$. 
\end{prf}

\subsection{Decomposable and degenerate Legendrian subvarieties}
\label{definition_decomposable}

\begin{defin}
Let $V_1$ and $V_2$ be two symplectic vector spaces and
let $X_1\subset \P(V_1)$ and $X_2\subset \P(V_2)$ be two Legendrian subvarieties.
Now assume $V:=V_1 \oplus V_2$ and $X:=X_1 * X_2 \subset \P(V)$, i.e.~$X$ is the join of $X_1$ and $X_2$
meaning the union of all lines from $X_1$ to $X_2$.
Now, clearly, the affine cone $\hat{X}$ is the product $\hat{X_1} \times \hat{X_2}$
(where $\hat{X_i}$ is the affine cone of $X_i$),
so by \S\ref{example_product_Lagrangian} $X$ is Legendrian.
In such a situation we say that $X$ is a \textbf{decomposable Legendrian subvariety}.
We say that a Legendrian subvariety in $\P(V)$ is \textbf{indecomposable},
if it is not of that form for any non-trivial symplectic decomposition $V = V_1 \oplus V_2$.
\end{defin}

The indecomposable Legendrian subvarieties have more consistent description 
of their projective automorphisms group (see Chapter~\ref{chapter_automorphisms}).
On the other hand, decomposable Legendrian varieties 
(which usually themselves are badly singular) 
will be used to provide some very interesting families of examples of smooth Legendrian 
varieties (see Chapter~\ref{chapter_hyperplane}).

\medskip

We say a subvariety of projective space is \textbf{degenerate} 
if it is contained in some hyperplane. Otherwise, 
we say it is \textbf{non-degenerate}.
The following easy proposition in some versions is well known. 
The presented version comes from \cite[Thm~3.4]{jabu06} but see also 
\cite[Prop.~17 (1)]{landsbergmanivel04}
or \cite[Tw.~3.16]{jabu_mgr}.

\begin{prop}\label{proposition_degenerate_legendrian}
Let $X \subset \P(V)$ be a Legendrian subvariety. 
Then the following conditions are equivalent:
\begin{itemize}
\item[(i)]
$X$ is degenerate.
\item[(ii)]
There exists a symplectic linear subspace $W' \subset V$ 
of codimension 2, such that $X' = \P(W') \cap X$ is a Legendrian subvariety 
in $\P(W')$ and $X$ is a cone over $X'$.
\item[(iii)]
$X$ is a cone over some variety $X'$.
\end{itemize}
In particular degenerate Legendrian subvarieties are decomposable.
\end{prop}
\noprf

\section{Contact manifolds}\label{section_contact}

\begin{defin}
Let $Y$ be a complex manifold or smooth algebraic variety and fix a short exact sequence 
$$
0\lra F \lra T Y \stackrel{\theta}{\lra} L \lra 0
$$
where $F\subset T Y$  is a corank $1$ subbundle of the tangent bundle.
We say that $Y$ is \textbf{a contact manifold} 
if the twisted 2-form 
\[
\ud \theta: \Wedge{2} F \lra L
\]
(see Proposition~\ref{properties_of_distribution}\ref{item_dtheta_well_defined})
is nowhere degenerate, so that $\ud \theta_y$ is a symplectic form on $F_y$ for every $y \in Y$. 
In such a case $F$ is called \textbf{the contact distribution} on $Y$
and $\theta$ is  \textbf{the contact form} on $Y$.
\end{defin}

\begin{ex}
By Proposition~\ref{proposition_dtheta_equals_omega}, the projective space with 
the contact distribution associated with a symplectic form is a contact manifold.
\end{ex}

The following properties are standard, well known (see for instance \cite{beauvillefano}):

\begin{prop}\label{properties_of_contact_manifolds}
We have the following properties of contact manifold $Y$:
\begin{enumerate}
\item 
The dimension of $Y$ is odd. 
\item\label{item_F_iso_F_dual_times_L}
Let $U \subset Y$ be an open subset,
let $\mu_F \in H^0(U, F)$ be any section 
and let $\phi_{\mu_F}: F|_U \to L|_U$ be a map of sheaves:
\[
\forall \nu \in  H^0(U,F) \quad \phi_{\mu_F}(\nu):=\theta\bigl([\mu_F, \nu]\bigl). 
\]
Then $\phi_{\mu_F}$ is a map of $\ccO_U$-modules and the assignment 
$\mu_F \mapsto \phi_{\mu_F}$ is an isomorphism of $\ccO_Y$-modules:
\[
F \simeq F^* \otimes L.
\]
\item\label{item_canonical_of_contact}
The canonical divisor $K_Y$ is isomorphic to $L^{\otimes(-n-1)}$.
In particular $Y$ is a Fano variety if and only if $L$ is ample.
\end{enumerate}
\end{prop}

\begin{prf}
We only prove \ref{item_F_iso_F_dual_times_L}, the other parts follow easily.
Map $\phi_{\mu_F}$ is a map of $\ccO_U$-modules by 
Proposition~\ref{properties_of_distribution}\ref{item_bracket_depends_only_on_vectors2}.
By Propsition~\ref{properties_of_distribution}\ref{item_bracket_depends_only_on_vectors}
we have equality:
\[
\phi_{\mu_F}(\nu) = \ud\theta(\mu_F,\nu).
\] 
Since $\ud\theta$ is non-degenerate, it follows that $\mu_F \mapsto \phi_{\mu_F}$ 
is indeed  an isomorphism.
\end{prf}

\subsection{Symplectisation}\label{section_sandwich}

The following construction is standard ---
 see for instance \cite{arnold}, \cite{4authors}, \cite{beauvillefano}. 

Let $\Ldot$ be the principal $\C^*$-bundle as in \S\ref{notation_Ldot}.
In \S\ref{section_automorphisms_of_distributions} and \S\ref{section_thetadot} 
we studied in details the properties of $\Ldot$ 
and an extension of the twisted form $\theta$ to $\Ldot$. 
We have an equivalence between contact structures on $Y$ 
and symplectic homogeneous weight 1 structures on $\Ldot$:

\begin{theo}\label{theorem_symplectisation}
  Let $Y$ be a complex manifold  or smooth algebraic variety 
  with a line bundle $L$  
  and the principal $\C^*$-bundle $\Ldot$ 
  as in \S\ref{notation_Ldot}.
  \begin{itemize}
    \item If $\theta: TY \lra L$ is a contact form,
          then $\ud \thetadot$ (see \S\ref{section_thetadot}) 
          is a homogeneous symplectic form on $\Ldot$ of weight 1.
    \item Conversely, assume $\omega$ is a symplectic form on $\Ldot$,
          which is homogeneous of weight 1. 
          Then there exists a unique contact form $\theta: TY \lra L$ on $Y$, 
          such that $\omega = \ud \thetadot$.
  \end{itemize}
\end{theo}

\begin{prf}
  See Proposition~\ref{prop_on_symplectisation}. 
\end{prf}

If $(Y,F)$ is a contact manifold, then
the symplectic manifold $(\Ldot, \ud \thetadot)$ from the theorem 
is called \textbf{the symplectisation of $(Y,F)$}.

Using the theorem and \S\ref{section_examples_of_symplectic} we have 
following examples of contact manifolds:

\begin{ex} 
  \label{example_adjoint_orbit_is_contact}
  Let $G$ be a simple group and 
  let $Y$ be the closed orbit in $\P(\gotg)$.
  Then $Y$ is a contact manifold (compare with Conjecture~\ref{con_contact}).
\end{ex}

\begin{ex}
  \label{example_projectivised_cotangent}
  If $Y\simeq \P(T^* M)$, then let $L=\ccO_{\P(T^* M)}(1)$ 
  and hence $\Ldot \simeq T^*M \backslash s_0$, 
  where $s_0$ is the zero section and $Y$ is a contact manifold.
\end{ex}

\begin{ex}\label{example_contact_Fano1}
If $Y$ is a contact Fano manifold, then 
\begin{align*}
Y     & \simeq \Proj \left( \bigoplus_{m \in \N} H^0(Y,L^m) \right),\\
\Ldot & \simeq \Spec \left( \bigoplus_{m \in \N} H^0(Y,L^m) \right) \setminus \set{0}
\end{align*}
where $0$ is the point corresponding to the maximal ideal $\bigoplus_{m \ge 1} H^0(Y,L^m)$
(see \cite[\S2.3]{ega2}).
\end{ex}

\begin{ex}\label{example_projective_space_as_contact}
  If $Y\simeq \P(V)$, 
  then by Proposition~\ref{properties_of_contact_manifolds}\ref{item_canonical_of_contact}
  we have 
  \[
    L \simeq \ccO_{\P(V)}(2). 
  \]
  Therefore
  $V \backslash \{0\}$ is a $2$ to $1$ unramified cover of $\Ldot$,
  see \S\ref{section_example_veronese2}. 
  In particular,
  from Theorem~\ref{theorem_symplectisation} 
  and Corollary~\ref{corolary_forms_on_Veronese_2}
  we conclude that every contact structure on $\P(V)$ is associated to  
  some constant symplectic form $\omega$ on $V$ (see \S\ref{section_projective_space}).
\end{ex}

Recall Theorem~\ref{theorem_KPSW_D} that every contact projective manifold 
$Y$ is either isomorphic to $\P(T^*M)$ or it is Fano with $b_2=1$. 
In the second case by 
Proposition~\ref{properties_of_contact_manifolds}\ref{item_canonical_of_contact}
and the Kobayashi-Ochiai characterisation of projective space \cite{kobayashi_ochiai}
either $Y\simeq \P(V)$ or $\Pic Y = \Z \cdot[L]$.

\subsection{Contact automorphisms}

Automorphisms of contact manifolds preserving the contact structure 
were also studied by LeBrun \cite{lebrun} and Beauville \cite{beauvillefano}. 
We use their methods to state slightly more general results about infinitesimal automorphisms.
In the end we globalise the automorphisms for projective contact manifolds.

Let $Y$ be a contact manifold 
and let $\pi: \Ldot \lra Y$ be the symplectisation as in \S\ref{section_sandwich}. 
Also let $\ccR_L$ be as in \S\ref{notation_Ldot}.

\begin{ex}\label{example_contact_Fano2}
  If $Y$ is a contact Fano manifold, then 
  \[
    H^0(Y,\ccR_L)= H^0(\Ldot, \ccO_{\Ldot}) = \left( \bigoplus_{m \in \N} H^0(Y,L^m) \right).
  \]
  Since $Y= \Proj \bigl( H^0(Y,\ccR_L) \bigr)$ 
  (see Example~\ref{example_contact_Fano1}), 
  all the structure of $Y$ as well as its global and local behaviour is determined by 
  this ring of global sections. 
  Hence in this case whatever is stated below for the sheaf $\ccR_L$ can be deduced from 
  the analogous statement about $ H^0(Y,\ccR_L)$ only.
\end{ex}

\pagebreak[2]

\renewcommand{\theenumi}{(\roman{enumi})}
\begin{cor}\label{corollary_Poisson_for_contact}\hfill\nopagebreak
\begin{enumerate}
\item
Let $f,g \in \ccO_{\Ldot}$ be two functions on $\Ldot$ homogeneous with respect to 
the action of $\C^*$. 
Then $\{f,g\}$ is also homogeneous and $\wt \{f,g\} = \wt f + \wt g - 1$
\item
The Poisson bracket descends to $\ccR_L$ and determines a bilinear map:
$$
H^0 (U, L^{m_1}) \times H^0 (U, L^{m_2}) \lra H^0 (U, L^{m_1+m_2-1}). 
$$
\end{enumerate}
\end{cor}

\begin{prf}
This follows from Corollary~\ref{lemma_homogeneous_symplectic_and_Poisson}.
See also \cite[Rem.~2.3]{lebrun}. 
\end{prf}

We will refer to the Lie algebra structure on $\ccR_L$ defined above 
also as Poisson structure and denote the bracket by $\{\cdot,\cdot\}$.
For $s \in H^0(U,L)$ let $\tilde{s}$ be the corresponding element in 
$H^0\left( \pi^{-1}(U), \Ldot \right) = \ccR_L(U)$.

By Corollary~\ref{corollary_Poisson_for_contact} 
the invertible sheaf $L$ has a Lie algebra structure
and it is crucial for our considerations, 
that it is isomorphic to the sheaf $\autinfF$ 
of infinitesimal automorphisms of $Y$ preserving $F$ 
(see \S\ref{section_automorphisms_of_distributions} for more details):
\[
   \autinfF(U) := \left\{ \mu \in H^0(U,TY) \mid [\mu,F] \subset F \right\}.
\]

\renewcommand{\theenumi}{(\arabic{enumi})}
\begin{theo}\label{theorem_contact_automorphisms}
  Let $Y$ be a contact manifold, 
  $F$ be the contact distribution, $\theta$ be the contact form 
  and let $U\subset Y$ be an open subset.
  Using the notation of \S\ref{section_inf_automorphisms} 
  and \S\ref{section_automorphisms_of_distributions}
  we have:
  \begin{enumerate}
    \item\label{item_contact_automorphisms}
       $TY = \autinfF \oplus F$ as sheaves of Abelian groups. 
    \item\label{item_theta_preserves_Lie_bracket}
       The map of sheaves $\theta|_{\autinfF}\colon \autinfF \lra L$ maps isomorphically 
       the Lie algebra structure of $\autinfF$ 
       onto the Lie algebra structure of $L$ given by the Poisson bracket.
    \item\label{item_natural_action}
       The following two Lie algebra representations of $\autinfF$ on $\ccO_{\Ldot}$ are equal:
       \begin{itemize}
         \item 
           The induced representation of $\autinfF$ on $\Ldot$ 
           (see \S\ref{section_automorphisms_of_distributions}).
         \item  
           The representation induced by the adjoint representation:
           $$
             \mu \in \autinfF(U),  f \in H^0(U,\ccO_{\Ldot})  
             \ \Longrightarrow  \ \mu.f := \bigl\{\widetilde{\theta(\mu)},f\bigr\}.
           $$
       \end{itemize}
  \end{enumerate}
\end{theo}

\begin{prf}
The following proof of \ref{item_contact_automorphisms}
is taken from  \cite[Prop.~1.1]{beauvillefano}, 
but see also \cite[Prop.~2.1]{lebrun}.

To prove~\ref{item_contact_automorphisms}, take any $\mu \in H^0(U,TY)$ 
and consider the map of sheaves:
$$
\begin{array}{rcl}
F|_U &\lra &L|_U\\
\nu &\longmapsto& \theta\bigl([\mu,\nu]\bigr).
\end{array}
$$
By Proposition~\ref{properties_of_distribution}\ref{item_bracket_depends_only_on_vectors2} 
the above map is a map of $\ccO_Y|_U$-modules, hence it is an element of 
$H^0(U,F^* \otimes L)$.
Let $\mu_F$ be the corresponding element of $H^0(U,F)$
(see Proposition~\ref{properties_of_contact_manifolds}\ref{item_F_iso_F_dual_times_L}). 
By the definition of the isomorphism $F^* \otimes L \simeq F$, we have 
\[
\theta\left([\mu_F, \nu]\right) = \theta\left([\mu, \nu]\right)
\]
for every $\nu \in F|_U$, hence $[\mu - \mu_F, \nu] \in F|_U$. 
Therefore $\mu - \mu_F \in \autinfF(U)$, so 
$$
\mu = \mu_F + (\mu - \mu_F)
$$
gives the required splitting.

\smallskip

For \ref{item_theta_preserves_Lie_bracket} see also 
\cite[Prop.~1.6]{beauvillefano} and \cite[Rem.~2.3]{lebrun}.
By \ref{item_contact_automorphisms}, 
the map  $\theta|_{\autinfF}$ is an isomorphism of sheaves of Abelian groups. 
So it is enough to prove that $\theta|_{\autinfF}$ preserves the Lie algebra structures.
For every $\mu, \nu \in \autinfF(U)$ denote by $\breve{\mu}$ and $\breve{\nu}$ the
induced infinitesimal automorphisms of $\Ldot$ 
(see \S\ref{section_automorphisms_of_distributions}). 
We have:
\begin{alignat*}{2}
\left\{ 
\widetilde{\theta(\mu)},
\widetilde{\theta(\nu)}
\right\} = & 
(\ud\thetadot)^{\vee}
\biggl(
\ud \left(\widetilde{\theta(\nu)}\right),
\ud \left(\widetilde{\theta(\mu)}\right)
\biggr) \\
=&
(\ud\thetadot)^{\vee}
\big(
\ud \thetadot(\breve{\nu}),
\ud \thetadot(\breve{\mu})
\big)
&& 
\text{by Prop.~\ref{proposition_on_mu_and_thetadot}}\\
=&
\ud\thetadot
\left(
\breve{\nu},
\breve{\mu}
\right)\\
=&
\widetilde{\theta([\mu, \nu])}
&&
\text{by Cor.~\ref{corollary_bracket_of_mu_nu}}.
\end{alignat*}
Hence $\theta|_{\autinfF}$ preserves the Lie algebra structures.

\smallskip

Part~\ref{item_natural_action} is local 
and since both representations satisfy the Leibniz rule 
(see Equation~\eqref{equation_define_mu_tilde} 
and Lemma~\ref{lemma_poisson_jacobi_and_leibnitz}),
 it is enough to check the equality for multiplicative generators of 
$\ccO_{\Ldot}$. Locally, these might be taken for instance as sections of $L$ 
and so~\ref{item_natural_action} follows from~\ref{item_theta_preserves_Lie_bracket}.

\end{prf}

We underline, that $\autinfF$, as a subsheaf of $TY$
is not a $\ccO_Y$-submodule (see \S\ref{section_automorphisms_of_distributions}).
So in particular the obtained splitting of the short exact sequence
of sheaves of Abelian groups
\[
0\lra F \lra T Y \stackrel{\theta}{\lra} L \lra 0
\]
is not a splitting of vector bundles.

Turning to global situation assume $Y$ is projective 
and let $\Aut(Y)$, $\Aut_F(Y)$ and $\gotaut(Y)$, $\gotaut_F(Y)$ 
denote, respectively, 
the group of automorphisms of $Y$,
the group of automorphisms of $Y$ preserving the contact structure 
and their Lie algebras.

LeBrun \cite{lebrun} and Kebekus \cite{kebekus_lines1} 
observed that in the case of projective contact Fano manifolds 
with Picard group generated by $L$,
the global sections of $L$ are isomorphic to $\gotaut(Y)$:

\renewcommand{\theenumi}{(\roman{enumi})}
\begin{cor}\label{corollary_AutY_eq_AutFY}
Let $Y$ be a projective contact manifold with contact distribution $F$.\begin{enumerate}
\item\label{item_theta_maps_autF_onto_H0L}
Then $\theta$ maps isomorphically $\gotaut_F(Y)$  onto $H^0(Y, L)$.
\item \label{item_H0L_tangent_to_AutY}
If moreover $Y$ is Fano with $\Pic(Y) = \Z [L]$, 
then $\Aut(Y) = \Aut_F(Y)$ and hence the Lie algebra $H^0(Y, L)$
is naturally isomorphic to $\gotaut(Y)$.
\end{enumerate}
\end{cor}

\begin{prf}
  By Corollary~\ref{corollary_Y_projective_infinitesimal_and_tangent_to_global}
  we have $\gotaut_F(Y)=\autinfF(Y)$, 
  so~\ref{item_theta_maps_autF_onto_H0L} follows from 
  Theorem~\ref{theorem_contact_automorphisms}~\ref{item_theta_preserves_Lie_bracket}.

\smallskip

  On the other hand \ref{item_H0L_tangent_to_AutY} follows from \cite[Cor.~4.5]{kebekus_lines1}.
\end{prf}

\section{Legendrian subvarieties in contact manifold}
\label{section_legendrian_automorphisms}

\begin{defin}
Let $Y$ be a complex contact manifold with a contact
 distribution $F$.
A subvariety $X\subset Y$ is \textbf{Legendrian} if $X$ is $F$-integrable
 (i.e., $TX \subset F$) and  $2\dim X +1 = \dim Y$ (i.e., $X$ has maximal
 possible dimension).
\end{defin}

If $Y\simeq \P^{2n+1}$,
then the above definition agrees with the definition in \S\ref{notation_Legendrian} by
Proposition~\ref{prop_equivalent_for_Legendrian_in_projective_space}. 
In general, we have analogous properties with $V$ replaced by $\Ldot$:

\renewcommand{\theenumi}{(\alph{enumi})}

\begin{prop}\label{proposition_isotropic_integrable}
Let $Y$ be a contact manifold with a contact
 distribution $F \subset T Y$ and with its symplectisation $\pi:\Ldot \to Y$.
Assume $X\subset Y$ is a subvariety. Then:
\begin{enumerate}
 \item
     \label{item_isotorpic_integrable}
   $X$ is $F$-integrable if and only if $\pi^{-1}(X) \subset \Ldot$ is  isotropic.
 \item
     \label{item_Legendrian_integrable}
   $X$ is Legendrian if and only if $\pi^{-1}(X) \subset \Ldot$ is Lagrangian.
\end{enumerate}
\end{prop}

\begin{prf}
  Part~\ref{item_isotorpic_integrable} follows from 
  Lemma~\ref{lemma_integrable_iff_dthetadot_vanish} and 
  Part~\ref{item_Legendrian_integrable} follows from 
  \ref{item_isotorpic_integrable}.
\end{prf}

In the case of subvarieties of a symplectic manifold, 
we have three important types of subvarieties (isotropic, Legendrian and coisotropic).
Also for subvarieties of contact manifold in addition to 
$F$-integrable and Legendrian subvarieties, it is useful to 
consider the subvarieties corresponding to the coisotropic case:

\begin{defin}
In the setup of Proposition~\ref{proposition_isotropic_integrable},
 we say that $X$ is \textbf{$F$-cointegrable}
if $\pi^{-1}(X) \subset \Ldot$ is coisotropic.
\end{defin}

\begin{ex}
  Assume $\widetilde{X} \subset \Ldot$ is irreducible and Lagrangian 
and let $X$ be the closure of $\pi(X) \subset Y$. Then $X$ is $F$-cointegrable. 
If moreover $\widetilde{X}$ is not $\C^*$-invariant, then $\dim X = \half (\dim Y +1)$.
\end{ex}

\begin{cor}\label{cor_Legendrian_in_cotangent}
  If $Y = \P(T^*M)$ for some smooth algebraic variety $M$ and $X$ is an algebraic 
  Legendrian subvariety, 
  then $X$ is the conormal variety $Z^{\#}$ to some algebraic subvariety $Z\subset M$. 
\end{cor}
\begin{prf}
  It follows from Proposition~\ref{proposition_isotropic_integrable}, 
  Example~\ref{example_projectivised_cotangent} and Lemma~\ref{lemma_Lagrangian_in_cotangent}.
\end{prf}

Let $\ccR_L=\pi_* \ccO_{\Ldot}$ be the sheaf of rings on $Y$ 
defined in \S\ref{notation_Ldot}.
For a subvariety $X\subset Y$, let $\widetilde{\I}(X) \ideal \ccR_L$ 
be the sheaf of ideals generated by those local sections of $L^m$ that vanish on $X$.
Then:
\begin{equation}
  \label{equation_pi_push_forward_ideal}
  \pi_*\I\left(\pi^{-1}(X)\right) = \widetilde{\I}(X)
\end{equation}
where $\I\left(\pi^{-1}(X)\right) \ideal \ccO_{\Ldot}$ is the ideal sheaf of $\pi^{-1}(X)$.
In this context, the meaning of Lemma~\ref{lemma_poisson_is_affine} is the following:

\begin{lem}\label{lemma_poisson_is_affine2}
  With the notation as above, let $\I\ideal \ccO_{\Ldot}$ be a coherent sheaf of ideals.
Then $\I$ is preserved by the Poisson bracket on $\ccO_{\Ldot}$ if and only if 
$\pi_*\I$ is preserved by the Poisson bracket on $\ccR_L$.
\end{lem}
\noprf

Hence we get the description of $F$-cointegrable subvarieties 
in terms of the Poisson bracket on $\ccR_L$:

\begin{prop}\label{proposition_cointegrable_Poisson}
With the assumptions as above, a subvariety $X\subset Y$ 
is  $F$-cointegrable if and only if $\widetilde{\I}(X)$ is preserved by the Poisson bracket on $\ccR_L$.
\end{prop}

\begin{prf}
The proposition combines 
Equation~\eqref{equation_pi_push_forward_ideal}, 
Theorem~\ref{Poisson_structure} 
together with Lemma~\ref{lemma_poisson_is_affine2}. 
\end{prf}

Given a subvariety $X\subset Y$,
we define $\autinfF(\cdot,X)$ to be the sheaf of Lie algebras 
of those infinitesimal automorphisms of $Y$, 
which preserve $X$ and contact distribution $F$ 
(see also \S\ref{section_automorphisms_of_distributions}):
\begin{multline*}
\autinfF(U,X) :=
\Bigl\{
  \mu \in H^0(U, TY) \mid  [\mu,F] \subset F   \text{ and } \\
  \forall f \in \I(X)|_U \ \ 
 (\ud f)(\mu) \in \I(X)|_U
\Bigr\}
\end{multline*}

Further, let $\widetilde{\I}(X)_1 \subset L$ be the degree $1$ part of the sheaf 
of homogeneous ideals $\widetilde{\I}(X)$.
Since $L$ is a line bundle with the action of $\autinfF$ 
(see \S\ref{section_automorphisms_of_distributions}),
 choosing a local trivialisation and using the gluing property of sheaves 
we can replace $\I(X)$ in the definition of  $\autinfF(\cdot,X)$
with $\widetilde{\I}(X)_1$:
\begin{multline}\label{equation_for_autinfFX}
\autinfF(U,X) =
\Bigl\{
  \mu \in H^0(U, TY) \mid   [\mu,F] \subset F  \text{ and}\\
  \mu.\widetilde{\I}(X)_1|_U\subset \widetilde{\I}(X)_1|_U
\Bigr\}
\end{multline}
where $.$~denotes the induced action of $\autinfF$ on $L$ described in 
\S\ref{section_automorphisms_of_distributions}.

The following theorem establishes a connection between the infinitesimal automorphisms 
of a Legendrian variety and its ideal:

\renewcommand{\theenumi}{\Alph{enumi}}
\renewcommand{\labelenumi}{\theenumi.}

\begin{theo}\label{theorem_isomorphisms_of_subvarieties}
Let $Y$ be a contact manifold with a contact
 distribution $F$ and let $\theta: TY\to L$ be the quotient map. 
Also let $U\subset Y$ be an open subset.
Assume $X\subset Y$ is an irreducible subvariety.
%and let $\widetilde{\I}(X)_1 \subset L$ be the degree $1$ part of the sheaf 
%of homogeneous ideals $\widetilde{\I}(X)$.
\begin{enumerate}
\item \label{item_inf_automorphisms_integrable}
If $X$ is $F$-integrable, then 
$\theta\left(\autinfF(U,X)\right) \subset H^0\left(U,\widetilde{\I}(X)_1\right)$.
\item \label{item_inf_automorphisms_cointegrable}
If $X$ is $F$-cointegrable, then 
$\theta\left(\autinfF(U,X)\right) \supset H^0\left(U,\widetilde{\I}(X)_1\right) $.
\item \label{item_inf_automorphisms_Legendrian}
If $X$ is Legendrian, then $
\theta\left(\autinfF(U,X)\right) =  H^0\left(U,\widetilde{\I}(X)_1\right) $.
\end{enumerate}
\end{theo}
\renewcommand{\labelenumi}{\theenumi}

\begin{prf}
In the case of \ref{item_inf_automorphisms_integrable},
choose arbitrary $\mu \in \autinfF(U,X)$. 
We must prove that $\theta(\mu) \in H^0(U, \widetilde{\I}_1(X))$ or, equivalently, 
that 
\[
\widetilde{\theta(\mu)} \in H^0\Bigl(\pi^{-1}(U), \I\bigl(\pi^{-1}(X)\bigl)\Bigr)
\] 
(recall that for  a section $s\in H^0(U,L)$ 
by $\tilde{s}$
we denote
the corresponding element in 
$H^0(\pi^{-1}(U), \ccO_{\Ldot})$).

By Equation~\eqref{equation_for_autinfFX} the action of $\mu$ preserves $\widetilde{\I}(X)|_U$ 
and hence also it preserves $\I\bigl(\pi^{-1}(X)\bigl)|_{\pi^{-1}(U)}$.
By Theorem~\ref{theorem_contact_automorphisms}\ref{item_natural_action} 
this means that 
\[
\left\{
\widetilde{\theta(\mu)} , \  \I\bigl(\pi^{-1}(X)\bigr)|_{\pi^{-1}(U)} 
\right\} \ \  
\subset  \ \ 
\I\bigl(\pi^{-1}(X)\bigr)|_{\pi^{-1}(U)}.
\]
Moreover $\pi^{-1}(X)$ is isotropic by Proposition~\ref{proposition_isotropic_integrable}.

By Lemma~\ref{lemma_isotropic_then_h_is_constant} 
function $\widetilde{\theta(\mu)}$ is constant on $\pi^{-1}(X)$. 
But $\widetilde{\theta(\mu)}$ is also a $\C^*$-homogeneous function of weight 1, 
so it  must vanish on 
$\pi^{-1}(X)$.
Therefore $\widetilde{\theta(\mu)} \in H^0\Bigl(\pi^{-1}(U), \I\bigl(\pi^{-1}(X)\bigl)\Bigr)$
as claimed.

\smallskip

To prove \ref{item_inf_automorphisms_cointegrable}
let $\mu\in \autinfF(U)$ be an infinitesimal automorphism such that 
$\theta(\mu) \in \widetilde{\I}(X)_1$.
By Proposition~\ref{proposition_cointegrable_Poisson} 
$$ 
\left\{\theta(\mu),\widetilde{\I}(X) \right\} \subset \widetilde{\I}(X)
$$ 
so by Theorem~\ref{theorem_contact_automorphisms}~\ref{item_natural_action}
we have
\[
\mu.\widetilde{\I}(X) \subset  \widetilde{\I}(X)
\]
(where $.$~denotes the induced representation of $\autinfF$ on $\Ldot$, 
see \S\ref{section_automorphisms_of_distributions}).
Hence by Equation~\eqref{equation_for_autinfFX} the infinitesimal automorphism 
$\mu$ is contained in $\autinfF(U,X)$ and 
$\ H^0\left(U,\widetilde{\I}(X)_1\right) \subset \theta\left(\autinfF(U,X)\right)$ as claimed.

\smallskip

Part~\ref{item_inf_automorphisms_Legendrian} is an immediate consequence of 
\ref{item_inf_automorphisms_integrable} 
and \ref{item_inf_automorphisms_cointegrable}.
\end{prf}

The following corollary says that in the case when $Y$ is projective 
also the global automorphisms of a Legendrian subvariety 
can be understood in terms of the ideal of the variety. 
In particular, in \ref{item_Y_projective_autF_eq_I} below, 
we generalise Theorem~\ref{theorem_ideal_and_group}.

\renewcommand{\theenumi}{(\roman{enumi})}
\begin{cor}
  \label{corollary_Y_projective_autF_eq_I}
Let $Y$ be a projective contact manifold,
let $F$ be the contact distribution
and let $X$ be a Legendrian subvariety. 
Let $\gotaut(Y,X)$ 
(resp.~$\gotaut_F(Y,X)$) 
be the Lie algebra of group of automorphisms of $Y$ preserving $X$ 
(resp.~preserving $X$ and $F$). 
Then:
\begin{enumerate}
\item 
    \label{item_Y_projective_autF_eq_I}
$\theta \bigl(\gotaut_F(Y,X)\bigr) = H^0\left(Y,\widetilde{\I}(X)_1\right)$;
\item
    \label{item_Y_not_proj_space_aut_X_eq_I}
If in addition $\Pic Y = \Z [L]$, 
then $\theta \left(\gotaut(Y,X)\right) = H^0\left(Y,\widetilde{\I}(X)_1\right)$.
\end{enumerate} 
\end{cor}

\begin{prf}
  It follows from Corollary~\ref{corollary_AutY_eq_AutFY}
  and
  Theorem~\ref{theorem_isomorphisms_of_subvarieties}\ref{item_inf_automorphisms_Legendrian}.
\end{prf}

In Chapter~\ref{chapter_automorphisms} we discuss the extension of 
Corollary~\ref{corollary_Y_projective_autF_eq_I}\ref{item_Y_not_proj_space_aut_X_eq_I} 
to $Y\simeq \P^{2n+1}$.

\medskip

The following corollary generalises \cite[Thm 5.8]{jabu06}:

\begin{cor}
  If $Y$ is a projective contact manifold and 
  $X\subset Y$ is an irreducible Legendrian subvariety such that 
  $\widetilde{\I}(X)$ is generated by $H^0\bigl(Y,\widetilde{\I}(X)_1\bigr)$,
  then $\Aut_F(Y,X)$ acts transitively on the smooth locus of $X$.
  In particular, if $X$ is in addition smooth, then $X$ is a homogeneous space.
\end{cor}

\begin{prf}
  If $S\subset X, S\neq X$ is a closed subvariety invariant under the action of $\Aut_F(Y,X)$,
  then by Theorem~\ref{theorem_contact_automorphisms}~\ref{item_natural_action} and 
  by Corollary~\ref{corollary_Y_projective_autF_eq_I}\ref{item_Y_projective_autF_eq_I}: 
  $$
    \forall f \in H^0\left(Y,\widetilde{\I}(X)_1\right) 
    \quad \left\{\widetilde{\I}(S),f\right\} \subset \widetilde{\I}(S).
  $$ 
  Hence by the Leibniz rule and since $\widetilde{\I}(X)$ 
  is generated by $H^0\bigl(Y,\widetilde{\I}(X)_1\bigr)$, we have:
  $$
     \biggl\{\I\Bigl(\pi^{-1}(S)\Bigr),\I\Bigl(\pi^{-1}(X)\Bigr)\biggr\} 
     \subset \I\Bigl(\pi^{-1}(S)\Bigr).
  $$
  So by Lemma~\ref{lemma_isotropic_and_its_singular_subset},
  variety  $S$ is contained in the singular locus of $X$.

  Now let $O\subset X$ be an orbit of a smooth point under the action of $\Aut_F(Y,X)$. 
  Then the closure $\overline{O}$ is not contained in the singular locus 
  so by above it must be equal to all of $X$. 
  Moreover $\overline{O} \setminus O$ is a closed subset  invariant under the action and 
  not equal to $X$, so it is contained in the singular locus. 
  So $O$ is the  whole smooth locus of $X$.
\end{prf}

We conclude this chapter by underlining that, unfortunately, 
the above results are proved only for automorphisms of $Y$, 
that preserve Legendrian subvariety $X$, not simply for automorphisms of $X$.

\chapter{Projective automorphisms of a Legendrian variety}
\label{chapter_automorphisms}

The content of this chapter is published in  \cite{jabu_toric}.

We are interested in the following conjecture:

\begin{con}\label{theorem_automorphisms_are_symplectic}
Let $X\subset \P^{2n-1}$ be an irreducible indecomposable 
Legendrian subvariety and let $G< \P\Gl_{2n}$ be a connected subgroup of linear automorphisms 
preserving $X$. 
Then $G$ is contained in the image of the natural map $\Sp_{2n} \ra \P\Gl_{2n}$.
\end{con}

It is quite natural to believe, that if a linear map preserves a form on a big number of linear
subspaces, then it actually preserves the form (at least up to scalar). 
With this approach, Janeczko and Jelonek
 proved the conjecture in the case
 where the image of $X$ under the Gauss map is 
non-degenerate in the Grassmannian of Lagrangian subspaces in $\C^{2n}$
 --- see \cite[Cor.~6.4]{janeczkojelonek}.
Unfortunately, this is not enough -  for example $\P^1 \times Q_1 \subset \P^5$ 
has a degenerate image under the Gauss map 
and this is one of the simplest examples of smooth Legendrian subvarieties.

In \S\ref{section_proof_of_the_conjecture} we prove 
Theorem~\ref{corollary_that_conjecture_is_true_for_smooth},
i.e.~that the conjecture holds for smooth $X$.
This theorem, combined with Corollary~\ref{corollary_Y_projective_autF_eq_I} 
gives us a good understanding of the group of projective automorphisms of 
a smooth Legendrian subvariety in $\P^{2n-1}$.

\section{Discussion of assumptions}

One obvious remark is that homotheties act trivially on $\P(V)$, 
but in general are not symplectic on $V$.
Therefore, it is more convenient to think of conformal symplectomorphisms 
(see \S\ref{section_symplectic_automorphisms}).

It is clear, that if we hope for a positive answer 
to the question whether a projective automorphism of a Legendrian subvariety 
necessarily preserves the contact structure, 
then  we must assume that our Legendrian variety is non-degenerate. 

Another natural assumption is that $X$ is irreducible 
--- one can also easily produce a counterexample if we skip this assumption.
Yet still this is not enough.

Let $X= X_1* X_2\subset \P(V_1 \oplus V_2)$ be a decomposable Legendrian variety. 
Then we can act via $\lambda_1 \Id_{V_1}$ on $V_1$ and 
via $\lambda_2 \Id_{V_2}$ on $V_2$ - such an action will preserve $X$ 
and in general it is not conformal symplectic. 
This explains why the assumptions of our conjecture \ref{theorem_automorphisms_are_symplectic}
are necessary.

%%%%%%%%%%%%%%%%%%%%%%%

%%%%%%%%%%%%%%%%%%%%%%%

%%%%%%%%%%%%%%%%%%%%%%%%%%%%%%

\section{Preservation of contact structure}

\label{section_proof_of_the_conjecture}

Let $X'\subset \P(V)$ be an irreducible, indecomposable Legendrian subvariety, 
let \mbox{$X$ be} the affine cone over $X'$  and $X_0$ be the smooth locus of $X$. 
Assume that $G$ is the maximal connected subgroup in $\Gl_{2n}$
preserving $X$. Let $\gotg< \gotgl_{2n}$ be the Lie algebra tangent to  $G$.
To prove the conjecture it would be enough to show that $\gotg$ is contained in the Lie algebra $\gotcsp_{2n}$
tangent to conformal symplectomorphisms, i.e.~the Lie algebra spanned by $\gotsp_{2n}$ and the identity
matrix $\Id_{2n}$.

Recall from \S\ref{section_symplectic_automorphisms} the notion of 
weks-symplectic matrices.

\begin{theo}\label{theorem_evidence_for_conjecture}
With the above notation the following properties hold:
\begin{itemize}
\item[I.]
The underlying vector space of $\gotg$ decomposes into symplectic and weks-symplectic part:
$$
\gotg = \big(\gotg \cap \gotsp(V)\big)  \oplus \big(\gotg \cap \gotasp(V)\big). 
$$ 
\item[II.]
If $g \in \gotg \cap \gotasp(V)$, then $g$ preserves every tangent space to $X$:
$$
\forall {x\in X_0} \quad g(T_x X) \subset T_x X
$$
and hence also
$$
\forall t\in \C \quad \forall x\in X_0 \qquad T_{\exp(t g)(x)} X = \exp(t g)(T_x X) = T_x X.
$$
\item[III.]
If $g \in \gotg \cap \gotasp(V)$ is semisimple,
then $g=\lambda \Id$ for some $\lambda \in \C$.
\item[IV.]
Assume $0 \ne g \in \gotg \cap \gotasp(V)$ is nilpotent and let $m \ge 1$ be an integer such that $g^{m+1} = 0$
and $g^{m}\ne 0$.
Then $g^m(X)$ is always non-zero and
is contained in the singular locus of $X$.
In particular, $X'$ is singular.
\end{itemize}
\end{theo}

In what follows we prove the four parts of Theorem~\ref{theorem_evidence_for_conjecture}.

\subsubsection{I.~Decomposition into symplectic and weks-symplectic part}
\begin{prf}
Take $g \in \gotg$ to be an arbitrary element. 
Then for every $x \in X_0$  one has 
$$g (x) \in T_x X$$
and therefore
$$
0 = \omega \big( g(x), x \big) = x^T g^T J x  = \half x^T \left( g^T J - J g \right) x.
$$

Hence the quadratic polynomial  $f(x):=x^T ( g^T J - J g) x$ is identically zero on $X$ and hence it is in the ideal of $X$.
Therefore by maximality of $G$ and Theorem~\ref{theorem_ideal_and_group}
the map $J \left(g^T J - J g\right)$ is also in $\gotg$. However,
$$
J \left(g^T J - J g\right) = J g^T J + g,
$$
so $Jg^T J \in \gotg$ and both symplectic and weks-symplectic components $g_+$ and $g_-$ are in $\gotg$.
\end{prf}

From the point of view of the conjecture, the symplectic part is fine.
We would only need to prove that $g_- = \lambda \Id$.
So from now on we assume $g = g_- \in \gotasp(V)$.

\subsubsection{II.~Action on tangent space}
\begin{prf}
Let $\gamma_t:=\exp(t g)$ for $t \in \C$.
Then $\gamma_t \in G$ and hence it acts on $X$.
Choose a point $x\in X_0$ and two tangent vectors in the same tangent space $u,v \in T_x X$.
Then clearly also $\gamma_t(u)$ and $\gamma_t(v)$ are contained in one tangent space, 
namely $T_{\gamma_t(x)} X$. Hence:
\begin{align*}
0=& \omega \left(\gamma_t(u),\gamma_t(v)\right)\\
= & \omega \bigl(
(\Id_{2n} + t g  +\ldots)u, 
(\Id_{2n} + t g  +\ldots)v
\bigr)  \\
= & \omega(u,v) + t \bigl( \omega(g u, v) + \omega(u, g v)\bigr) + t^2(\ldots).
\end{align*}
In particular the part of  the expression linear in $t$ vanishes, hence $\omega(g u, v) + \omega(u, g v) =0$.
Combining this with Equation~\eqref{weks_symplectic_equation} we get that:
$$
\omega(g u, v) = \omega(u, g v) = 0.
$$
However, this implies that $g u \in (T_x X)^{\perp_{\omega}} = T_x X$.
Therefore $g$ preserves the tangent space at every smooth point of $X$ 
and hence also $\gamma_t$ preserves that space.
\end{prf}

\subsubsection{III.~Semisimple part}
Since $G$ is an algebraic subgroup in $\Gl(V)$, 
hence $\gotg$ has the natural Jordan decomposition inherited from $\gotgl(V)$, 
i.e.~if we write the Jordan decomposition for $g= g_s +g_n$, then $g_s, g_n \in \gotg$ 
(see \cite[Thm~15.3(b)]{humphreys}).
Therefore by Proposition~\ref{properties_of_asp}(i),
proving that
for $g\in \gotg \cap \gotasp(V)$ we have $g_s = \lambda \Id_{2n}$ and $g_n=0$
 would be enough to establish the conjecture.

Here we deal with the semisimple part.

\begin{prf}
Argue by contradiction.
Let $V_1$ be an arbitrary eigenspace of $g$ and let $V_2$ be the sum of the other eigenspaces. 
If $g \ne \lambda \Id_{2n}$, then both $V_1$ and $V_2$ are non-zero and by 
Proposition~\ref{properties_of_asp}(ii) and (iii) they are $\omega$-perpendicular, 
complementary symplectic subspaces of $V$. Let $x\in X_0$ be any point. 
Since $g$ preserves $T_x X$ by part II it follows
that  $T_x X= (T_x X\cap V_1) \oplus (T_x X \cap V_2)$. 
But then both $(T_x X\cap V_i) \subset V_i$ are Lagrangian subspaces, hence have constant 
(independent of $x$) dimensions.
Hence $T_x X_0 = (T_x X_0 \cap V_1) \oplus (T_x X_0 \cap V_2)$ is a sum of two vector bundles 
and  from Proposition~\ref{proposition_decomposable} we get that 
$X$ is a product of two Lagrangian subvarieties, 
contradicting our assumption that $X'$ is indecomposable.
\end{prf}

%%%%%%%%%%%%%%%%%

\subsubsection{IV.~Nilpotent part --- $X'$ is singular}

\begin{lem}\label{limit_of_exptgv}
Assume $X' \subset \P(V)$ is any closed subvariety preserved by the action of $\exp(t g)$
for some nilpotent endomorphism $g\in \gotgl(V)$.
If $v$ is a point of the affine cone over $X'$ and $m$ is an integer such
that $g^{m+1}(v)=0$ and $g^{m}(v) \ne 0$, then $[g^{m}(v)]\in \P(V)$ is in $X'$.
\end{lem}

\begin{prf}
Point $[g^{m}(v)]  \in \P(V)$ is just the limit of $[\exp(t g)(v)]$ as $t$ 
goes to $\infty$.
\end{prf}

\begin{lem}\label{smooth_point_implies_linear}
Assume $g\in \gotgl(V)$ is nilpotent and
 $g^{m+1}=0$, $g^{m} \ne 0$ for an integer $m\ge 1$.
Let $X\subset V$ be an affine cone over some irreducible projective subvariety in $\P(V)$,
which is preserved by the action of $\exp(t g)$, but is not contained in the set of the fixed points. 
Assume that this action preserves the tangent space $T_x X$ at every smooth point $x$ of $X$.
If there exists a non-zero vector in $V$ which is a smooth point of $X$ contained in $g^{m}(X)$,
then $X$ is a linear subspace.
\end{lem}

\begin{prf}
\emph{Step~0 - notation.}
We let $Y$ to be the closure of $g^m(X)$, so in particular $Y$ is irreducible.
By Lemma~\ref{limit_of_exptgv},
we know that $Y \subset X$.
Let $y$ be a general point of $Y$.  
Then by our assumptions $y$ is a smooth point of both $X$ and $Y$.

Next denote by 
$$
W_y:= (g^m)^{-1} (\C^* y).
$$  
You can think of $W_y$ as union of those lines in $V$
(or points in the projective space $\P(V)$), 
which under the action of $\exp(t g)$
converge 
to the line spanned by $y$ (or $[y]$)\footnote{
This statement is not perfectly precise, though it is true on an open dense subset.
There are some other lines, which converge to $[y]$, namely those generated by $v \in \ker g^m$, 
but $g^k(v) = \lambda y$ for some $k<m$.
We are not interested in those points.
}
as $t$ goes to $\infty$ .
We also note that the closure $\overline{W_y}$ is a linear subspace spanned by an arbitrary element 
$v \in W_y$ and $\ker g^m$.

Also we let $F_y:= W_y \cap X$, so that:
$$
F_y:= (g^m|_X)^{-1} (\C^* y).
$$
Finally, $v$ from now on will always denote an arbitrary point of $F_y$.

\medskip

\emph{Step~1 - tangent space to $X$ at points of $F_y$.}
Since $y$ is a smooth point of $X$ also $F_y$ consists of smooth points of $X$. 
This is because the set of singular points is closed and $\exp(t g)$ invariant. 
By our assumptions  $\exp(t g)$ preserves every tangent space to $X$ 
and thus for every $v \in F_y$ we have:
$$
T_v X = T_{\frac{1}{t^m}\exp(t g)(v)} X = T_{\lim_{t\ra \infty} \left(\frac{1}{t^m}\exp(t g)(v)\right)} = T_y X.
$$
So the tangent space to $X$ is constant over the $F_y$ and in particular $F_y \subset T_y X$.

\medskip

\emph{Step~2 - dimensions of $Y$ and  $F_y$.}
From the definitions of $Y$ and $y$ and by Step~1 we get that for any point $v\in F_y$:
$$
T_y Y \ = \ \im (g^m|_{T_v X}) =   \im (g^m|_{T_y X}).
$$
Hence $\dim Y = \dim T_y Y = \rk (g^m|_{T_y X})$.

Since $y$ was a general point of $Y$, we have that:
$$
\dim Y + \dim F_y = \dim X +1.
$$
So $\dim F_y = \dim \ker(g^m|_{T_y X}) +1$.

\medskip

\emph{Step~3 - the closure of $F_y$ is a linear subspace.}
From the definition of $F_y$ and Step~1
we know that $F_y \subset T_y X \cap W_y$ and 
$$
 T_y X \cap \overline{W_y} =
 T_y X \cap \sspan\{v, \ker g^m \} = 
\sspan\{v, \ker (g^m|_{T_y X}) \}.
$$
Hence $\dim F_y = \dim T_y X \cap W_y$,
so the closure of $F_y$ is exactly $ T_y X \cap \overline{W_y}$ and clearly this closure is contained in $X$.
In particular $\ker (g^m|_{T_y X}) \subset X$. 

\medskip

\emph{Step~4 - $Y$ is contained in $\ker (g^m|_{T_y X})$.}
Let $Z$ be $X \cap \ker g^m$.
By Step~3 we know that $\ker (g^m|_{T_y X}) \subset Z$. 
Now we calculate the local dimension of $Z$ at $y$:
$$
\dim \ker (g^m|_{T_y X}) \le
\dim_y Z \le \dim T_y Z \le 
\dim (T_y X \cap \ker g^m) = \dim \ker (g^m|_{T_y X}). 
$$
Since the first and the last entries are identical, we must  have all equalities. 
In particular the local dimension of $Z$ at $y$ is equal to the dimension of the tangent space to $Z$ at $y$.
So $y$ is a smooth point of $Z$ and therefore there is a unique component of $Z$ passing through $y$, 
namely the linear space $\ker (g^m|_{T_y X})$. 
Since $Y$ is contained in $Z$ (because $\im g^m \subset \ker g^m$) and $y \in Y$, 
we must have $Y \subset  \ker (g^m|_{T_y X})$.

\medskip

\emph{Step~5 - vary $y$.}
Recall, that by Step~1 the tangent space to $X$ is the same all over $F_y$. 
So also it is the same on every smooth point of $X$, 
which falls into the closure of $F_y$. 
But by Step~4, $Y$ is a subset of $ \ker (g^m|_{T_y X})$, 
which is in the closure of $F_y$ by Step~3. 
So the tangent space to $X$ is the same for an open subset of points in $Y$.
Now apply again Step~1 for different $y$'s in this open subset 
and we get that $X$ has constant tangent space on a dense open subset of $X$.
This is possible if and only if $X$ is a linear subspace, which completes the proof of the lemma.
\end{prf}

Now part IV of the theorem follows easily:

\begin{prf}
By the assumptions of the theorem $X$ is  not contained in any hyperplane, 
so in particular $X$ is not contained in $\ker g^m$. 
So by Lemma~\ref{limit_of_exptgv} the image $g^m(X)$ contains points other than $0$.
Next by Lemma~\ref{smooth_point_implies_linear} and part II of the theorem, since $X$ cannot be a linear subspace,
there can be no smooth points of $X$ in  $g^m(X)$.  
\end{prf}

\subsubsection{Smooth case}
We conclude that parts I, III and IV of Theorem~\ref{theorem_evidence_for_conjecture}
together with Proposition~\ref{properties_of_asp}(i) and \cite[Thm.~15.3(b)]{humphreys}
imply Theorem~\ref{corollary_that_conjecture_is_true_for_smooth}. 
We only note that a smooth Legendrian subvariety is either a linear subspace or it is indecomposable.

%%%%%%%%%%%%%%%%%%%

%%%%%%%%%%%%%%%%%%%%%
%%%%%%%%%%%%%%%%%%%%%
%%%%%%%%%%%%%%%%%%%%%
%%%%%%%%%%%%%%%%%%%%%

\section{Some comments}\label{section_comments}

Conjecture~\ref{theorem_automorphisms_are_symplectic} is now reduced  
to the following special case not covered by Theorem~\ref{theorem_evidence_for_conjecture}:

\begin{con}
Let $X'\subset \P(V)$ be an irreducible Legendrian subvariety.
Let $g \in \gotasp(V)$ be a nilpotent endomorphism  and $m$ be  an integer such that 
\mbox{$g^m\ne 0$} and $g^{m+1}=0$. 
Assume that the action of $\exp(t g)$ preserves $X'$.
Assume moreover, that $X'$ is singular at points of the image of the rational 
map $g^m(X')$. 
Then $X'$ is decomposable.
\end{con}

We also note the improved relation between projective automorphisms 
of a Legendrian subvariety and quadratic equations satisfied by its points:

\begin{cor}
Let $X\subset \P(V)$ be an irreducible Legendrian subvariety 
for which Conjecture~\ref{theorem_automorphisms_are_symplectic} holds
(for example $X$ is smooth).
If $G < \P\Gl(V)$ is the maximal subgroup preserving $X$, 
then $\dim G = \dim \I_2(X)$, where $\I_2(X)$ is the space of homogeneous quadratic 
polynomials vanishing on $X$.
\end{cor}
\begin{prf}
It follows immediately from the statement of the conjecture and 
Theorem~\ref{theorem_ideal_and_group}.
\end{prf}

Finally,
it is important to note, that Theorem~\ref{theorem_evidence_for_conjecture} part III does not imply 
that every torus acting on an indecomposable, but singular Legendrian variety  $X'$ 
is contained in the image of $\Sp(V)$. 
It only says that the intersection of such a torus with the weks-symplectic part is always finite. 
Therefore if there is a non-trivial torus acting on $X'$, 
there is also some non-trivial connected subgroup of $\Sp(V)$ acting on $X'$ 
and also some quadratic equations in the ideal of $X'$.

\chapter{Toric Legendrian subvarieties in projective space}\label{chapter_toric}

The content of this chapter is published in  \cite{jabu_toric}.

We apply Theorem~\ref{corollary_that_conjecture_is_true_for_smooth}
to classify smooth toric Legendrian subvarieties. 
We choose appropriate coordinates to 
reduce this problem to some combinatorics (for surface case --- see \S\ref{section_toric_surfaces})
and some elementary geometry of convex bodies (for higher dimensions --- 
see \S\ref{section_higher_dimensional}). 
Eventually the proof of Theorem \ref{theorem_classification_of_smooth_toric_legendrian}
is obtained in Corollaries~\ref{smooth_toric_surfaces} and \ref{smooth_toric_varieties}.

\section{Classification of toric Legendrian varieties}

Within this chapter $X$ is a toric subvariety of dimension $n-1$ in a projective space of dimension $2n-1$.
We assume it is embedded torically, so that the action of $T:=(\C^*)^{n-1}$ on $X$ 
extends to an action on the whole $\P^{2n-1}$,
but we do not assume that the embedding is projectively normal. 
The notation is based on \cite{sturmfels} though we also use techniques of \cite{oda}.
We would like to understand when $X$ can be Legendrian with respect to some contact structure
on $\P^{2n-1}$ and in particular, when it can be a smooth toric Legendrian variety.

There are two reasons for considering non projectively normal toric varieties here. 
The first one is that the new example we find is not projectively normal. 
The second one is the conjecture \cite[Conj.~2.9]{sturmfels}, 
which says that a smooth, toric, projectively normal variety is defined by quadrics. 
We do not expect to produce a counterexample to this conjecture 
and on the other hand all smooth Legendrian varieties defined by quadrics are known
to be just the subadjoint varieties
(see \cite[Thm~5.11]{jabu06}).

\medskip

In addition we assume that either $X$ is smooth or at least the following condition is satisfied:
\begin{itemize}\label{star_condition}
\item[($\star$)]
The action of the torus $T$ on $\P^{2n-1}$ preserves the contact structure on $\P^{2n-1}$.
In other words, the image of $T\ra \P\Gl_{2n}$ is contained in the image of $\Sp_{2n} \ra \P\Gl_{2n}$.
\end{itemize}

In the case where $X$ is smooth, the ($\star$) condition is always satisfied
by Theorem~\ref{corollary_that_conjecture_is_true_for_smooth}. 
But for some statements below we do not need non-singularity, so we only assume ($\star$).

\begin{theo}\label{theorem_toric_legendrian}
Let $X\subset \P^{2n-1}$ be a toric (in the above sense) non-degenerate Legendrian subvariety
satisfying ($\star$).
Then there exists a choice of symplectic coordinates on $\C^{2n}$
and coprime integers $a_0 \ge a_1 \ge \ldots \ge a_{n-1} > 0$
such that $X$ is the closure of the image of the following map:
\begin{multline*}
T \ni (t_1, \ldots, t_{n-1}) 
\mapsto
[-a_0 t_1^{a_1} t_2^{a_2} \ldots t_{n-1}^{a_{n-1}}, \ \
a_1 t_1^{a_0}, \  a_2 t_2^{a_0},\  
\ldots,  \ a_{n-1}t_{n-1}^{a_0},\\
 t_1^{-a_1} t_2^{-a_2} \ldots t_{n-1}^{-a_{n-1}}, \ \ 
 t_1^{-a_0}, \ t_2^{-a_0},\  
\ldots, \  t_{n-1}^{-a_0}] \in \P^{2n-1}.
\end{multline*}
In other words, $X$ is the closure of the orbit of a point 
$$
[-a_0, a_1, a_2, \ldots, a_{n-1}, 1, 1, \ldots 1]\in \P^{2n-1}
$$
under the torus action with weights 
\begin{gather*}
w_0:=(a_1, a_2,\ldots, a_{n-1}),
\\
w_1:=(a_0,0,\ldots, 0), \
w_2:=(0,a_0,0,\ldots, 0), \
\ldots, \
w_{n-1}:=(0,\ldots, 0, a_0) 
\\
\text{and } 
-w_0, -w_1, \ldots, -w_{n-1}.
\end{gather*}

Moreover every such $X$ is a non-degenerate toric Legendrian subvariety.
\end{theo}

We are aware that for many choices of the $a_i$'s from the theorem, the action of the torus on 
$X$ (and on $\P^{2n-1}$) is not faithful, 
so that for such examples a better choice of coordinates could be made. 
However,
 we are willing to pay the price of taking a quotient of $T$ to get a uniform description.
One advantage of the description given in the theorem is that a part of it is almost 
independent of the choice of the $a_i$'s. 
This part is the $(n-1)$-dimensional ``octahedron'' 
$\conv\{w_1,\ldots w_{n-1}, -w_1, \ldots -w_{n-1}\}\subset \Z^{n-1} \otimes \R$.

\begin{prf}
Assume $X$ is Legendrian with respect to a symplectic form $\omega$, that $X$ is non-degenerate,
that the torus $T$ acts on $\P^{2n-1}$ preserving $X$ and satisfies ($\star$).
Replacing if necessary $T$ by some covering we may assume that $T\ra \P\Gl_{2n}$
factorises through a maximal torus $T_{\Sp_{2n}} \subset \Sp_{2n}$:
$$
T\ra T_{\Sp_{2n}} \subset \Sp_{2n} \ra \P\Gl_{2n}.
$$

This implies, that for an appropriate symplectic basis 
the variety $X$ is the closure of the image of the map $T\ra \P^{2n-1}$ given by:
$$
T \ni t \mapsto [x_0 t^{w_0},x_1 t^{w_1}\ldots, x_{n-1}t^{w_{n-1}}, t^{-w_0}, t^{-w_1}\ldots, t^{-w_{n-1}}] \in \P^{2n-1}
$$
where $x_i \in \C$, $w_i \in \Z^{n-1}$ and for
 $v= (v_1,\ldots v_{n-1})\in \Z^{n-1}$ we let \mbox{$t^v := t_1^{v_1} \ldots t_{n-1}^{v_{n-1}}$}.
 This means that $X$ is the closure of the $T$-orbit of the point%
\footnote{Note
that usually one assumes that this point is 
just [1,\ldots,1]. 
In our case we would have to consider non-symplectic coordinates. 
We prefer to deal with a point with more complicated coordinates.
} 
$[x_0, \ldots x_{n-1}, 1, \ldots, 1]$
where $T$ acts with  weights 
$w_0, \ldots w_{n-1}, -w_0,\ldots, -w_{n-1}$.

Since $X$ is non-degenerate, the weights are pairwise different. 
Also the weights are not contained in any hyperplane in $\Z^{n-1}\otimes \R$,
because the dimension of $T$ is equal to the dimension of $X$ and we assume $X$ has an open orbit of
 the $T$-action. 
So there exists exactly one (up to scalar) linear relation:
$$
 - a_0 w_0 + a_1 w_1 + \ldots + a_{n-1} w_{n-1} = 0 .
$$
We assume that the $a_i$'s are coprime integers.
Permuting coordinates appropriately we can assume that $|a_0| \ge |a_1| \ge \ldots \ge |a_{n-1}| \ge 0$.
After a symplectic change of coordinates, we can assume without loss of generality that all the $a_i$'s are non negative by exchanging 
$w_i$ with $-w_i$ (and $x_i$ with $-\frac{1}{x_i}$) if necessary. 
Clearly not all the $a_i$'s are zero, so in particular $a_0>0$ and hence
$$
w_0 = \frac{a_1 w_1 + \ldots + a_{n-1} w_{n-1}}{a_0}.
$$
Therefore, if we set $e_i:= \frac{w_i}{a_0}$ for $i \in \{1, \ldots, n-1 \}$, 
the points $e_i$ form a basis of a lattice $M$ containing all $w_i$'s. 
The lattice $M$ might be finer than the one generated by the $w_i$'s. 
Replacing again $T$ by a finite cover, 
we can assume that the action of $T$ is expressible in the terms of weights in $M$.
Then:
\begin{align*}
w_0     & =  a_1 e_1 + \ldots + a_{n-1} e_{n-1},\\
w_1     & =  a_0 e_1,\\
        & \vdots,\\
w_{n-1} & =  a_0 e_{n-1}.
\end{align*}

It remains to prove three things: that $a_{n-1}>0$, 
that the $x_i$'s might be chosen as in the statement of the theorem 
and finally that every such variety is actually Legendrian.
We will do all three together. 

The torus acts symplectically on the projective space, 
thus the tangent spaces to the affine cone are Lagrangian if and only if 
just one tangent space at a point of the open orbit is Lagrangian.
So take the point $[x_0, \ldots x_{n-1}, 1, \ldots, 1]$. 
The affine tangent space is spanned by the following vectors:
\begin{alignat*}{9}
v:=      &(&   x_0,&&   x_1,&&x_2,&\ldots,&x_{n-1},&\ & 1,&&       1,&&    1,&    \ldots,& 1),\\
u_1:=    &(&x_0a_1,&&x_1a_0,&&0,  &\ldots,&      0,&  &-a_1,&&     -a_0,&& 0,&    \ldots,& 0),\\
u_2:=    &(&x_0a_2,&&     0,&&x_2a_0,&\ldots,& 0,  &  &-a_2,&&      0,&&  -a_0,&  \ldots,&  0),\\
&\vdots \\
u_{n-1}:=&(&x_0a_{n-1},&&  0,&& 0,&\ldots,& x_{n-1}a_0,&  &-a_{n-1},&&  0,&&    0,& \ldots,& -a_0).\\
\end{alignat*}

Now the products are following:
\begin{align*}
\omega(u_i, u_j) &=0;  \\
\omega(u_i, v)   &= 2(x_0 a_i + x_i a_0).
\end{align*}
Therefore the linear space spanned by $v$ and the $u_i$'s is Lagrangian if and only if:
$$
x_i = -x_0 \frac{a_i}{a_0}.
$$
In particular, since $x_i \ne 0$, the $a_i$ cannot be zero either.
After another conformal symplectic base change,
we can assume that $x_0= -a_0$ and then $x_i = a_i$.
 On the other hand, the above equation is satisfied for the variety in the
theorem. Hence the theorem is proved.
\end{prf}

Our next goal is to determine for which values of the $a_i$'s the variety $X$ is smooth. 
The curve case is not interesting at all and also very easy, so we start from $n=3$, 
i.e.~Legendrian surfaces.

\section{Smooth toric Legendrian surfaces}\label{section_toric_surfaces}

We are interested in knowing when the toric projective surface with weights of torus action
\begin{align*}
w_0:=&(a_1,a_2),&  w_1:=&(a_0, 0), & w_2:=&(0, a_0),
\\
-w_0 =& (-a_1, -a_2), & -w_1=&(-a_0,0),&-w_2 =& (0,-a_0) 
\end{align*}
  is smooth.
Our assumptions on the $a_i$'s are following:
\begin{equation}\label{inequalities_for_surface}
a_0\ge a_1\ge a_2 >0
\end{equation}
and $a_0, a_1,a_2$ are coprime integers.

\begin{figure}[htb]
\centering
\includegraphics[width=0.3\textwidth]{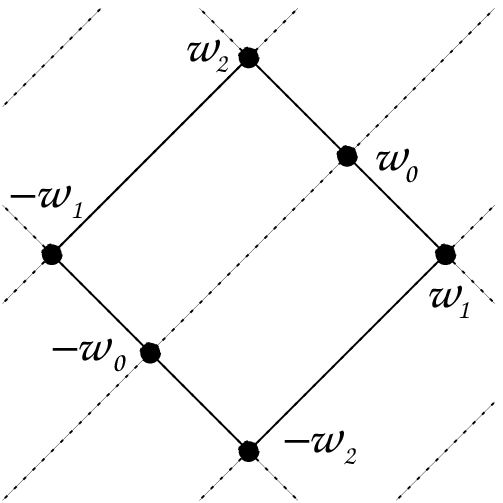}
\hspace{1cm}
\includegraphics[width=0.3\textwidth]{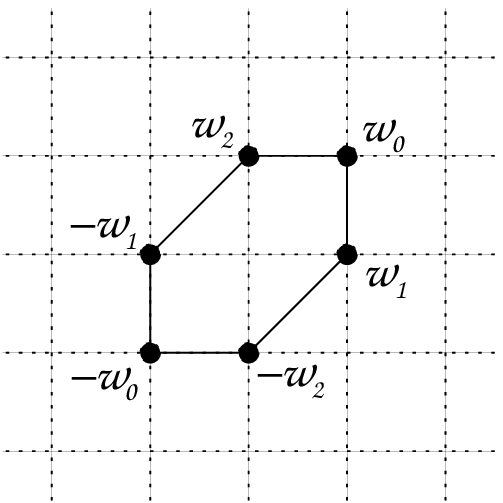}
\caption{
\footnotesize
The two examples of weights giving smooth toric Legendrian surfaces. 
}
 \label{figure_surface_examples}
\end{figure}

\begin{ex} \label{example_line_times_conic}
Let $a_0=2$ and  $a_1=a_2=1$ (see Figure~\ref{figure_surface_examples}). Then $X$ is the product of $\P^1$ and a quadric plane curve $Q_1$.
\end{ex}

\begin{ex}\label{example_plane_blown_up}
Let $a_0=a_1=a_2=1$ (see Figure~\ref{figure_surface_examples}).
Although the embedding is not projectively normal (we lack the weight $(0,0)$ in the middle),
the image is smooth anyway.
Then $X$ is the blow up of $\P^2$ in three non-colinear points.
\end{ex}

We will prove there is no other smooth example.

\medskip

\begin{figure}[htb]
\centering
\includegraphics[width=0.4\textwidth]{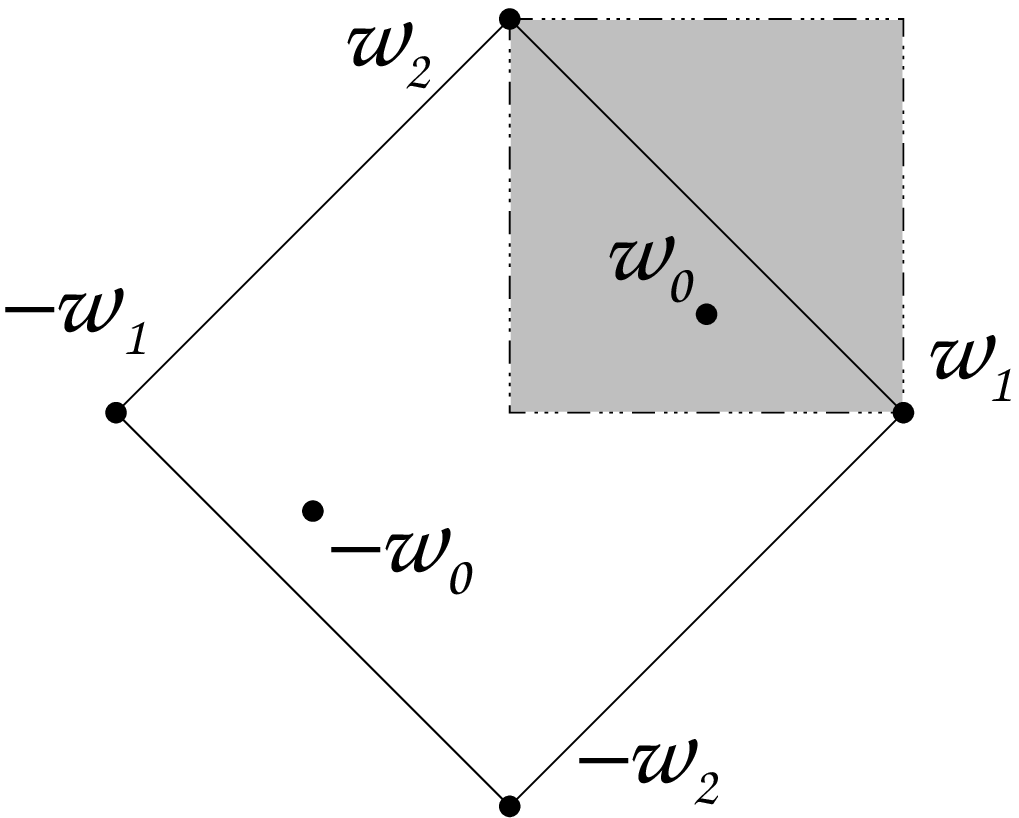}
\hspace{1cm}
\includegraphics[width=0.4\textwidth]{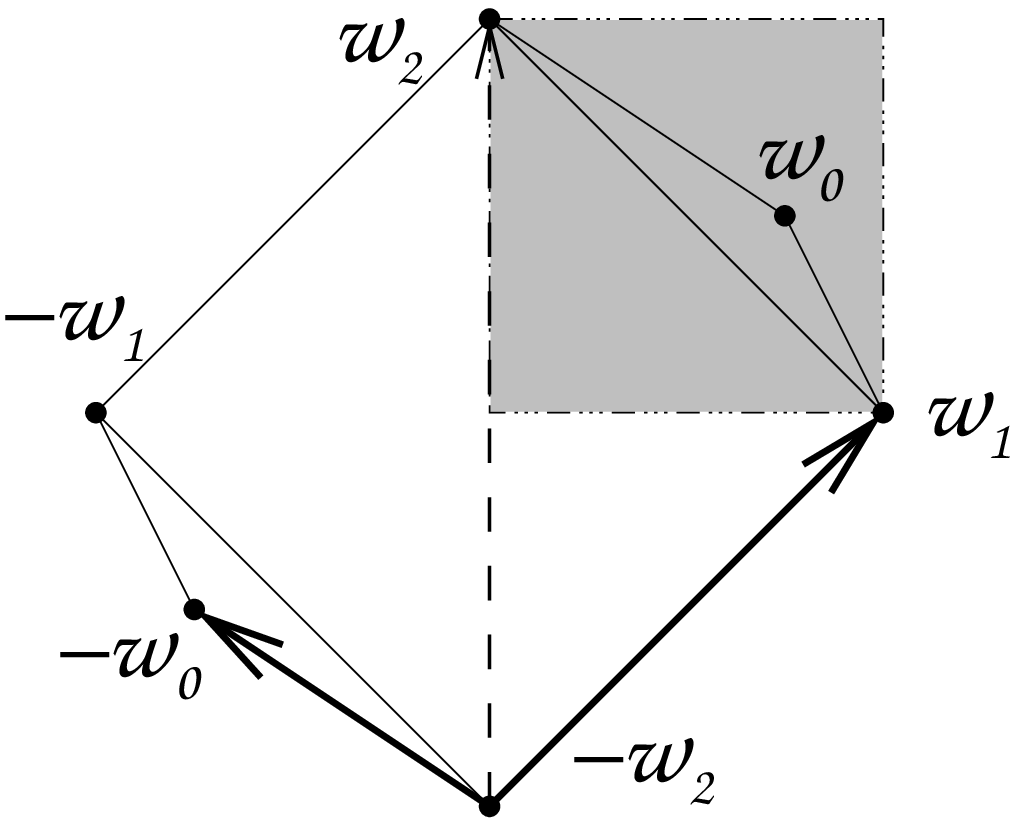}
\caption{
\footnotesize
Due to the inequalities $a_0 \ge a_1> 0$  and $a_0 \ge a_2> 0$,  
the weight $w_0$ is located somewhere in the grey square. 
The two cases we consider are if $w_0$ is also inside the square $\conv \{ w_1,w_2, -w_1, -w_2 \}$ 
(left figure) or it is outside (right figure). 
In the second case, a necessary condition to get a smooth variety, is that the two bold vectors generate a lattice containing all the weights. 
In particular the dashed vector can be obtained as an integer combination of the bold ones.
}
 \label{figure_surface_proof}
\end{figure}

We must consider two cases (see Figure~\ref{figure_surface_proof}):
either $a_0 > a_1 +a_2 $ (which means that $w_0$ is in the interior of the square 
$\conv \{ w_1,w_2, -w_1, -w_2 \}$) or 
$a_0 \le a_1+a_2$ (so that $w_0$ is outside or on the border of the square).

\smallskip

Geometrically, case $a_0 > a_1 +a_2 $ means, that the normalisation of $X$ is $\P^1\times \P^1$. 
It is just an easy explicit verification 
that $X$ is not smooth with these additional weights in the interior.

\smallskip

In the other case, 
for a vertex $v$ of the polytope 
$$
\conv \{ w_0,w_1,w_2, -w_0,-w_1, -w_2 \},
$$ 
we define the sublattice $M_v$ to have the origin at $v$ and to be generated by
$$
\{ w_0-v,w_1-v,w_2-v, -w_0-v,-w_1-v, -w_2-v\} .
$$
Since $X$ is smooth, for every vertex $v$
the vectors of the edges meeting at $v$ 
must form a basis of $M_v$ 
(compare with \cite[Prop.~2.4 \& Lemma~2.2]{sturmfels}). 
In particular, if $v=-w_2$ 
(it is immediate from Inequalities~\eqref{inequalities_for_surface} that $v$ is indeed a vertex), 
then $w_2-(-w_2) = (0, 2 a_0)$ can be expressed as an integer combination of 
$w_1+w_2= (a_0,a_0)$ and $-w_0+w_2 = (-a_1, a_0-a_2)$ 
(see the right hand side of Figure~\ref{figure_surface_proof}).
So write:
\begin{equation}\label{equation_for_surface_case}
(0, 2 a_0) = k (a_0,a_0) + l (-a_1, a_0-a_2)
\end{equation}
for some integers $k$ and $l$.
 It is obvious that $k$ and $l$ must be strictly positive, 
since $w_2$ is in the cone generated by $w_1+w_2$ and  $-w_0+w_2$ with the vertex at $-w_2$.
But then (since $a_0-a_2 \ge 0$) 
from Equation~\eqref{equation_for_surface_case} on the second coordinate we get that either
 $k=1$ or $k=2$.

\smallskip

If $k=1$, then we easily get that:
$$
\left\{
\begin{array}{l}
a_0 = l a_1 \\
a_0 = a_1 +a_2.
\end{array}
\right.
$$
Hence $(l-1)a_1 = a_2$ and by Inequalities~(\ref{inequalities_for_surface}) we get $l=2$ 
and therefore (since the $a_i$'s are coprime)
$(a_0,a_1, a_2) =(2,1,1)$, which is  Example~\ref{example_line_times_conic}.

\smallskip

On the other hand, if $k=2$, then 
$$
a_0 = a_2 
$$
and hence by Inequalities~(\ref{inequalities_for_surface}) and since the $a_i$'s are coprime,
 we get $(a_0,a_1, a_2) =(1,1,1)$, which is Example~\ref{example_plane_blown_up}.

\begin{cor}\label{smooth_toric_surfaces}
If $X\subset \P^5$ is smooth toric Legendrian surface, 
then it is either $\P^1 \times Q_1$ 
or $\P^2$ blown up in three non-colinear points or plane $\P^2\subset \P^5$.
\end{cor}
\noprf

\section{Higher dimensional toric Legendrian varieties}
\label{section_higher_dimensional}

In this section we assume that $n\ge 4$. 
By means of the geometry of  convex bodies we will prove there is only one smooth 
toric non-degenerate Legendrian variety   
in dimension $n-1=3$ and no more in higher dimensions.
We use Theorem~\ref{theorem_toric_legendrian} so that we have a toric variety with weights:
\begin{align*}
&w_0:=(a_1,a_2,\ldots, a_{n-1}),\\
&w_1:=(a_0,0,\ldots 0),\\
&\vdots\\
&w_{n-1}:=(0,\ldots 0, a_0),\\
&-w_0, -w_1, \ldots, -w_{n-1}
\end{align*}
 
where the $a_i$'s are coprime positive integers with $a_0\ge a_1 \ge \ldots \ge a_{n-1}$.

\begin{figure}[htb]
\centering
\includegraphics[scale=0.7]{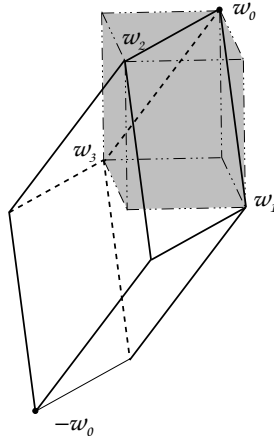}
\caption{
\footnotesize
The smooth example in dimension 3: $(a_0,a_1,a_2,a_3)= (1,1,1,1)$.
}
 \label{octahedron_case4}
\end{figure}

\begin{ex}\label{example_line_cubed}
Let $n=4$ and $(a_0,a_1,a_2,a_3)= (1,1,1,1)$. 
Then the related toric variety is $\P^1 \times \P^1 \times \P^1$ 
(see Figure~\ref{octahedron_case4}).
\end{ex}

Further, let $A$ be the polytope defined by the weights:
$$
A:=\conv \{ w_0,w_1,\ldots , w_{n-1}, -w_0, -w_1, \ldots, -w_{n-1}\} \subset \Z^{n-1} \otimes \R.
$$

\begin{lem}\label{edges_of_A}
Let $I,J \subset \{1, \ldots, n-1\}$ be two complementary subsets of indexes. 
\begin{itemize}
\item[(a)]
Assume $i_1, i_2 \in I$ and $ i_1 \ne i_2$. If 
$$
\left| \sum_{i \in I} a_i -\sum_{j \in J} a_j \right| < a_0 ,
$$
then the interval $(w_{i_1},w_{i_2})$ is an edge of $A$.
\item[(b)]
Assume $k \in I$ and $l \in J$. If
$$
\sum_{i \in I} a_i -\sum_{j \in J} a_j  > a_0, 
$$
then both intervals $(w_0,w_k)$ and $(w_0, -w_l)$ are edges of $A$.
\item[(c)]
If $k,l \in \{1, \ldots, n-1\}$ and $k \ne l$, then  $(w_k, -w_l)$ is an edge of $A$.
\end{itemize} 
\end{lem}

\begin{prf}
Fix $\epsilon >0$ small enough,
set $\alpha:= \sum_{i \in I} a_i -\sum_{j \in J} a_j$
and define the following hyperplanes in $\Z^{n-1} \otimes \R $:
$$ 
H_a:= \left\{ 
\sum_{i \in I} x_i - (1-\epsilon)  \sum_{j \in J} x_j  = a_0 
\right\},
$$
$$
H_b:=
\left\{
\left(a_0-a_k\right)\left(\sum_{i \in I} x_i - \sum_{j \in J} x_j -\alpha \right)
+(\alpha -a_0)\left(x_k-a_k\right) =  0
\right\},
$$
$$
H'_b:=
\left\{
\left(a_0+a_l\right)\left(\sum_{i \in I}  x_i - \sum_{j \in J} x_j - \alpha  \right)
+(\alpha -a_0)\left(x_l+a_l\right) =  0
\right\}
$$
$$
\textrm{and }
H_c:=
\left\{
x_k-x_l =  a_0
\right\}.
$$

Assuming the inequality of (a), 
$H_a \cap A$ is equal to $\conv\{w_i \mid i \in I \}$ and 
the rest of $A$ lies on one side of $H_a$.
So $H_a$ is a supporting hyperplane for the face \mbox{$\conv\{w_i \mid i \in I \}$},
which is a simplex of dimension $(\#I-1)$ 
and therefore all its edges are also edges of $A$ as claimed in (a).

\smallskip

Next assume that the inequality of (b) holds. 
Then $H_b$ (respectively $H'_b$) is a supporting hyperplane for the edge $(w_0,w_k)$
(respectively $(w_0,-w_l)$).  

\smallskip

Similarly, in the case of (c), $H_c$ is a supporting hyperplane for $\{w_k, -w_l\}$.

\end{prf}

\begin{theo}\label{classification_toric_higher_dimensional}
Let $X\subset \P^{2n-1}$ be a toric non-degenerate Legendrian variety of dimension $n-1$
satisfying ($\star$) (see page~\pageref{star_condition}).
If $n\ge 4$ and normalisation of $X$ has at most quotient singularities,
 then $n=4$ and $X=\P^1 \times \P^1 \times \P^1$.
\end{theo}

\begin{prf}
Since the normalisation of $X$ has at most quotient singularities,
it follows that the polytope $A$ is simple, 
i.e.~every vertex has exactly $n-1$ edges (see \cite{fulton} or \cite[\S2.4, p.~102]{oda}).
We will prove this is impossible, unless $n=4$ and $(a_0, a_1, a_2,a_3)= (1,1,1,1)$.
%We must consider several positions of $w_0$ relative to the ``octahedron'' 
%$B:=\conv \{w_1,\ldots , w_{n-1}, -w_1,\ldots -w_{n-1}\}$.

\smallskip

If  $w_0 \in B:=\conv \{w_1,\ldots , w_{n-1}, -w_1,\ldots -w_{n-1}\}$,
then $A$ is just equal to $B$ 
 and clearly in such a case every vertex of $A$ has  $2(n-2)$ edges. 
Hence more than $n-1$ for $n\ge 4$.

\smallskip
 
Thus from now on we can assume that $a_1 + \ldots + a_{n-1} > a_0$. 
yB Lemma~\ref{edges_of_A}(b), $(w_0, w_i)$ is an edge for every $i \in \{1, \ldots, n-1\}$.

Choose any $j \in \{1, \ldots, n-1\}$ and set $I:=\{1, \ldots,j-1, j+1, \ldots, n-1\}$.

If either  
$$
\left|\left(\sum_{i \in I} a_i\right) \ -  \ a_j \right| \ < \ a_0 \quad \textrm{ or}
$$
$$
\left(\sum_{i \in I} a_i\right) \ -  \ a_j \  > \ a_0, 
$$
then using Lemma~\ref{edges_of_A} we can count the edges at either $w_i$ or $w_0$ 
and see that there is always more than $n-1$ of them. 
We note that $a_j -\left(\sum_{i \in I} a_i\right) \ge a_0$ never happens due to our assumptions on the $a_i$'s.

\smallskip

Therefore the remaining case to consider is
$$
\left(\sum_{i \in I} a_i\right) \  - \  a_j \   = \  a_0, 
$$
where the equality holds for every $j \in \{1, \ldots,n-1\}$.
This implies:
$$
a_1 =a_2 =\ldots = a_{n-1} = \frac{1}{n-3} a_0.
$$
Since the $a_i$'s are positive integers and coprime, we must have 
$$
(a_0,a_1, \ldots, a_{n-1}) = (n-3, 1,\ldots, 1)
$$ 
which is exactly Example~\ref{example_line_cubed} for $n=4$.
Otherwise, if  $n\ge 5$ we can take $J:=\{j_1,j_2\}$ for any two different 
$j_1,j_2 \in \{1, \ldots, n-1\}$ and set $I$ to be the complement of $J$.
Then $\# I \ge 2$ and  by Lemma~\ref{edges_of_A}(a) and (c) there are too many edges 
at the $w_i$'s.
\end{prf}

\begin{cor}\label{smooth_toric_varieties}
If $X\subset \P^{2n-1}$ is a smooth toric Legendrian subvariety 
and $n\ge 4$,
then it is either a linear subspace or $n=4$ and $X= \P^1 \times \P^1 \times \P^1$.
\end{cor}
\noprf

\chapter{Examples of quasihomogeneous Legendrian varieties}
\label{chapter_sl}

The content of this chapter is published in \cite{jabu_sl}.

We construct a family of examples of Legendrian subvarieties in
projective spaces.
Although most of them are singular,
a new example of a smooth Legendrian variety in dimension 8 is in this family.
The 8-fold has interesting properties:
it is a compactification of the special linear group,
a Fano manifold of index 5 and Picard number 1
 (see Theorem~\ref{theorem_classify_invertible}(b)).
 Also we show how  this
construction generalises to give new smooth examples in dimensions 5 and 14 
(see \S\ref{section_other_examples}).

% and
%finally we announce a result which will produce plenty of such examples
%(see \S\ref{section_hyperplane}).

\medskip

In \S\ref{section_basics2} we introduce the notation for this chapter.
In \S\ref{results} we formulate the results and make some comments on possible generalisations.
In \S\ref{section_action} we study the structure of a group action related to the problem.
In \S\ref{main_part} we finally prove the results.

\section{Notation and definitions}\label{section_basics2}

For this chapter we fix an integer $m\ge 2$.

\subsubsection{Vector space $V$}\label{notation_V}
\label{notation_begins}
Let $V$ be a vector space over complex numbers $\C$ of dimension $2m^2$,
 which we interpret as a space of pairs of $m \times m$ matrices. 
The coordinates are: $a_{i j}$ and $b_{i j}$ for $i,j \in \{1,\ldots m\}$. 
By $A$ we denote the matrix $(a_{i j})$ and similarly for $B$ and $(b_{i j})$.

\medskip

Given two $m\times m$ matrices $A$ and $B$, by $(A,B)$ we denote the point of the vector space $V$, 
while by $[A,B]$ we denote the point of the projective space $\P(V)$.

\medskip

Sometimes, we will represent some linear maps $V\lra V$ and
some 2-linear forms $V\otimes V \lra \C$ as $2m^2 \times 2m^2$ matrices. In such a case
we will assume the coordinates on $V$ are given in the lexicographical
order:
$$
a_{11}, \ldots, a_{1m}, a_{21}, \ldots, a_{mm}, b_{11}, \ldots, b_{1m},
b_{21}, \ldots, b_{mm}.
$$

\subsubsection{Symplectic form $\omega$}
\label{notation_omega}
On $V$ we consider the standard symplectic form 
%$\omega:= \sum \ud a_{ij} \wedge \ud b_{ij}$  so that:
\begin{equation}
 \label{properties_of_inner_product}
\omega\big((A,B), (A', B')\big):= 
\sum_{i,j} (a_{i j}b'_{i j} -a'_{i j}b_{i j}) = \tr \left( A(B')^T - A' B^T \right).
\end{equation}
Further we set $J$ to be the
matrix of $\omega$:
$$
J:=M(\omega)= 
\left[ 
\begin{array}{cc}
0          & \Id_{m^2}\\
-\Id_{m^2} & 0
\end{array}
\right].
$$

\subsubsection{Varieties $Y$, $\Xinv(m)$ and $\Xdeg(m,k)$}
\label{notation_Y}

We consider the following subvariety of $\P(V)$:
\begin{equation}
\label{defining_equations_of_Y}
Y:=\left\{[A,B]\in \P(V)  \mid AB^T = B^T A =
\lambda^2 \Id_m \textrm{ for some } \lambda \in \C\right\}.
\end{equation} 
The square at $\lambda$ seems to be irrelevant here, but it slightly simplifies the notation in the 
proofs  of Theorem~\ref{theorem_classify_invertible}(b)
and Proposition~\ref{proposition_orbits_of_G}(ii).
Although it is not essential for the content of this chapter, 
we note that $Y$ is $F$-cointegrable.

\label{notation_Xinv_Xdeg}

Further we define two types of subvarieties of $Y$:
\begin{gather*}
\Xinv(m):=\overline{
\bigg\{\left[g,\left(g^{-1}\right)^T\right] \in \P(V) \mid \det g = 1 \bigg\}},\\
\Xdeg(m,k):=\Big\{[A,B] \in \P(V) \mid  AB^T=B^TA=0, \ \rk A
  \le k, \ \rk B \le m-k  \Big\},
\end{gather*}
where $k \in {0,1,\ldots m}$.
The varieties $\Xdeg(m,k)$ have been also studied by \cite{strickland}
and \cite{mehta_trivedi}. $\Xinv(m)$ (especially $\Xinv(3)$) 
is the main object of this chapter.

\subsubsection{Automorphisms $\psi_{\mu}$} \label{notation_psi}

For any $\mu \in \C^*$ we let $\psi_{\mu}$ be the following linear
automorphism of $V$:
$$
\psi_{\mu}\big((A,B)\big) := (\mu A, \mu^{-1} B).
$$

Also the induced automorphism of $\P(V)$ will be written in the same
way:
$$
\psi_{\mu}\big([A,B]\big) := [\mu A, \mu^{-1} B].
$$

\medskip

%The notation introduced so far is sufficient to state the results of this
%chapter (see \S\ref{results}), but to prove them we need a few more
%notions. 

\label{notation_middle_end}

\subsubsection{Groups $G$ and $\widetilde{G}$, Lie algebra $\gotg$ and
   their representation}
\label{notation_G}

\label{notation_middle_begin}

We set $\widetilde{G}:=\Gl_m \times \Gl_m$ and let it act on $V$ by:
\begin{gather*}
(g,h) \in \widetilde{G}, \ g,h \in \Gl_m, \ (A,B) \in V
\\
(g,h) \cdot (A,B) := (g^T A h, g^{-1} B (h^{-1})^T).
\end{gather*}
This action preserves the symplectic form $\omega$.

We will mostly consider the restricted action of 
$G:=\Sl_m \times \Sl_m < \widetilde{G}$.

We also set $\gotg:= \gotsl_m \times \gotsl_m$ to be the Lie algebra of $G$ and we have the
tangent action of $\gotg$ on $V$:
$$
(g,h) \cdot (A,B) = (g^T A + A h, - g B - B h^T).
$$
Though we denote the action of the groups $G$, $\widetilde{G}$ 
and the Lie algebra $\gotg$ by the same $\cdot$ 
we hope it will not lead to any confusion. 
Also the induced action of $G$ and $\widetilde{G}$ on $\P(V)$ will
be denoted by $\cdot$.

\subsubsection{Orbits $\Inv^m$ and $\cDeg^m_{k,l}$}
\label{notation_Inv_Deg}

We define the following sets:
\begin{align*}
\Inv^m:=&
\bigg\{\left[g,\left(g^{-1}\right)^T\right] \in \P(V) \mid \det g = 1 \bigg\},
\\
\cDeg^m_{k,l}:=& \Big\{[A,B] \in \P(V) \mid  AB^T=B^TA=0, \ \rk A = k, \
\rk B = l \Big\},
\end{align*}
so that $\Xinv(m)=\overline{\Inv^m}$ and $\Xdeg(m,k)= \overline{\cDeg^m_{k,m-k}}$. 

Clearly, if $k+l>m$, then $\cDeg^m_{k,l}$ is empty, so whenever
we are considering
$\cDeg^m_{k,l}$ we will assume $k+l \le m$.

\subsubsection{Elementary matrices $E_{i j}$ and points $p_1$ and $p_2$}
\label{notation_Eij}

Let $E_{i j}$ be the elementary $m \times m$ matrix with
unit in the $i^{\textrm{th}}$ row and the $j^{\textrm{th}}$ column and zeroes
elsewhere.

\label{notation_p1_p2}

We distinguish two points 
$p_1 \in \cDeg^m_{1,0}$ and
$p_2 \in \cDeg^m_{0,1}$:
$$
p_1:=
\left[ E_{mm}, 0 \right]
\
\textrm{ and }
\
p_2:=
\left[ 0, E_{mm} \right]
$$

These points will be usually chosen as nice representatives of the
closed orbits
$\cDeg^m_{1,0}$ and $\cDeg^m_{0,1}$.

\subsubsection{Submatrices - extracting rows and columns}\label{notation_submatrices}

Assume $A$ is an $m\times m$ matrix and $I,J$ are two sets of indices of
cardinality $k$ and $l$ respectively:
\begin{align*}
I:=&\set{i_1,i_2, \ldots, i_k \mid 1 \le i_1 <  i_2 < \ldots < i_k \le m},\\
J:=&\set{j_1,j_2, \ldots, j_l \mid 1 \le j_1 < j_2 < \ldots < j_l \le m}.
\end{align*}
Then we denote by $A_{I,J}$ the $(m-k) \times (m-l)$ submatrix of $A$ obtained by 
removing rows of indices $I$ and columns of indices $J$.
Also for a set of indices $I$ we denote by $I'$ the set of $m-k$ indices complementary to $I$.

\smallskip

We will also use a simplified version of the above notation when we
remove only a single column and single row:
$A_{i j}$ denotes the $(m-1) \times (m-1)$ submatrix of $A$ obtained by removing $i$-th row
and $j$-th column, i.e. $A_{i j} = A_{\{i\},\{j\}}$

\smallskip

Also in the simplest situation where we remove only the last row and
the last column, we write $A_m$, so that $A_m = A_{mm} = A_{\{m\},\{m\}}$.

\label{notation_ends}

\section{Main results} \label{results}

In this chapter we give a classification%
\footnote{This problem was suggested by Sung Ho Wang.} 
of Legendrian subvarieties
in $\P(V)$ that are contained in $Y$.

\begin{theo}
\label{classification_theorem}
Let projective space $\P(V)$, varieties $Y$, $\Xinv(m)$,
 $\Xdeg(m,k)$ and automorphisms $\psi_{\mu}$ be defined as in
 \S\ref{notation_begins}.
Assume $X\subset\P(V)$ is an irreducible subvariety.
Then $X$ is Legendrian and contained in $Y$ if and only if 
$X$ is one of the following varieties:
\begin{itemize}
\item[1.]
$X = \psi_{\mu}(\Xinv(m))$ for some $\mu \in \C^*$ or 
\item[2.]
$X=\Xdeg(m,k)$ for some $k \in \{0,1,\ldots m\}$.
\end{itemize}
\end{theo}

The idea of the proof of Theorem~\ref{classification_theorem} is based on
the observation that every Legendrian subvariety that is contained in $Y$
must be invariant under the action of the group $G$. This is
explained in \S\ref{section_action}. A proof of the theorem is
presented in \S\ref{section_classification}.

Also we analyse which of the above varieties appearing in 1.~and 2.~are smooth:

\begin{theo}
    \label{theorem_classify_invertible}
  With the definition of $\Xinv(m)$ as in \S\ref{notation_Xinv_Xdeg},
  the family $\Xinv(m)$ contains the following varieties:
  \begin{itemize}
    \item[(a)] 
      $\Xinv(2)$ is a linear subspace.
    \item[(b)] 
      $\Xinv(3)$ is smooth, its  Picard group is generated by a hyperplane section.
      Moreover $\Xinv(3)$ is a compactification of $\Sl_3$ and it is
      isomorphic to a hyperplane section of Grassmannian $Gr(3,6)$. 
      The connected component of $\Aut(\Xinv(3))$ is equal to $G=\Sl_3 \times \Sl_3$ 
      and $\Xinv(3)$ is not a homogeneous space.
    \item[(c)]
      $\Xinv(4)$ is the 15 dimensional spinor variety $\mathbb{S}_6$.
    \item[(d)] 
      For $m\ge 5$, the variety $\Xinv(m)$ is singular.
  \end{itemize}
\end{theo}

A proof of the theorem is explained in \S\ref{section_invertible}.

Variety $\Xinv(3)$ is an original example of \cite{jabu_sl}.
Also it is the first described example of a 
smooth non-homogeneous Legendrian variety of dimension
bigger than $2$ (see \S\ref{section_intro_legendrian}).
This example is very close to a homogeneous one, namely it is isomorphic to a hyperplane section of $Gr(3,6)$, a well known subadjoint variety.
So a natural question arises, 
whether  general hyperplane sections of other Legendrian varieties admit a Legendrian embedding. 
The answer is yes and we explain it (as well as many conclusions from this surprisingly simple observation)
in Chapter~\ref{chapter_hyperplane}.

\begin{theo}
\label{theorem_smooth_degenerate}
With the definition of $\Xdeg(m)$ as in \S\ref{notation_Xinv_Xdeg},
variety $\Xdeg(m,k)$ is smooth if and only if $k=0$ , $k=m$ or $(m,k)=(2,1)$. 
In the first two cases, $\Xdeg(m,0)$ and $\Xdeg(m,m)$ are linear spaces, 
while $\Xdeg(2,1) \simeq \P^1 \times \P^1 \times \P^1 \subset \P^7$.
\end{theo}

A proof of the theorem is presented in \S\ref{section_degenerate}.

\medskip

\subsection{Generalisation: Representation theory and further examples}
    \label{section_other_examples}

The interpretation of 
Theorem~\ref{theorem_classify_invertible} (b) and (c) 
can be following:
We take the exceptional Legendrian variety $Gr(3,6)$, slice it with a linear section 
and we get a description, which generalised to matrices of bigger size, gives 
the bigger exceptional Legendrian variety $\bS_6$.
A similar connection can be established between other exceptional Legendrian varieties
(see \S\ref{section_intro_legendrian}).

\medskip

For instance, assume that $V^{sym}$ is a vector space of dimension $2 \binom{m+1}{2}$,
which we interpret as the space of pairs of $m \times m$ symmetric matrices $A,B$.
Now in $\P(V^{sym})$ consider the subvariety $\Xinv^{sym}(m)$, which is the closure of the following set:
$$
\{[A,A^{-1}]\in \P(V^{sym}) | A=A^T \textrm{ and } \det A=1\}.
$$
\begin{theo}\label{theorem_classify_sym}
All the varieties $\Xinv^{sym}(m)$ are Legendrian and we have:
\begin{itemize}
\item[(a)] $\Xinv^{sym}(2)$ is a linear subspace.
\item[(b)] $\Xinv^{sym}(3)$ is smooth and it is isomorphic to a
	   hyperplane section of Lagrangian Grassmannian $Gr_L(3,6)$.
\item[(c)] $\Xinv^{sym}(4)$ is smooth and it is Grassmannian variety $Gr(3,6)$.
\item[(d)] For $m\ge 5$, the variety $\Xinv^{sym}(m)$ is singular.
\end{itemize}
\end{theo}

The proof is exactly as the proof of Theorem~\ref{theorem_classify_invertible}.

\medskip

Similarly, we can take $V^{skew}$ to be a vector space of dimension $2 \binom{2m}{2}$,
which we interpret as the space of pairs of $2m \times 2m$ skew-symmetric matrices $A,B$.
Now in $\P(V^{skew})$ consider subvariety $\Xinv^{skew}(m)$, which is the closure of the following set:
$$
\{[A, - A^{-1}]\in \P(V^{skew}) | A= - A^T \textrm{ and } \pf A=1\}.
$$
\begin{theo}\label{theorem_classify_skew}
All the varieties $\Xinv^{skew}(m)$ are Legendrian and we have:
\begin{itemize}
\item[(a)] $\Xinv^{skew}(2)$ is a linear subspace.
\item[(b)] $\Xinv^{skew}(3)$ is smooth and it is isomorphic to a hyperplane section of the spinor variety $\bS_6$.
\item[(c)] $\Xinv^{skew}(4)$ is smooth and it is the 27 dimensional $E_7$ variety.
\item[(d)] For $m\ge 5$, the variety $\Xinv^{skew}(m)$ is singular.
\end{itemize}
\end{theo}

Here the only difference is that we replace the determinants by the
Pfaffians of the appropriate submatrices
and also for the previous cases we will be picking some diagonal matrices
as nice representatives.
Since there is no non-zero skew-symmetric diagonal matrix,
we must modify our calculations a little bit, 
but there is essentially no difference in the technique.

\medskip

Prior to \cite{jabu_sl} neither $\Xinv^{sym}(3)$ nor $\Xinv^{skew}(3)$ 
have been identified as smooth Legendrian subvarieties.

\medskip

Therefore we have established a connection 
between the subadjoint varieties of the 4 exceptional
groups $F_4$, $E_6$, $E_7$ and $E_8$.
A similar connection was obtained by \cite{landsbergmanivel02}.

We note that 
$m\times m$ symmetric matrices,
$m \times m $ matrices and
$2m \times 2m$ skew-symmetric matrices naturally correspond to 
$m \times m$ Hermitian matrices with coefficients in $\F \otimes_{\R} \C$, 
where $\F$ is the field of, respectively, 
real numbers $\R$, 
complex numbers $\C$ and
quaternions $\H$  (see \cite{landsbergmanivel01} and references therein).
An algebraic relation 
(analogous to Parts~(c) of 
Theorems~\ref{theorem_classify_invertible}, \ref{theorem_classify_sym} and \ref{theorem_classify_skew})
 between Lie algebras of types $E_6$, $E_7$ and $E_8$
and $4\times 4$ Hermitian matrices with coefficients in $\F \otimes_{\R} \C$ 
is described in \cite{brylinski_konstant}.

%Also it can be interesting to understand why this construction fails for the last step.
%For $Gr_L (3,6)$, $Gr(3,6)$, $\bS_6$, we just took a rank 4 Jordan algebra (corresponding respectively 
%to reals, complex numbers and quaternions) 
%and related to it the appropriate variety.
%So what happens if we take a $4\times 4$ octonionic hermitian matrices?
%One obstruction for defining the variety is that the determinants of degree 3 are not well defined
%(because octonions are not associative),
%so it is hard to say what $A^{-1}$ should be. 
%But note that to define $\Xinv^{sym}(4)$, $\Xinv(4)$ and  $\Xinv^{skew}(4)$, 
%we only use quadratic equations which can be expressed in a form of some matrix multiplication
%or as some combination of some $2 \times 2$ minors (or $4\times 4$ Pfaffians).
%So the only obstruction now is the noncommutativity 
%and somehow  it is not an obstruction for the quaternion case.

\section{$G$-action and its orbits}\label{section_action}

Recall the definition of $Y$ in \S\ref{notation_Y}.

The following polynomials are in 
the homogeneous ideal of $Y$ 
(the indices $i,j$ below run through $\{1,\ldots, m\}$, $k$ is a summation index):
\begin{subequations}\label{equations_of_Y}
\begin{align}
&\sum_{k=1}^m a_{i k} b_{i k}  - \sum_{k=1}^m a_{1k} b_{1k},
\label{equations_of_Y1}\\
&\sum_{k=1}^m a_{i k} b_{j k}  \ \textrm{ for $i \ne j$},
\label{equations_of_Y2}\\
&\sum_{k=1}^m a_{k i} b_{k i}  - \sum_{k=1}^m a_{k1} b_{k1},
\label{equations_of_Y3}\\
&\sum_{k=1}^m a_{k i} b_{k j}  \ \textrm{ for $i \ne j$}.
\label{equations_of_Y4}
\end{align}
\end{subequations}

These equations simply come from eliminating $\lambda$ from the
defining equation of $Y$ --- see Equation~\eqref{defining_equations_of_Y}.

For the statement and proof of the following proposition, recall our
notation of 
\S\ref{notation_Eij}.
\begin{prop}
\label{proposition_action_of_G}
Let $X\subset \P(V)$ be a Legendrian subvariety.
If $X$ is contained in $Y$, 
then $X$ is preserved by the induced action of $G$ on $\P(V)$.
\end{prop}
\begin{prf}
Let $\widetilde{\I}(X)_2$ be as in the Theorem~\ref{theorem_ideal_and_group} 
and define  $\widetilde{\I}(Y)_2$ analogously.
Clearly $\widetilde{\I}(Y)_2 \subset \widetilde{\I}(X)_2$.
Also let $\rho$ be the map described in Theorem~\ref{theorem_ideal_and_group}. 
By Theorem~\ref{theorem_ideal_and_group} it is enough to show that 
$\gotg \subset \rho\left(\widetilde{\I}(Y)_2\right)$ or that the images of the 
quadrics \eqref{equations_of_Y1}--\eqref{equations_of_Y4} under $\rho$ generate $\gotg$.

We write out the details of the proof only for $m=2$. 
There is no difference between this case and the general
one, except for the complexity of notation%
%(which should be 3-dimensional at this point to be clear enough).
.

Let us take the quadric
$$
q_{i j}:=\sum_{k=1}^m a_{i k} b_{j k} = a_{i1} b_{j1} + a_{i2} b_{j2} 
$$
for any $i,j \in \{1,\ldots, m\} = \{1,2\}$.
Also let $Q_{i j}$ be the $2m^2 \times 2m^2$ symmetric matrix
corresponding to $q_{i j}$. For instance:
$$
Q_{12} =
\left[ 
\begin{array}{cccccccc}
0&0 &0&0 &0&0 &\half & 0 \\
0&0 &0&0 &0&0 &0 & \half \\
0&0 &0&0 &0&0 &0 & 0 \\
0&0 &0&0 &0&0 &0 & 0 \\
0&0 &0&0 &0&0 &0 & 0 \\
0&0 &0&0 &0&0 &0 & 0 \\
\half&0 &0&0 &0&0 &0 & 0 \\
0&\half &0&0 &0&0 &0 & 0
\end{array}
\right] .
$$

Choose an arbitrary $(A,B) \in V$ and at the
moment we want to think of it as of a single vertical $2 m^2$-vector:
$(A,B) = [a_{11}, a_{12},a_{21}, a_{22},b_{11}, b_{12},b_{21},
b_{22}]^T$, so that the following multiplication makes sense:
\begin{multline*}
\rho(q_{12}) = 
2 J \cdot Q_{12} \cdot (A,B) = 
\\
= \left[ 
\begin{array}{cccccccc}
0&0 &0&0 &1&0 &0 & 0 \\
0&0 &0&0 &0&1 &0 & 0 \\
0&0 &0&0 &0&0 &1 & 0 \\
0&0 &0&0 &0&0 &0 & 1 \\
-1&0 &0&0 &0&0 &0 & 0 \\
0&-1 &0&0 &0&0 &0 & 0 \\
0&0 &-1&0 &0&0 &0 & 0 \\
0&0 &0&-1 &0&0 &0 & 0
\end{array}
\right]
\left[ 
\begin{array}{cccccccc}
0&0 &0&0 &0&0 &1 & 0 \\
0&0 &0&0 &0&0 &0 & 1 \\
0&0 &0&0 &0&0 &0 & 0 \\
0&0 &0&0 &0&0 &0 & 0 \\
0&0 &0&0 &0&0 &0 & 0 \\
0&0 &0&0 &0&0 &0 & 0 \\
1&0 &0&0 &0&0 &0 & 0 \\
0&1 &0&0 &0&0 &0 & 0
\end{array}
\right]
\left[ 
\begin{array}{c}
a_{11}\\
a_{12}\\
a_{21}\\
a_{22}\\
b_{11}\\
b_{12}\\
b_{21}\\
b_{22}
\end{array}
\right]
=
\\
=
\left[ 
\begin{array}{cccccccc}
0&0 &0&0 &0&0 &0 & 0 \\
0&0 &0&0 &0&0 &0 & 0 \\
1&0 &0&0 &0&0 &0 & 0 \\
0&1 &0&0 &0&0 &0 & 0 \\
0&0 &0&0 &0&0 &-1 & 0 \\
0&0 &0&0 &0&0 &0 & -1 \\
0&0 &0&0 &0&0 &0 & 0 \\
0&0 &0&0 &0&0 &0 & 0 
\end{array}
\right]
\left[ 
\begin{array}{c}
a_{11}\\
a_{12}\\
a_{21}\\
a_{22}\\
b_{11}\\
b_{12}\\
b_{21}\\
b_{22}
\end{array}
\right]
=
\\
=
\left[ 
\begin{array}{c}
0\\
0\\
a_{11}\\
a_{12}\\
-b_{21}\\
-b_{22}\\
0\\
0
\end{array}
\right]
\stackrel{\textrm{back to the matrix notation} }{=}
\left(
\left[ 
\begin{array}{cc}
0& 0 \\
a_{11} & a_{12}
\end{array}
\right] , \
\left[ 
\begin{array}{cc}
-b_{21} & -b_{22}\\
0& 0
\end{array}
\right]
\right)
= 
\\
=
\left(
\left[ 
\begin{array}{cc}
0 &1 \\
0 &0
\end{array}
\right]^T
\left[ 
\begin{array}{cc}
a_{11} & a_{12}\\
a_{21} & a_{22}\\
\end{array}
\right] , \ 
- \left[ 
\begin{array}{cc}
0& 1 \\
0 &0
\end{array}
\right]
\left[ 
\begin{array}{cc}
b_{11} & b_{12}\\
b_{21} & b_{22}
\end{array}
\right]
\right)
 \ = \ (E_{12}^T A, \  -E_{12} B).
\end{multline*}

Similar calculations show that:

%\begin{equation} \label{sl_elementary_action}
$$
2 J \cdot Q_{i j} \cdot (A,B)  \ = \ (E_{i j}^T A, \  -E_{i j} B).
$$
%\end{equation}

Next in the ideal of $Y$ we have the following quadrics: $q_{i j}$ for
$i \ne j$ 
(see Equation~\eqref{equations_of_Y2})
 and $q_{ii}-q_{11}$ 
(see Equation~\eqref{equations_of_Y1}).
By taking images under $\rho$ of the linear combinations of those quadrics we can get an arbitrary
traceless matrix $g\in \gotsl_m$ acting on $V$ in the following way:
$$
g \cdot (A,B) = (g^T A, -g B).
$$
Exponentiate this action of $\gotsl_m$ to get the action of $\Sl_m$:
$$
g \cdot (A,B) = (g^T A, g^{-1} B).
$$

This proves that the action of subgroup $\Sl_m \times 0 < G = \Sl_m \times \Sl_m$
preserves $X$ as claimed in the proposition. 
The action of the other component $0 \times \Sl_m$ is calculated in the same way, but
using quadrics \eqref{equations_of_Y3}--\eqref{equations_of_Y4}.
\end{prf}

\subsection{Invariant subsets}\label{section_orbits}

Here we want to decompose $Y$ into a union of some $G$-invariant
subsets, most of which are orbits.

\begin{prop}
\label{proposition_orbits_of_G}
\hfill
\begin{itemize}
\item[(i)]
The sets $\Inv^m$, $\psi_{\mu}(\Inv^m)$ and $\cDeg^m_{k,l}$ 
are $G$-invariant and they are all contained in $Y$.
\item[(ii)]
$Y$ is equal to the union of all $\psi_{\mu}(\Inv^m)$ (for
  $\mu \in \C^*$) and all $\cDeg^m_{k,l}$ (for integers $k,l\ge
  0$, $k+l\le m$).
\item[(iii)]
Every $\psi_{\mu}(\Inv^m)$ is an orbit of the action of $G$. 
If $m$ is odd, then $\Inv^m$ is isomorphic (as algebraic variety) to
$\Sl_m$. Otherwise if $m$ is even, then $\Inv^m$ is isomorphic to 
$(\Sl_m / \Z_2) $. In both cases 
\[
\dim \psi_{\mu}(\Inv^m) =  \dim
\Inv^m = m^2 -1.
\] 
\end{itemize}
\end{prop}

\begin{prf}
The proof of part (i) is an explicit verification from the definitions
in \S\ref{section_basics2}.

\medskip

To prove part (ii), assume $[A,B]$ is a point of $Y$, so $A B^T=B^TA=\lambda^2\Id_m$. 
First assume that the ranks of both matrices are maximal: 
$$
\rk A = \rk B =m. 
$$
Then $\lambda$ must be non-zero and $B=\lambda^2 (A^{-1})^T$.
Let $d:=(\det A)^{-\frac{1}{m}} $ so that 
$$
\det (d A)=1
$$
and let 
$\mu:= \frac{1}{d \lambda} $.
Then we have:
\begin{align*}
[A,B] =&
 \left[A,\lambda^2 \left(A^{-1}\right)^T\right] =
 \left[ \frac{d A}{d \lambda }, d\lambda \left((d A)^{-1}\right)^T\right] =\\
=& \left[ \mu (d A), \mu^{-1} \left((d A)^{-1}\right)^T\right]= 
 \psi_{\mu} \left(\left[ (d A), \left((d A)^{-1}\right)^T\right]\right).
\end{align*}
Therefore $[A,B] \in \psi_{\mu}(\Inv^m)$.

Next, if either of the ranks is not maximal:
$$
\rk A < m \textrm{ or } \rk B < m,
$$
then by Equation~\eqref{defining_equations_of_Y} we must have $AB^T = B^TA =
0$. 
So $[A,B] \in \cDeg^m_{k,l}$ for $k=\rk A$ and $l=\rk B$. 

\medskip

Now we prove (iii).
The action of $G$ commutes with $\psi_{\mu}$:
$$
(g,h)\cdot \psi_{\mu}\big([A,B]\big) = \psi_{\mu}\big((g,h) \cdot [A,B]\big).
$$
So to prove $\psi_{\mu}(\Inv^m)$ is an orbit it is enough to prove
that $\Inv^m$ is an orbit, which follows from the definitions of the
action and $\Inv^m$.

We have the following epimorphic map:
$$
\begin{array}{rcl}
\Sl_m &\lra & \Inv^m\\
g &\longmapsto& [g,(g^{-1})^T].
\end{array}
$$
If $[g_1,(g_1^{-1})^T] = [g_2,(g_2^{-1})^T]$, then we must have
$g_1 = \alpha g_2$ and $g_1 = \alpha^{-1} g_2$ for some 
$\alpha \in \C^*$. 
Hence $\alpha^2=1$ and $g_1 = \pm g_2$. If $m$ is odd and 
$g_1\in \Sl_m$, then $-g_1 \notin \Sl_m$ so $g_1 = g_2$. So $\Inv^m$ is 
either isomorphic to $\Sl_m$ or to $\Sl_m /\Z_2$ as stated.
\end{prf}

From Proposition~\ref{proposition_orbits_of_G}(ii) we conclude that 
$\Xinv(m)$ is an equivariant compactification of $\Sl_m$ (if $m$ is odd) 
or $\Sl_m / \Z_2$ (if $m$ is even).
See \cite{timashev} and references therein 
for the theory of equivariant compactifications. 
In the setup of \cite[\S8]{timashev}, 
this is the compactification corresponding to the representation $W \oplus W^*$,
where $W$ is the standard representation of $\Sl_m$. 
Therefore some properties of $\Xinv(m)$ could also be read from the general description 
of group compactifications.

\begin{prop}
\label{degenerate_orbits_of_G}
\hfill
\begin{itemize}
\item[(i)]
The dimension of $\cDeg^m_{k,l}$ is $(k+l)(2m-k-l)-1$.
In particular, if $k+l=m$, then the dimension is equal to $m^2-1$.
\item[(ii)]
$\cDeg^m_{k,l}$ is an orbit of the action of $G$, unless $m$ is even
  and $k=l= \half m$. 
\item[(iii)]
If $m\ge 3$, then there are exactly two closed orbits of the action of
$G$: $\cDeg^m_{1,0}$ and $\cDeg^m_{0,1}$.
\end{itemize}
\end{prop}

\begin{prf}
Part (i) follows from \cite[Prop.~2.10]{strickland}.

\medskip

For part (ii) 
let $[A,B] \in \cDeg^m_{k,l}$ be any point.
By Gaussian elimination and elementary linear algebra,
we can prove that there exists 
$(g,h) \in G$ such that $[A',B']:=(g,h) \cdot [A,B]$ is a pair of
diagonal matrices. Moreover, if $k+l < m$, then we can choose $g$ and
$h$ such that:
$$
A':=\diag(\underbrace{1,\ldots, 1}_{k},\underbrace{0,\ldots, 0}_{l},
\underbrace{0,\ldots, 0}_{m-k-l}),
$$\nopagebreak
$$
B':=\diag(\underbrace{0,\ldots, 0}_{k},\underbrace{1,\ldots, 1}_{l},
\underbrace{0,\ldots, 0}_{m-k-l}).
$$
Hence $\cDeg^m_{k,l} = G\cdot [A',B']$  and this finishes the proof in the case
$k+l < m$. 

So assume $k+l=m$. Then we can choose $(g,h)$ such that:
$$
A':=\diag(\underbrace{1,\ldots, 1}_{k},\underbrace{0,\ldots, 0}_{l}),
$$
$$
B':=\diag(\underbrace{0,\ldots, 0}_{k},\underbrace{d,\ldots, d}_{l}),
$$
for some $d \in \C^*$. 
If $k\ne l$, then set $e:=d^{\frac{1}{l-k}}$ and let
$$
g':=\diag(\underbrace{e^l,\ldots, e^l}_{k},\underbrace{e^{-k},\ldots, e^{-k}}_{l}).
$$
Clearly $\det (g')=1$ and:
$$
(g', \Id_m) \cdot [A',B'] = 
\left[
\diag(\underbrace{e^l,\ldots, e^l}_{k},\underbrace{0,\ldots,0}_{l}),
\diag(\underbrace{0,\ldots, 0}_{k},\underbrace{d e^k,\ldots,d e^k}_{l})
\right]
$$
where
$$
d e^k = d^{1+\frac{k}{l-k}} = d^{\frac{l}{l-k}} = e^l.
$$
So rescaling we get:
$$
(g', \Id_m) \cdot [A',B'] = 
\left[
\diag(\underbrace{1,\ldots, 1}_{k},\underbrace{0,\ldots,0}_{l}),
\diag(\underbrace{0,\ldots, 0}_{k},\underbrace{1,\ldots,1}_{l})
\right]
$$
and this finishes the proof of (ii).

\medskip

For part (iii),
denote by $W_1$ (respectively, $W_2$) the standard representation of 
the first (respectively, the second) component of $G= \Sl_m \times \Sl_m$ 
with the trivial action of the other component.
Then $V$ as a representation of $G$ is isomorphic to $(W_1 \otimes W_2) \oplus (W_1^*\otimes W_2^*)$. 
For $m\ge 3$ the representation $W_i$ is not isomorphic to $W_i^*$ and
therefore $V$ is a union of two irreducible non-isomorphic
representations. 
It is a standard fact from the representation theory, 
that the closed orbits of the action of a semisimple group on a projectivised representation
are in one to one correspondence with the irreducible subrepresentations.
So there are exactly two closed orbits of the action of $G$ on $\P(V)$. 
These orbits are simply $\cDeg^m_{1,0}$ and $\cDeg^m_{0,1}$.
\end{prf}

\subsection{Action of $\widetilde{G}$}

The action of $\widetilde{G}$ extends the action of $G$, but it does
not preserve $\Xinv(m)$. 
So we will only consider the action of $\widetilde{G}$ when speaking of
$\Xdeg(m,k)$.

We have properties analogous to 
Proposition~\ref{degenerate_orbits_of_G}(ii) and  (iii) but with no exceptional cases:

\begin{prop}
\label{properties_of_Gtilde}
\hfill
\begin{itemize}
\item[(i)]
Every $\cDeg^m_{k,l}$ is an orbit of the action of $\widetilde{G}$.
\item[(ii)]
For every  $m$ there are exactly two closed orbits of the action of
$\widetilde{G}$: $\cDeg^m_{1,0}$ and $\cDeg^m_{0,1}$.
\end{itemize}
\end{prop}

\begin{prf}
This is exactly as the proof of 
Proposition~\ref{degenerate_orbits_of_G}(ii) and (iii).
\end{prf}

\section{Legendrian varieties in $Y$}\label{main_part}

In this section we prove the main results of the chapter.

\subsection{Classification}\label{section_classification}

We start by proving the Theorem~\ref{classification_theorem}. 

\begin{prf}
First assume $X$ is Legendrian and contained in $Y$.
If $X$ contains a point $[A,B]$ where both
 $A$ and $B$ are invertible,  
then by Proposition~\ref{proposition_action_of_G}
 it must contain the orbit of $[A,B]$, which by 
Proposition~\ref{proposition_orbits_of_G}(ii) and (iii) is equal to 
$\psi_{\mu}(\Inv^m)$ for some $\mu \in \C^*$. But the dimension of $X$ is $m^2-1$ which is
 exactly the dimension of $\psi_{\mu}(\Inv^m)$ 
(see Proposition~\ref{proposition_orbits_of_G}(iii)), 
so
$$
X =\overline{\psi_{\mu}(\Inv^m)} = \psi_{\mu}(\Xinv(m)).
$$

\medskip

On the other hand, 
if $X$ does not contain any point $[A,B]$ where both $A$ and $B$ are invertible,
then in fact $X$ is contained in the locus $Y_0:=\{[A,B]:AB^T=B^T A=0\}$. 
This locus is just the union of all $\cDeg^m_{k,l}$ and its irreducible components are
the closures of $\cDeg^m_{k,m-k}$, which are exactly $\Xdeg(m,k)$.
 So in particular every irreducible component has
 dimension $m^2-1$ (see Proposition~\ref{degenerate_orbits_of_G}(i)) 
and hence $X$ must be one of these components.

\medskip

Therefore it remains to show that all these varieties are Legendrian. 

The fact that $\Xdeg(m,k)$ is a Legendrian variety follows from
\cite[pp524--525]{strickland}. Strickland proves there that the affine
cone over $\Xdeg(m,k)$ (or $W(k, m-k)$ in the notation of
\cite{strickland}) is the closure of a conormal bundle. Conormal bundles
are classical examples of Lagrangian varieties 
(see Example~\ref{example_conormal_Lagrangian}).

\medskip

Since $\psi_{\mu}$ preserves the symplectic form $\omega$, it is
enough to prove that $\Xinv(m)$ is Legendrian.

The group $G$ acts symplectically on $V$ and the action has an open
orbit on $\Xinv(m)$ --- see Proposition~\ref{proposition_orbits_of_G} (iii).
Thus the tangent spaces to the affine cone over $\Xinv(m)$ are Lagrangian if
 and only if just one tangent space at a point of the open orbit is
 Lagrangian.

So we take $[A,B]:=[\Id_m, \Id_m]$.
Now the affine tangent space to $\Xinv(m)$ at $[\Id_m,\Id_m]$ 
is the linear subspace of $V$ spanned by $(\Id_m,\Id_m)$ 
and the image of the tangent action of the Lie algebra $\gotg$. 
We must prove that for every four traceless matrices
$g,h, g', h'$ we have:
\begin{subequations}
\begin{align}
\label{equation_for_Inv_Legendrian1}
\omega \bigl((g,h) \cdot (\Id_m,\Id_m), \ (g',h') \cdot
(\Id_m,\Id_m)\bigr) &= 0 \ \textrm{ and}\\
\label{equation_for_Inv_Legendrian2}
\omega \bigl((\Id_m,\Id_m),  \ (g,h)\cdot (\Id_m,\Id_m)\bigr)& = 0.
\end{align}
\end{subequations}

Equality~\eqref{equation_for_Inv_Legendrian1} is true without the
assumption on the trace of the matrices:
\begin{multline*}
\omega \big((g,h) \cdot (\Id_m,\Id_m),  \ (g',h') \cdot
(\Id_m,\Id_m)\big) 
\\
= \ \omega\Big( \left( g^T+h, \ -(g+h^T)\right), \quad
\left((g')^T+h', \ -(g' +(h')^T)\right)\Big)
\\
\stackrel{\textrm{by \eqref{properties_of_inner_product}}}{=}
 \   \tr \Big( -\left(g^T+h\right) \left((g')^T + h'\right)   \  +
 \  \left(g+h^T\right) \left( g'+(h')^T\right)   \Big)  \
= \ 0.
\end{multline*}

For Equality~\eqref{equation_for_Inv_Legendrian2} we calculate:
\begin{multline*}
\omega \big((\Id_m,\Id_m),  \ (g,h) \cdot (\Id_m,\Id_m)\big) 
\\
= \ \omega\Big((\Id_m,\Id_m) , \quad
\left(g^T+h, \ -(g +h^T)\right)\Big) 
\\
\stackrel{\text{by \eqref{properties_of_inner_product}}}{=}
- \tr (g^T+h)  - \tr (g+h^T) \ = \ 0.
\end{multline*}

Hence we have proved that the closure of $\Inv^m$ is Legendrian.
\end{prf}

\subsection{Degenerate matrices}\label{section_degenerate}

By \cite[Prop.~1.3]{strickland} the ideal of $\Xdeg(m,k)$ is generated by
the coefficients of $AB^T$, the coefficients of $B^TA$, 
the $(k+1)\times(k+1)$-minors of $A$ and the $(m-k+1)\times(m-k+1)$-minors of $B$.
In short, we will say that the equations of $\Xdeg(m,k)$ are given by:

\begin{equation}\label{equations_of_Xdeg}
AB^T=0, \ B^TA=0, \quad \rk(A) \le k, \ \rk(B) \le m-k.
\end{equation}

\begin{lem}\label{lemma_tangent_cone_to_Xdeg}
Assume $m\ge 2$ and $1\le k \le m-1$. Then:
\begin{itemize}
\item[(i)]
The tangent cone to $\Xdeg(m,k)$ at $p_1$ is a product of a linear space of dimension 
$(2m-2)$ and the affine cone over $\Xdeg(m-1,k-1)$.
\item[(i')]
The tangent cone to $\Xdeg(m,k)$ at $p_2$ is a product of a linear space of dimension 
$(2m-2)$ and the affine cone of $\Xdeg(m-1,k)$.
\item[(ii)]
$\Xdeg(m,k)$ is smooth at $p_1$ if and only if $k=1$.
\item[(ii')]
$\Xdeg(m,k)$ is smooth at $p_2$ if and only if $k=m-1$.
\end{itemize}
\end{lem}

\begin{prf}
We only prove (i) and (ii), while (i') and (ii') follow in the same way by exchanging $a_{i j}$ and $b_{i j}$.
Consider Equations~\eqref{equations_of_Xdeg} of $\Xdeg(m,k)$ 
restricted to the affine neighbourhood of $p_1$ obtained by substituting $a_{mm}=1$. 
Taking the lowest degree part of these equations we get some of the equations of the tangent cone at $p_1$
(recall our convention on the notation of submatrices  --- see 
\S\ref{notation_submatrices}):
$$
b_{i m}=b_{mi}=0, \ A_m B_m^T =0, \ B_m^T A_m = 0, 
$$
$$ 
\rk A_m \le k-1, \rk B_m \le m-k.
$$
These equations define the product of the linear subspace  $A_m=B_m=0, b_{i m}=b_{mi}=0$ 
and the affine cone over $\Xdeg(m-1,k-1)$ embedded in the set of those pairs of matrices,
 whose last row and column are zero: $a_{i m}=a_{mi}=0, b_{i m}=b_{mi}=0$. 
So the variety defined by those equations is irreducible and its dimension is equal to
$(m-1)^2 + 2m-2 = m^2-1 = \dim \Xdeg(m,k)$. 
Since this contains the tangent cone we are interested in
 and by \S\ref{properties_of_tangent_cone}(1),
they must coincide as claimed in (i).

\medskip

Next (ii) follows immediately, since for $k=1$ the equations above reduce to 
$$
b_{i m}=b_{mi}=0, \ \textrm{ and } A_m=0
$$
and hence the tangent cone is just the tangent space, so  $p_1$ is a smooth point of $\Xdeg(m,1)$.
Conversely,  if $k>1$, then $\Xdeg(m-1, k-1)$ is not a linear space, so by (i) 
the tangent cone is not a linear space either and $X$ is singular at $p_1$ --- see 
\S\ref{properties_of_tangent_cone}(3).
\end{prf}

Now we can prove Theorem~\ref{theorem_smooth_degenerate}:

\begin{prf}
It is obvious from the definition of $\Xdeg(m,k)$, 
that $\Xdeg(m,0)= \{A=0\}$ and $\Xdeg(m,m) = \{B=0\}$,
so these are indeed linear spaces.

Therefore assume $1\le k \le m-1$. 
But $\Xdeg(m,k)$ is
$\widetilde{G}$ invariant (see Proposition~\ref{properties_of_Gtilde}(i)) and so is its singular locus $S$. 
Hence $\Xdeg(m,k)$ is singular if and only if $S$ contains a closed orbit of $\widetilde{G}$.

So $\Xdeg(m,k)$ is smooth, 
if and only if it is smooth at both $p_1$ and $p_2$
(see Proposition~\ref{properties_of_Gtilde}(ii)), which 
by Lemma~\ref{lemma_tangent_cone_to_Xdeg}(ii) and (ii')
holds if and only if $k=1$ and $m=2$.

\smallskip
 
To finish the proof, 
it remains to verify what kind of variety is $\Xdeg(2,1)$.
Consider the following map:
\begin{align*}
\P^1 \times \P^1 \times \P^1 \lra &\P(V) \simeq \P^7
\\
[\mu_1, \mu_2],[\nu_1,\nu_2],[\xi_1,\xi_2] \longmapsto&
\left[
\xi_1
\left(
\begin{array}{cc}
\mu_1\nu_1 & \mu_1 \nu_2\\
\mu_2\nu_1 & \mu_2 \nu_2
\end{array}
\right),
\xi_2
\left(
\begin{array}{cc}
\mu_2\nu_2 & -\mu_2 \nu_1\\
-\mu_1\nu_2 & \mu_1 \nu_1
\end{array}
\right)
\right]
\end{align*}
Clearly this is the Segre embedding in appropriate coordinates.
The image of this embedding is contained in $\Xdeg(2,1)$ 
(see Equation~\eqref{equations_of_Xdeg}) 
and since dimension of $\Xdeg(2,1)$ 
is equal to the dimension of  $\P^1 \times \P^1 \times \P^1$ 
we conclude the above map gives an isomorphism of 
$\Xdeg(2,1)$ and  $\P^1 \times \P^1 \times \P^1$.
\end{prf}

\subsection{Invertible matrices}\label{section_invertible}

We wish to determine some of the equations of $\Xinv(m)$.
Clearly the equations of $Y$ (see Equation~\eqref{equations_of_Y}) 
are quadratic equations of $\Xinv(m)$.
To find other equations, we recall that 
$$
\Xinv(m):=\overline{
\bigg\{\left[g,\left(g^{-1}\right)^T\right] \in \P(V) \mid \det g = 1 \bigg\}}.
$$
However, for a matrix $g$ with determinant 1 we know that the entries of
$(g^{-1})^T$
consist of the appropriate minors (up to sign)
of $g$.
Therefore we get many inhomogeneous equations satisfied by every pair
$\left(g, (g^{-1})^T\right) \in V$
(recall our convention on the notation of submatrices  --- see 
\S\ref{notation_submatrices}):
$$
\det (A_{i j}) = (-1)^{i+j} b_{i j} \  \textrm{ and } \  a_{kl} = (-1)^{k+l} \det(B_{kl})  
$$

To make them homogeneous, multiply two such equations appropriately:
\begin{equation}\label{minors_1_of_X}
\det(A_{i j}) a_{kl} = (-1)^{i+j+k+l} b_{i j} \det(B_{kl}).
\end{equation}
These are degree $m$ equations, 
which are satisfied by the points of $\Xinv(m)$ 
and we state the following theorem:

\begin{theo}
\label{theorem_Xinv3_is_smooth}
Let $m=3$. Then the quadratic polynomials \eqref{equations_of_Y1}--\eqref{equations_of_Y4} and
the cubic polynomials \eqref{minors_1_of_X} generate a homogeneous ideal $\I$ in $\C[V]$ 
which defines $\Xinv(3)$ as a subscheme of $\P(V)$. Moreover $\Xinv(3)$ is smooth.
\end{theo}

\begin{prf}
It is enough to prove that the scheme $X$ defined by $\I$ is smooth, because the  reduced subscheme
of $X$ coincides with $\Xinv(3)$.

The scheme $X$ is $G$ invariant, 
hence as in the proof of Theorem~\ref{theorem_smooth_degenerate} and by
Proposition~\ref{degenerate_orbits_of_G}(iii) it is
 enough to verify smoothness 
at $p_1$ and $p_2$.
Since we have the additional symmetry here (exchanging $a_{i j}$'s with $b_{i j}$'s)
it is enough to verify the smoothness at $p_1$.

Now we calculate the tangent space to $X$ at $p_1$ by taking 
linear parts of the polynomials evaluated at $a_{33}=1$. 
From polynomials \eqref{equations_of_Y} we get that 
$$
b_{31}=b_{32}=b_{33}=b_{23}=b_{13}=0.
$$
Now from Equations~\eqref{minors_1_of_X} for $k=l=3$ and $i,j \ne 3$  we get the following evaluated equations:
$$
a_{i' j'} - a_{i'3}a_{3j'} = \pm b_{i j} B_{33} 
$$
(where $i'$ is either $1$ or $2$, which ever is different than $i$ and analogously for $j'$)
so the linear part is just $a_{i' j'}=0$. Hence by varying $i$ and $j$ we can get 
$$
a_{11}= a_{12}=a_{21}=a_{22}=0.
$$
Therefore the tangent space has codimension at least $9$, which is exactly the codimension of 
$\Xinv(3)$ --- see Proposition~\ref{proposition_orbits_of_G}(iii). 
Hence $X$ is smooth (in particular reduced) and $X=\Xinv(3)$. 
\end{prf}

\begin{rem} \label{remark_ideal_of_Xinv3_is_saturated}
  The ideal $\I$ from Theorem~\ref{theorem_Xinv3_is_smooth} could possibly be unsaturated
  (in other words, the affine subscheme in $V$ defined by this ideal could be not reduced at $0$).
  Computer calculations show that this is not the case and $\I$ is saturated. 
  We will only need that the saturated ideal of $\Xinv(3)$ has no other quadratic polynomials than $\I$ 
  and we will prove it later.
\end{rem}

To describe $\Xinv(m)$ for $m > 3$ we must find more equations. 

There is a more general version of the above  property of an inverse
of a matrix with determinant 1, which is less commonly known.

\begin{prop}\label{proposition_minors_2}
\hfill
\begin{itemize}
\item[(i)]
Assume $A$ is a $m \times m$ matrix of determinant $1$ and $I,J$ are
two sets of indices, both of cardinality $k$
(again recall our convention on indices and submatrices  --- see
\S\ref{notation_submatrices}). Denote by $B:=(A^{-1})^T$. 
Then the appropriate minors are equal (up to sign):
$$
\det A_{I,J} = (-1)^{\Sigma I + \Sigma J} \det B_{I', J'}.
$$
\item[(ii)]
A coordinate free way to express these equalities is following: 
Assume $W$ is a vector space of dimension $m$, 
$f$ is a linear automorphism of $W$
and $k \in \set{0,\ldots, m}$.
Let $\Wedge{k} f$ be the induced automorphism of $\Wedge{k} W$.
If $\Wedge{m} f = \Id_{\Wedge{m} W}$,  then:
$$
\Wedge{m-k} f  = \Wedge{k} \left( \Wedge{m-1} f\right) .
$$
\item[(iii)]
Consider the induced action of $G$ on the polynomials on $V$.
Then the vector space spanned by the set of equations of (i) for a fixed
	  $k$ is $G$ invariant.
\end{itemize}
\end{prop}

\begin{prf}
Part (ii) follows immediately from (i), since 
if $A$ is a matrix of $f$, then the terms of  the matrices of  the maps 
$\Wedge{m-k} f$ and $\Wedge{k} ( \Wedge{m-1} f)$
 are exactly the appropriate minors of $A$ and $B$.

Part (iii) follows easily from (ii).

As for (i), 
we only sketch the proof, leaving the details to the reader.
Firstly, reduce to the case when $I$ and $J$ are just $\{1, \ldots, k\}$ and the determinant of $A$
is possibly $\pm 1$ (which is where the sign shows up in the equality).
Secondly if both determinants $\det A_{I,J}$ and $\det B_{I', J'}$ are 
zero, then the equality is clearly satisfied. Otherwise assume for example $\det A_{I,J} \ne 0$.
Then performing the appropriate row and column operations we can change $A_{I,J}$ into a diagonal matrix, 
$A_{I', J}$ and $A_{I, J'}$ into the zero matrices 
and all these operations can be done without changing $B_{I',J'}$ 
nor $\det A_{I,J}$. 
Then the statement follows easily.
\end{prf}

In particular we get:

\begin{cor}
\label{corollary_equations_Xinv}
Assume $k$, $I$ and $J$ are as in Proposition~\ref{proposition_minors_2}(i).
\begin{itemize}
\item[(a)]
If $m$ is even and $k=\half m$, then the equation
$$
\det A_{I,J} = (-1)^{\Sigma I + \Sigma J} \det B_{I', J'}
$$
is homogeneous of degree $\half m$ and it is satisfied by points of $\Xinv(m)$.
\item[(b)]
If $0 \le k < \half m $ and $l= m - 2k$, then
$$
\left(\det A_{I,J}\right)^2 = \left(\det B_{I',J'}\right)^2 \cdot (a_{11}b_{11} + \ldots + a_{1m}b_{1m})^l
$$
is a homogeneous equation of degree $2(m-k)$ satisfied by points of $\Xinv(m)$.
\end{itemize}
\end{cor}

\begin{prf}
Clearly both equations are homogeneous.
If $\det A =1$ and $B=(A^{-1})^T$, then the following equations are satisfied: 
\begin{equation}\label{equation_for_cor1}
\det A_{I,J} = (-1)^{\Sigma I + \Sigma J} \det B_{I',J'},
\end{equation}
\begin{equation}\label{equation_for_cor2}
1 = (a_{11}b_{11} + \ldots  a_{1m}b_{1m})^l 
\end{equation}
(Equation~\eqref{equation_for_cor1} follows from Proposition~\ref{proposition_minors_2}(i)
and Equation~\eqref{equation_for_cor2} follows from $AB^T=\Id_m$).
Equation in (b) is just Equation~\eqref{equation_for_cor1} squared multiplied side-wise 
by Equation~\eqref{equation_for_cor2}.

So both equations in (a) and (b) are satisfied by every pair 
$\left( A, (A^{-1})^T\right)$ 
and by homogeneity also by
$\left(\lambda A, \lambda (A^{-1})^T \right)$. 
Hence (a) and (b) hold on an open dense subset of $\Xinv(m)$, so also on whole $\Xinv(m)$.
\end{prf}

We know enough equations of $\Xinv(m)$ to prove the Theorem~\ref{theorem_classify_invertible}:

\subsubsection{Case $m=2$ --- linear subspace}

\begin{prf}
To prove (a) just take the linear equations from  Corollary~\ref{corollary_equations_Xinv}(a) for $k=1$:
$$
a_{i j} = \pm b_{i' j'},
$$
where $\{i,i'\} = \{j,j'\} = \{1,2\}$.
\end{prf}

\subsubsection{Case $m=3$ --- hyperplane section of $Gr(3,6)$}

\begin{prf}
For (b), $\Xinv(3)$ is smooth by Theorem~\ref{theorem_Xinv3_is_smooth} 
and it is a compactification of $\Inv^3 \simeq \Sl_3$ by 
Proposition~\ref{proposition_orbits_of_G}(i) and (iii).

\paragraph{Picard group of $\Xinv(3)$.}

The complement of the open orbit 
$$
D:=\Xinv(3) \backslash \Inv^3
$$
must be a union of some orbits of $G$, each of them must have
dimension smaller than $\dim \Inv^3 =8$. 
So by Propositions~\ref{proposition_orbits_of_G}(ii), (iii), 
\ref{degenerate_orbits_of_G}(i) and (ii) 
the only candidates are 
$\cDeg^3_{1,1}$, $\cDeg^3_{0,1}$, $\cDeg^3_{1,0}$,  $\cDeg^3_{0,2}$ and $\cDeg^3_{2,0}$.
We claim that only $\cDeg^3_{1,1}$, $\cDeg^3_{0,1}$ and $\cDeg^3_{1,0}$ are contained in $\Xinv(3)$. 
To exclude  $\cDeg^3_{2,0}$ consider the point:
\[
\left[
\left(
\begin{array}{ccc}
1  &0  &0\\
0  &1  &0\\
0  &0  &0\\
\end{array}
\right),
0
\right].
\]
This point is in $\cDeg^3_{2,0}$ and does not satisfy Equation~\eqref{minors_1_of_X} for $i=j=1$ and $k=l=3$ 
and therefore it is not contained in $\Xinv(3)$. So $\cDeg^3_{2,0}$ is disjoint from $\Xinv(3)$.

Analogously, $\cDeg^3_{0,2}$ is disjoint from $\Xinv(3)$. 

It only remains to prove
that $\cDeg^3_{1,1} \subset \Xinv(3)$, since the other orbits are in
the closure of $\cDeg^3_{1,1}$. Take the curve in $\Xinv(3)$
parametrised by:
$$
\left[
\left(
\begin{array}{ccc}
t  &0  &0\\
0  &1  &0\\
0  &0  &t^{-1}\\
\end{array}
\right),
\left(
\begin{array}{ccc}
t^{-1} &0  &0\\
0      &1  &0\\
0      &0  &t\\
\end{array}
\right)
\right].
$$
For $t=0$ the curve meets $\cDeg^3_{1,1}$, and therefore  $\cDeg^3_{1,1}$ is contained in  $\Xinv(3)$ 
and that finishes the proof of the claim.

Since $\dim \cDeg^3_{1,1} = 7$  
(see Proposition~\ref{degenerate_orbits_of_G}(i)), $D$ is a prime divisor.
We have $\Pic(\Sl_3) = 0$ and by \cite[Prop.~II.6.5(c)]{hartshorne} the Picard group of 
 $\Xinv(3)$ is isomorphic to $\Z$ with the ample generator $[D]$.

Next we check that $D$ is linearly equivalent 
(as a divisor on $\Xinv(3)$) to a hyperplane section $H$ of $\Xinv(3)$. 
Since we already know that $\Pic(\Xinv(3))=\Z \cdot [D]$, we must have
$H \stackrel{lin}{\sim} k D$ for some positive integer $k$.
But there are lines contained in $\Xinv(3)$ (for example those
contained in $\cDeg^3_{1,0} \simeq \P^2 \times \P^2$)\footnote{
Actually, the reader could also easily find explicitly some lines (or
even planes) which intersect the open orbit 
and conclude that $\Xinv(3)$ is covered by lines.
}.
So let $L\subset \Xinv(3)$ be any line and we intersect:
$$
D\cdot L = \frac{1}{k} H \cdot L = \frac{1}{k}.
$$
But the result must be an integer, so $k=1$ as claimed.

\paragraph{Complete embedding.}

Since $D$ itself is definitely not a hyperplane section of $\Xinv(3)$, 
the conclusion is that the Legendrian embedding of $\Xinv(3)$ is not given by a complete linear system.
The natural guess for a better embedding is the following:
$$
X' := \overline{
\bigg\{\left[1,g,\Wedge{2} g\right] \in \P^{18} = \P(\C \oplus V)
\mid \det g = 1 \bigg\}},
$$
(we note that $\Wedge{2} g = (g^{-1})^T$ for $g$ with $\det g=1$)
and one can verify that the projection from the point $[1,0,0]\in \P^{18}$
restricted to $X'$ gives an isomorphism with $\Xinv(3)$.

The Grassmannian $Gr(3,6)$ in its Pl\"u{}cker embedding can be described as the closure of:
$$
\bigg\{\left[1,g,\Wedge{2} g, \Wedge{3} g\right] \in \P^{19} = \P(\C \oplus V \oplus \C)
\mid g \in M_{3 \times 3}\bigg\}
$$
and we immediately identify $X'$ as the section $H:=\left\{ \Wedge{3} g=1 \right\}$ of the Grassmannian. 

Though it is not essential, we note that $H^1(\ccO_{Gr(3,6)}) = 0$ 
(see Kodaira vanishing theorem \cite[Thm~4.2.1]{lazarsfeld};
alternatively, it follows from the fact that $b_1=0$ for Grassmannians)
and hence the above embedding of $\Xinv(3)$ 
is given by the complete linear system.

\paragraph{Quadratic equations of $\Xinv(3)$.}

In order to calculate the automorphism group of $\Xinv(3)$
we need to know that the quadratic part of its saturated ideal 
has no other polynomial than linear combinations of 
polynomials \eqref{equations_of_Y1}--\eqref{equations_of_Y4}
 (see also Remark~\ref{remark_ideal_of_Xinv3_is_saturated}).

The ideal of grassmannian $Gr(3,6) \subset \P^{19}$ is generated by 35 quadratic equations, 
which in the above coordinates take form: 
\begin{align*}
[\lambda_0, A,B,\lambda_3] &\in
  \P\left(\End\left(\Wedge{0} \C^3\right) \oplus 
          \End\left(\Wedge{1} \C^3\right) \oplus 
          \End\left(\Wedge{2} \C^3\right) \oplus   
          \End\left(\Wedge{3} \C^3\right)\right)\\ 
  \Wedge{2}A& =  \lambda_0B, \\
  \Wedge{2}B& =  \lambda_3A, \\
  A B^T &=  \lambda_0 \lambda_3 \Id_m, \\  
  B^T A &=  \lambda_0 \lambda_3 \Id_m. 
\end{align*}
Although at the first glance there are 36 equations above, 
the trace equality $\tr (AB^T) = \tr (B^TA)$ makes one of them redundant.

The homogeneous ideal $\I(X')$ is the same but with $\lambda:=\lambda_0 = \lambda_3$ 
and  the ideal of $\Xinv(3)$ arises as the elimination ideal of $\I(X')$ with
the $\lambda$ eliminated. In particular the quadratic polynomials in $\I\left(\Xinv(3)\right)$ 
are exactly those quadrics in $\I(X')$ that do not contain any term with $\lambda$. 
No term in $\lambda$ other than $\lambda^2$ appears more than once in the equations listed above.
Therefore the quadratic part of the ideal $\I\left(\Xinv(3)\right)$ arises by eliminating $\lambda$ from:
\begin{align*}
  A B^T = & \lambda^2 \Id_m, \\  
  B^T A = & \lambda^2 \Id_m. 
\end{align*}
Therefore polynomials \eqref{equations_of_Y} span all the quadratic polynomials in $\Xinv(3)$.

\paragraph{Automorphism group.}

It remains to calculate 
$\Aut\left(\Xinv(3)\right)^0$ 
--- the connected component of the automorphism group.

The tangent Lie algebra of the group of automorphisms of a complex projective 
manifold is equal to the global sections of the tangent bundle, see 
Theorem~\ref{theorem_automorphisms_of_projective_variety}.
A vector field on $\Xinv(3)$ is also a section of $T Gr(3,6)|_{\Xinv(3)}$ 
and we have the following short exact sequence:
$$
0 \lra T Gr(3,6)(-1) \lra T Gr(3,6) \lra T Gr(3,6)|_{\Xinv(3)} \lra 0
$$
The homogeneous vector bundle $T Gr(3,6)(-1)$ is isomorphic to 
$U^* \otimes Q \otimes \Wedge{3} U$, where $U$ is the universal subbundle in  $Gr(3,6) \times \C^6$ 
and $Q$ is the universal quotient bundle.
This bundle corresponds to an irreducible module of the parabolic subgroup in $\Sl_6$.
Calculating explicitly its highest weight and applying Bott formula 
\cite{ottaviani} we get that $H^1\big(T Gr(3,6)(-1)\big) = 0$. 
Hence every section of $T\Xinv(3)$ extends to a section of $T Gr(3,6)$.
In other words, if $P < \Aut(Gr(3,6)) \simeq \P\Gl_6$ 
is the subgroup preserving $\Xinv(3) \subset Gr(3,6)$,
then the restriction map 
\[
P \lra \Aut\bigl(\Xinv(3)\bigr)^0
\]
is epimorphic.

The action of $\Sl_6$ on $\Wedge{3} \C^6$ preserves 
the natural symplectic form $\omega'$:
$$
\omega' : \Wedge{2}\left(\Wedge{3} \C^6\right) \lra  \Wedge{6} \C^6 \simeq \C.
$$
Since the action of $P$ on $\P\left(\Wedge{3} \C^6\right)$ 
preserves the hyperplane $H$ containing $\Xinv(3)$, 
it must also preserve $H^{\perp_{\omega'}}$, 
i.e.~$P$ preserves $[1,0,0,1]\in \P^{19} = \P(\C\oplus V \oplus \C)$.
Therefore $P$ acts on the quotient $H/(H^{\perp_{\omega'}}) = V$ 
and hence the restriction map factorises:
$$
P \lra \Aut\bigl(\P(V), \Xinv(3)\bigr)^0 \epi \Aut\bigl(\Xinv(3)\bigr)^0.
$$

By Theorem~\ref{corollary_that_conjecture_is_true_for_smooth},
group $\Aut(\P(V), \Xinv(3))^0$ 
is contained in the image of $\Sp(V) \to \P\Gl(V)$,
so by Theorem~\ref{theorem_ideal_and_group},
Proposition~\ref{proposition_action_of_G}
and the above considerations about quadratic polynomials:
\[
\Aut\bigl(\P(V), \Xinv(3)\bigr)^0= G.
\]
In particular $\Xinv(3)$ cannot be homogeneous as it contains 
more than one orbit of the connected component of automorphism group.

\end{prf}

We note that the  fact that $\Xinv(3)$ is not homogeneous 
can be also proved without calculating the automorphism group.
Since $\Pic \Xinv(3) \simeq \Z$, 
it follows from Theorem~\ref{theorem_LM_legendrian_homogeneous_b2_one_then_subadjoint},
that $\Xinv(3)$ could only be one of the subadjoint varieties.
But none of them has Picard group $\Z$ and dimension 8.

\subsubsection{Case $m=4$ --- spinor variety $\mathbb{S}_6$}

\begin{prf}
Similarly to Grassmannian $Gr(3,6)$ the spinor variety can be described in terms of matrices.
Let $\Wedge{2} \C^6$ be the vector space of all $6\times 6$ skew-symmetric matrices. 
Consider a map:
\begin{align} \label{equation_parametrisation_of_S6}
\Wedge{2} \C^6 \lra &\P\left(\Wedge{0} \C^6 \oplus 
                             \Wedge{2} \C^6 \oplus
                             \Wedge{4} \C^6 \oplus
                             \Wedge{6} \C^6 \oplus\right) \nonumber\\
g \longmapsto & \left[1, \ g, \ g \wedge g, \  g \wedge g \wedge g \right]. 
\end{align}
In coordinates, $g \wedge g$ is just the matrix of all $4\times 4$ Pfaffians of $g$ and 
$ g \wedge g \wedge g$ is the Pfaffian of $g$. The spinor variety $\mathbb{S}_6$ is the closure of the image of this map.

For a skew-symmetric matrix $g$ by $Pf_{ij}$ denote the Pfaffian of $g$ 
with $i^{\textrm{th}}$ and $j^{\textrm{th}}$ rows and columns removed 
and by $Pf(g)$ denote the Pfaffian of $g$.
Now consider the following map:
\begin{gather}
\Wedge{2} \C^6 \lra \P(V)
\nonumber
\\
g:=\left(
\begin{array}{cccccc}
0        &g_{12}  &g_{13}  &g_{14}  &g_{15}  &g_{16}  \\
-g_{12}  &0       &g_{23}  &g_{24}  &g_{25}  &g_{26}  \\
-g_{13}  &-g_{23} &0       &g_{34}  &g_{35}  &g_{36}  \\
-g_{14}  &-g_{24} &-g_{34} &0       &g_{45}  &g_{46}  \\
-g_{15}  &-g_{25} &-g_{35} &-g_{45} &0       &g_{56}  \\
-g_{16}  &-g_{26} &-g_{36} &-g_{46} &-g_{56} &0       \\
\end{array}
\right) \longmapsto  
\label{equation_parametrisation_of_Xinv(4)}
\\
\left[
\left(
\begin{array}{cccc}
 Pf_{14}  &-Pf_{15}  & Pf_{16}  &g_{23}   \\
-Pf_{24}  & Pf_{25}  &-Pf_{26}  &g_{13}   \\
 Pf_{34}  &-Pf_{35}  & Pf_{36}  &g_{12}   \\
  g_{56}  &  g_{46}  &  g_{45}  &1 
\end{array}
\right),
\left(
\begin{array}{cccc}
 g_{14}  &  g_{15}  & g_{16} & Pf_{23}   \\
 g_{24}  &  g_{25}  & g_{26} &-Pf_{13}   \\
 g_{34}  &  g_{35}  & g_{36} & Pf_{12}   \\
Pf_{56}  &-Pf_{46}  &Pf_{45} &-Pf(g)
\end{array}
\right),
\right].
\nonumber
\end{gather}

The image of this map is in an open neighbourhood of $p_1$ and satisfies the equations of $\Xinv(4)$:
    30 quadratic equations $Y$ as in \eqref{equations_of_Y} 
and 36 quadratic equations from Corollary~\ref{corollary_equations_Xinv} (a).
So the image is contained in $\Xinv(4)$.
Moreover, Map~\eqref{equation_parametrisation_of_Xinv(4)}
is just a linear coordinate change different from Map~\eqref{equation_parametrisation_of_S6}.
Therefore, since $\dim \Xinv(4) = \dim \mathbb{S}_6$ and both varieties are irreducible it follows that 
(up to linear change of coordinates): 
\[
\Xinv(4) = \mathbb{S}_6,
\]
as claimed in the theorem.
\end{prf}

\subsubsection{Case $m \ge 5$ --- singular varieties} 

\begin{prf}
Finally we prove (d).
We want to prove, that for $m\ge 5$ variety $\Xinv(m)$ is singular at $p_1$. 
To do that, we calculate the reduced tangent cone 
$$
T:=\big(TC_{p_1} \Xinv(m)\big)_{red}.
$$ 
From Equations~\eqref{equations_of_Y} we easily get the following 
linear and quadratic equations of $T$
(again we suggest to have a look at \S\ref{notation_submatrices}):
$$
b_{i m} = b_{mi}=0, \quad A_m B_m^T = B_m^T A_m = \lambda^2 \Id_{m-1}
$$
for every $i\in \{1,\ldots m\}$  and some $\lambda\in \C^*$.

Next assume $I$ and $J$  are two sets of indices, both of cardinality $k = \left\lfloor \half m \right\rfloor$
and such that neither $I$ nor $J$ contains $m$. 
Consider the equation of $\Xinv(m)$ as in Corollary~\ref{corollary_equations_Xinv}(b):
$$
\left(\det A_{I,J}\right)^2 = \left(\det B_{I',J'}\right)^2 \cdot (a_{11}b_{11} + \ldots  a_{1m}b_{1m})^l.
$$
To get an equation of $T$, we evaluate at $a_{mm}=1$ and take the lowest degree part,
which is simply 
$\left(\det \left((A_m)_{I,J}\right)\right)^2=0$.
Since $T$ is reduced, by varying $I$ and $J$
we get that: 
$$
\rk A_m \le m - 1 - k - 1 = \left\lceil \half m \right\rceil - 2 
$$
and therefore also:
$$
 A_m B_m^T = B_m^T A_m = 0.
$$

Hence $T$ is contained in the product of the linear space 
$W:=\{A_m=0, B=0\}$ and the affine cone $\hat{U}$ over the union of $\Xdeg(m-1,k)$ for 
$k\le \left\lceil \half m \right\rceil - 2$. We claim that $T= W \times \hat{U}$.
By Proposition~\ref{degenerate_orbits_of_G}(i), every component of $W \times \hat{U}$ 
has dimension $2m-2+(m-1)^2 = m^2-1=\dim \Xinv(m)$,
so by \S\ref{properties_of_tangent_cone}(1) the tangent cone must be a union of some of the components.
Therefore to prove the claim it is enough to find for every $k \le \left\lceil \half m \right\rceil - 2$ 
a single element of $\cDeg^{m-1}_{k, m-k-1}$ that is contained in the tangent cone. 

So take $\alpha$ and $\beta$ to be two strictly  positive integers such that 
$$
\alpha =  \left( \half m - k - 1 \right) \beta
$$
and consider the curve in $\P(V)$ with the following parametrisation:
$$
\left[
\diag\{
\underbrace{t^{\alpha},\ldots, t^{\alpha}}_{k},
\underbrace{t^{\alpha+\beta},\ldots,t^{\alpha+\beta}}_{m-k-1}, 1
\},
\diag\{
\underbrace{t^{\alpha+\beta},\ldots, t^{\alpha+\beta}}_{k},
\underbrace{t^{\alpha},\ldots,t^{\alpha}}_{m-k-1}, t^{2\alpha + \beta}
\}
\right].
$$
It is easy to verify that this family is contained in $\Inv^m$ for $t\ne 0$
and as $t$ converges to $0$, it gives rise to a tangent vector 
(i.e.~an element of the reduced tangent cone --- see the point-wise definition in
\S\ref{properties_of_tangent_cone}) that belongs to $\cDeg^{m-1}_{k, m-k-1}$.

So indeed $T= W \times \hat{U}$, which for $m\ge 5$ contains more than 1 component, 
hence cannot be a linear space. 
Therefore by \S\ref{properties_of_tangent_cone}(3) variety $\Xinv(m)$ is singular at $p_1$.
\end{prf}

\chapter{Hyperplane sections of Legendrian subvarieties}
\label{chapter_hyperplane}

The content of \S\ref{section_hyperplane} and \S\ref{section_extend} 
of this chapter is published in \cite{jabu_hyperplane}.

The Legendrian  variety $\Xinv(3)$ constructed in Chapter~\ref{chapter_sl} 
is isomorphic to a hyperplane section of another Legendrian variety $Gr(3,6)$. 
In this chapter we prove that general hyperplane sections of other Legendrian varieties 
also admit  a Legendrian embedding. 
This gives numerous new examples of smooth Legendrian subvarieties.

\section{Hyperplane section}\label{section_hyperplane}

\subsection{Construction}\label{section_construction}
The idea of
the construction is built on the concept of symplectic reduction 
(see \S\ref{notation_symplectic_reduction}).
Let $H \in \P(V^*)$ be a hyperplane in $V$. 
By 
$$
h:= H^{\perp_{\omega}} \subset V
$$
we denote the $\omega$-perpendicular
to $H$ subspace of $V$, which in this case is a line contained in $H$.
We think of $h$ both as a point in the projective space $\P(V)$ and a line in $V$.
We define 
$$
\pi: \P(H) \backslash \{ h \}  \lra \P(H \slash h)
$$
 to be the projection map and for a given Legendrian subvariety $X\subset \P(V)$ we let
 $\widetilde{X}_H := \pi(X \cap H)$.

We have the  natural symplectic structure $\omega'$  on
$H \slash h$  determined by $\omega$ (see \S\ref{notation_symplectic_reduction}).
Also $\widetilde{X}_H$ is always 
Legendrian by Proposition~\ref{proposition_symplectic_reduction}
and Lemma~\ref{lemma_algebraic_map_submersion}.

Note that so far we have not used any smoothness condition on $X$.

\subsection{Proof of smoothness}\label{section_smooth}

Hence to prove Theorem~\ref{theorem_hyperplane} it is enough to prove that for  a general 
$H \in \P(V^*)$, the map $\pi$ gives an isomorphism of  
the smooth locus of $X \cap H$ onto its image, an open subset in 
$\widetilde{X}_H$.

Recall the definition of secant variety from \S\ref{section_definition_of_dual_variety}.

\begin{lem}
Let $Y\subset \P^m$, choose such a point $y \in \P^m$ that 
$y \notin \sigma(Y)$ and let 
$\pi: \P^m \backslash \{ y \} \lra \P^{m-1}$
be the projection map. 
\begin{itemize}
\item[(a)]
If $Y$ is smooth, then $\pi$ gives an isomorphism of $Y$ and $\pi(Y)$.
\item[(b)]
In general, $\pi$ is 1 to 1 and $\pi$ is an isomorphism of the smooth
	  part of $Y$ onto its image. In particular, the dimension of
	  singular locus of $Y$ is greater or equal to the dimension of
	  singular locus of $\pi(Y)$.
\end{itemize}
\end{lem}

\begin{prf}
See \cite[Prop.~IV.3.4]{hartshorne} 
--- the proof is identical in general as for the curve case.
We only note that if $Y$ is smooth, then the secant variety $\sigma(Y)$ contains all the
 embedded tangent spaces of $Y$. 
They arise when limits $y_2$ approaches $y_1$.
\end{prf}

Now we can prove Theorem~\ref{theorem_hyperplane}:

\begin{prf}
By the lemma and the construction in 
 \S\ref{section_construction} it is enough to prove that there
 exists $h \in \P(V)$ s.t.~$ h \notin \sigma(X\cap h^{\perp_{\omega}})$.

Given two different points $x_1$ and $x_2$ in a projective space we denote
by  $\langle x_1, x_2 \rangle$ the projective line through $x_1$ and $x_2$.
Let
\begin{align*}
\tilde{\sigma}(X)& \subset X\times X \times \P(V),  \\
 \tilde{\sigma}(X)&:=\overline{\{ (x_1, x_2 , p)  | \   p \in \langle x_1, x_2 \rangle  \} },\\
\intertext{%
so that $\tilde{\sigma}(X)$ is the incidence variety for the secant variety of $X$.
Obviously, $\dim(\tilde{\sigma}(X)) = 2 \dim X + 1= \dim(\P(V))$ 
and $\tilde{\sigma}(X)$ is irreducible.
Also we  let:}
\kappa(X) &\subset \tilde{\sigma}(X),  
\\
\kappa(X) &:= \overline{\{ (x_1,
 x_2 , h)  |  \   h \in \langle x_1, x_2 \rangle \textrm{ and }  x_1,
 x_2 \in h^{\perp_{\omega}}  \} },
\end{align*}
so that the image of the projection of $\kappa(X)$ onto the last
 coordinate is the locus of `bad' points. 
More precisely,
for a point $h \in \P(V)$ there exist $(x_1, x_2)$ such that $(x_1,x_2, h) \in \kappa(X)$ 
if and only if 
$h \in \sigma(X \cap h^{\perp_{\omega}})$.

We claim that the image of  $\kappa(X)$ under the projection is not the whole $\P(V)$.
To see this note that the condition defining $\kappa(X)$, i.e.,
$h \in \langle x_1, x_2 \rangle , \  x_1, x_2 \in h^{\perp_{\omega}}$ is equivalent to
$h \in \langle x_1, x_2 \rangle$ and $\langle x_1, x_2\rangle$ is an isotropic subspace of $V$.
Now either $X$ is a linear subspace and then both the claim and the theorem are obvious
or there exist two points $x_1, x_2 \in X$ such that $\omega(\hat{x_1}, \hat{x_2}) \ne 0$
where by $\hat{x_i}$ we mean some non-zero point in the line $x_i \subset V$.
Therefore $\kappa(X)$ is strictly contained in $\tilde{\sigma}(X)$ and 
$$
\dim(\kappa(X)) < \dim(\tilde{\sigma}(X)) = \dim \P(V),
$$
so the image of $\kappa(X)$ under the projection cannot be equal to $\P(V)$%
\footnote{%
The inequality on the dimensions, although simple, is essential for the proof.
An analogous construction for Lagrangian subvarieties in symplectic manifolds 
is known as symplectic reduction (see \S\ref{notation_symplectic_reduction}
for linear algebra baby version of this),
but does not produce smooth Lagrangian subvarieties.
}.
\end{prf}

\begin{cor}\label{corollary_linear_section}
Let $X\subset \P(V)$ be an irreducible Legendrian subvariety whose
 singular locus has dimension at most $k-1$. 
Let $F$ be the contact distribution on $\P(V)$.
If $H\subset \P(V)$ is a general $F$-cointegrable linear subspace of
 codimension $k$, then
$H$ does not intersect the singular locus of $X$ and 
$\widetilde{X}_H: = X \cap H$ is smooth and admits a Legendrian embedding via
an appropriate subsystem of linear system $\ccO_{\widetilde{X}_H}(1)$ into $\P\left(H/H^{\perp_{\omega}}\right)$.
\end{cor}
\noprf

We sketch some proofs of Examples~\ref{many_examples}:

\begin{prf}
$K3$ surfaces of (a) arise as codimension 4 linear sections of
 Lagrangian Grassmannian  $Gr_L(3,6)$. 
Since the canonical divisor $K_{Gr_L(3,6)} = \ccO_{Gr_L(3,6)}(-4)$ (in
 other words $Gr_L(3,6)$ is Fano of index 4),
 by the adjunction formula, the canonical divisor of the section is
 indeed trivial. 
On the other hand, by \cite[Prop.~9]{landsbergmanivel04} it must have
 genus~9.
Although we take quite special (F-cointegrable) sections,
 they fall into the 19 dimensional family of Mukai's $K3$-surfaces of genus 9
 \cite{mukai_K3} and they form a 13
 dimensional subfamily.

The other families of surfaces as in (b) arise as sections of the other
 exceptional subadjoint varieties: $Gr(3,6)$, $\bS_6$ and $E_7$.  
They are all Fano of index 6, 10 and 18 respectively and their
 dimensions are 9, 15 and 27 hence taking successive linear sections we get to Calabi-Yau
 manifolds as stated in (c). 
Further the canonical divisor is very
 ample, so we have examples of general type as stated in (b) and (d).

The Fano varieties arise as intermediate steps, before coming down to
the level of  Calabi-Yau manifolds.
Also $\P^1\times Q^{n}$ is a subadjoint variety and its hyperplane
 section is the blow up of a quadric $Q^n$ in a codimension 2 linear
 section.
The Del Pezzo surfaces are the hyperplane sections of the blow up of $Q^3$ in a conic curve.
\end{prf}

\section{Linear sections of decomposable Legendrian varieties}
\label{linear_section_decomposable}

Assume $m_1$ and $m_2$ are two positive integers, $m_1\ge m_2$. 
Let $V_1\simeq\C^{2m_1+2}$ and $V_2\simeq\C^{2m_2+2}$ be two symplectic vector spaces, 
and let $X_1\subset \P(V_1)$ and $X_2\subset \P(V_2)$ be 
two smooth, irreducible, non-degenerate, Legendrian subvarieties.
In this setup $\dim X_i = m_i$. 
Consider the decomposable variety $X_1 * X_2\subset \P(V_1\oplus V_2)$.
Clearly $\Sing (X_1 * X_2) =X_1 \sqcup X_2$, 
hence $\dim \bigl(\Sing (X_1 * X_2)\bigr) = m_1$, while 
\[
  \dim (X_1 * X_2)=m_1 + m_2+1.
\]
Therefore a general codimension $m_1+1$ $F$-cointegrable linear section of $X_1*X_2$ will be smooth of $\dim m_2$
and admit a Legendrian embedding.
The purpose of this section is to explain that these newly constructed varieties have essentially different properties 
than those of $X_1$ and $X_2$.
Hence our method also allows to produce new examples without dropping the dimension.

Let $L$ be the following line bundle on $X_1 \times X_2$:
\[
L:=\ccO_{X_1}(1) \boxtimes \ccO_{X_2}(-1).
\]
Also let $(X_1* X_2)_0$ be the smooth locus of $X_1*X_2$. 

\begin{lem}\label{lemma_X1*X2_eq_Ldot}
  $(X_1* X_2)_0$ is isomorphic to $\Ldot$, the total space of the 
  $\C^*$-bundle associated to $L$ (see \S\ref{notation_Ldot}).
\end{lem}

\begin{prf}
  Let $\C^*$ act on $V_1\oplus V_2$ with weight $-1$ on $V_1$ and weight $1$ on $V_2$.
  Then
  \[
    \Bigl(\P(V_1 \oplus V_2) \setminus \bigl(\P(V_1)\sqcup \P(V_2)\bigr)\Bigr) \Big/ \C^*
    = \P(V_1) \times \P(V_2)
  \]
  and the quotient map:
  \[
    \Bigl(\P(V_1 \oplus V_2) \setminus \bigl(\P(V_1)\sqcup \P(V_2)\bigr)\Bigr) 
    \stackrel{\Big/ \C^*}{\lra} 
    \P(V_1) \times \P(V_2)
  \]
  is a principal $\C^*$-bundle obtained by removing the zero section 
  from the total space of the line bundle $\ccO_{\P(V_1) \times \P(V_2)}(d_1,d_2)$
  for some integers $d_1$ and $d_2$.
  We have,
  \begin{align*}
    \Pic \Bigl(\P(V_1 \oplus V_2) \setminus \bigl(\P(V_1)\sqcup \P(V_2)\bigr)\Bigr)
    =& \Pic \P(V_1 \oplus V_2) = \Z[\ccO_{\P(V_1 \oplus V_2)}(1)]\\ 
     &\text{(by \cite[Prop.~II.6.5(b)]{hartshorne}).}
    \\
  \intertext{On the other hand,}
    \Pic \Bigl(\P(V_1 \oplus V_2) \setminus \bigl(\P(V_1)\sqcup \P(V_2)\bigr)\Bigr) =& 
    \Pic \bigl(\P(V_1)\times \P(V_2)\bigl) 
    \big/ \left\langle \ccO_{\P(V_1)\times \P(V_2)}(d_1,d_2) \right\rangle \\
     & \text{(by Lemma~\ref{lemma_on_Pic_of_Cstar_bundles}}).
  \end{align*}
  Moreover via the isomorphism
  \[
    \Pic \bigl(\P(V_1)\times \P(V_2)\bigl) 
    \big/ \left\langle \ccO_{\P(V_1)\times \P(V_2)}(d_1,d_2) \right\rangle 
    \simeq 
    \Z[\ccO_{\P(V_1 \oplus V_2)}(1)] 
  \]  
  the class of line bundle $\ccO_{\P(V_1 \oplus V_2)} (e_1,e_2)$ is mapped to 
  $\ccO_{\P(V_1 \oplus V_2)} (e_1+e_2)$. 
  Hence $(d_1,d_2) =(1,-1)$ or $(-1,1)$. 
  In both cases the total spaces of the line bundles 
  are the same after removing the zero sections (the difference is only in 
  the sign of the weights of the $\C^*$-action, which we ignore at this point).

  To finish the proof just note that:
  \[
    (X_1* X_2)_0 = (X_1* X_2) \cap 
    \Bigl(\P(V_1 \oplus V_2) \setminus \bigl(\P(V_1)\sqcup \P(V_2)\bigr)\Bigr)
  \] 
  and the image of $(X_1* X_2)_0$ under the quotient map is equal 
  to $X_1\times X_2$.
\end{prf}
  
Hence by Lemma~\ref{lemma_on_Pic_of_Cstar_bundles} we have:
\[
\Pic \left(X_1\times X_2\right) \epi \Pic \left(X_1* X_2\right)_0 = \Cl \left(X_1* X_2\right)
\]
and the kernel of the epimorphic map is generated by $L$.
If $L_1 \in \Pic X_1$ and $L_2 \in \Pic X_2$, by $[L_1 \boxtimes L_2]$ 
we will denote a line bundle on $(X_1* X_2)_0$ which represents the image 
of $L_1\boxtimes L_2$ under the epimorphic map.

\begin{theo}
  Let $m_1$, $m_2$, $X_1$, $X_2$ be as above. 
  Let $F$ be the contact distribution on $\P(V_1 \oplus V_2)$ and 
  let $H\subset \P(V_1 \oplus V_2)$ 
  be a general $F$-cointegrable linear subspace of codimension $m_1+1$.
  Then 
  $X:=(X_1*X_2) \cap H$ is smooth, admits a Legendrian embedding 
  and has the following properties:
  \renewcommand{\theenumi}{(\alph{enumi})}
  \begin{enumerate}
    \item \label{item_degrees}
      $\deg X = \deg X_1 \cdot \deg X_2$;
    \item \label{item_canonical}
      $K_X \simeq [K_{X_1} \boxtimes K_{X_2}]|_X \otimes \ccO_X(m_1+1)$;
    \item \label{item_Pic}
      We have the restriction map on the Picard groups:
      \[
        i^* : \Pic (X_1 \times X_2) \big/ 
        \left\langle L \right\rangle
        \lra \Pic X.
      \]
      If $m_2 \ge 3$, then $i^*$ is an isomorphism. 
      If $m_2=2$, then $i^*$ is injective.
  \end{enumerate}
  In particular, we have:
  \begin{enumerate}
    \setcounter{enumi}{3}
    \item \label{item_non_projectively_isomorphic}
      $X$ is not projectively isomorphic to neither $X_1$ or $X_2$.
    \item \label{item_canonical_for_subcanonical_embedding}
      If $K_{X_1} \simeq \ccO_{X_1}(d_1)$ and $K_{X_2} \simeq \ccO_{X_2}(d_2)$,
      then $K_X \simeq \ccO_X(d_1+d_2+m_1+1)$;
    \item \label{item_canonical_for_almost_subcanonical_embedding}
       If $K_{X_1}\simeq \ccO_{X_1}(d_1)\otimes E_1$ 
       and $K_{X_2}\simeq \ccO_{X_2}(d_2)\otimes E_2$,
       where the $E_i$'s are line bundles 
       corresponding to some effective divisors, then 
       \[
         K_{X}\simeq \ccO_{X}(d_1+d_2+m_1+1)\otimes E
       \]  
       for some $E$ corresponding to an effective divisor;
    \item \label{item_Pic_for_Picard_number_1}
       If $m_2 \ge 3$, $\Pic X_1=\Z[\ccO_{X_1}(1)]$,
       $\Pic X_2=\Z[\ccO_{X_2}(1)]$ and  either $X_1$ or $X_2$ is simply connected 
       (for example Fano),
       then $\Pic X=\Z[\ccO_{X}(1)]$.
  \end{enumerate}
\end{theo}

\begin{prf}
  Since $X_1$ and $X_2$ are smooth, the singular locus of $X_1*X_2$ has dimension $m_1$.
  By Corollary~\ref{corollary_linear_section}, 
  a general coisotropic linear section of codimension $m_1+1$ (such as $H$) does not intersect the singular locus
  of $X_1*X_2$. Therefore $X$ is smooth and
  has natural Legendrian embedding into $\P\left(H/H^{\perp_{\omega}}\right)$.

  Part~\ref{item_degrees} is immediate, since $\deg (X_1*X_2) = \deg X_1 \cdot \deg X_2$ 
  and neither linear section nor projection changes the degree.
 
  \smallskip
 
  Part~\ref{item_canonical} follows from Lemma~\ref{lemma_X1*X2_eq_Ldot},
  \S\ref{section_on_Pic_of_Cstar_bundles} and the adjunction formula
  (see \cite[Prop.~II.8.20]{hartshorne}).

  \smallskip

  Part~\ref{item_Pic} follows from the following generalisation of Grothendieck-Lefschetz theorem
  due to Ravindra and Srinivas \cite{ravindra_srinivas}:
  \begin{theo}
     Let $Y$ be a subvariety of $\P^m$ with singular locus of codimension at least 2. 
     Let $Y'$ be its general hyperplane section and
     let $Y_0$ and $Y'_0$ be the smooth locus of $Y$ and $Y'$ respectively. 
     Then the restriction map $\Pic(Y_0) \to \Pic(Y'_0)$ is an isomorphism if dimension of $Y'$ is at least $3$
     or is  injective if $\dim Y' = 2$.
  \end{theo}
  Since our variety $X$ arise as repeated hyperplane section and projection of $X_1*X_2$ 
  and the projection part is always an isomorphism of smooth parts, 
  we repeatedly apply the theorem of Ravindra and Sriniva to conclude our claim.

  \smallskip
  Part~\ref{item_non_projectively_isomorphic} follows from \ref{item_degrees}:
  $X_1$ and $X_2$ are not linear space, hence  $\deg X_i >1$. 
  Therefore $\deg X > \deg X_i$ for $i =1,2$ and the varieties cannot be projectively isomorphic.

  \smallskip

  Parts~\ref{item_canonical_for_subcanonical_embedding} 
  and \ref{item_canonical_for_almost_subcanonical_embedding}
  are immediate consequences of \ref{item_canonical} 
  and \ref{item_Pic}, since $[\ccO_{X_1} (d_1) \boxtimes \ccO_{X_2} (d_2)]$
  is isomorphic to $\ccO_{(X_1*X_2)_0} (d_1+d_2)$.

  \smallskip

  Finally, Part~\ref{item_Pic_for_Picard_number_1} follows from \ref{item_Pic} and
  from \cite[Ex.~III.12.6]{hartshorne}.
\end{prf}

To conclude we give a further series of examples:

\begin{ex}
  Apply the theorem to both $X_1$ and $X_2$ equal to the $E_7$-variety.
  As a result we get $X$ which we denote by $(E_7)^{*2}$, 
  a smooth Legendrian Fano variety of dimension 27,
  Picard group generated by a hyperplane section
  and of index 8. Now apply the theorem to $X_1$ being the $E_7$-variety again and
  $X_2=(E_7)^{*2}$. The result, $(E_7)^{*3}$ again has 
  the Picard group generated by a hyperplane section 
  and $K_{(E_7)^{*3}} = \ccO_{(E_7)^{*3}}(2)$, hence is very ample.
  Analogously we construct $(E_7)^{*k}$ and
  combining this result with Corollary~\ref{corollary_linear_section}, 
  we get infinitely many families of smooth Legendrian varieties of general type with 
  Picard group generated by a very ample class in every dimension $d$, 
  where $3\le d \le 27$.
\end{ex}

\begin{ex}
  Let $X_1=\P^1 \times Q^{m_1-1}$ and $X_2$ be arbitrary.
  If $m_1 \ge 3$ and $\dim X_2 \ge 3$, then $X$ has  Picard group isomorphic to 
  $\Pic X_2 \oplus \Z$. 
  Hence we can get a smooth Legendrian variety with arbitrarily big Picard rank.
\end{ex}

\begin{ex}
  Let $X_1=X_2=\P^1 \times Q^{m-1}$. 
  Let the resulting $X$ be called $(\P^1 \times Q^{m-1})^{*2}$.
  Then $K_{X_i}=\ccO_{X_i}(-m)\otimes E_i$, where $E_i$ is effective.
  Hence 
  \[
    K_{(\P^1 \times Q^{m-1})^{*2}} = \ccO_{(\P^1 \times Q^{m-1})^{*2}}(-m+1) \otimes E
  \]
  for an effective $E$.
  Construct analogously $(\P^1 \times Q^{m-1})^{*k}$ by taking the section of 
  \[
    \Bigl(\left(\P^1 \times Q^{m-1}\right)^{*(k-1)}\Bigr) * \Bigl(\P^1 \times Q^{m-1}\Bigr).
  \]
  We get that 
  \[
    K_{(\P^1 \times Q^{m-1})^{*k}} = \ccO_{(\P^1 \times Q^{m-1})^{*k}}(-m-1+k) \otimes E
  \] 
  and for $k>m+1$ we get that the canonical divisor can be written 
  as an ample plus an effective, so it is big.
  Hence in every dimension, it is possible to construct many smooth Legendrian varieties
  with the maximal Kodaira dimension.
\end{ex}

\section{Extending Legendrian varieties}\label{section_extend}

Our motivation is the example of Landsberg and Manivel \cite[\S4]{landsbergmanivel04},
a Legendrian embedding of a Kummer K3 surface blown up in 12 points.
It can be seen, that this embedding is given by a codimension 1 linear system.
We want to find a Legendrian 3-fold in $\P^7$ whose hyperplane section is this example.
Unfortunately, we are not able to find a smooth 3-fold with these properties,
but we get one with only isolated singularities.

We recall the setup for the construction of the example.
Let $Z^\sharp \subset \P(T^*\P^n) \subset \P^n \times \check{\P}^n$ be 
\textbf{the conormal variety}, i.e., the closure of the union of 
projectivised conormal spaces over smooth points of $Z$.
Landsberg and Manivel study in details an explicit birational map
\[
\varphi:=\varphi_{H_0, p_0}:\P(T^*\P^n) \dashrightarrow \P^{2n-1}
\]
which depends on a hyperplane $H_0$ in $\P^n$ and on a point $p_0 \in H_0$.
After Bryant \cite{bryant} they observe that $\overline{\varphi(Z)}$ 
(if only makes sense) is always a Legendrian subvariety, but usually singular.
Next they study conditions under which $\overline{\varphi(Z)}$ is smooth.
In particular, they prove that  the conditions are satisfied
 when $Z$ is a Kummer quartic surface in $\P^3$ 
in general position with respect to $p_0$ and $H_0$
and this  gives rise to their example.

We want to modify the above construction just a little bit to obtain our 3-fold.
Instead of considering $Z^\sharp$ as a subvariety in 
$$
\P(W)\times \P(W^*) = (W\setminus\{0\})\times(W^*\setminus\{0\}) \slash \C^* \times \C^*,
$$
we consider a subvariety $X$ in
$$
\P^{2n+1} = \P(W\oplus W^*) = (W\times W^*) \setminus \{0\} \slash \C^*
$$
such that the underlying affine cone of $X$ in $W\times W^*$ is the same 
as the underlying affine pencil of $Z^\sharp$.
In other words, we take $X$ to be the closure of preimage of $Z^\sharp$
under the natural projection map:
$$
p: \P(W\oplus W^*) \dashrightarrow \P(W)\times \P(W^*).
$$

Both $\P(W)$ and $\P(W^*)$ are naturally embedded into $\P(W\oplus W^*)$. 
Let $H$ be a hyperplane in $\P(W\oplus W^*)$ which does not contain $\P(W)$ nor $\P(W^*)$.
Set $H_0:= \P(W) \cap H$ and $p_0$ to be the point in $\P(W)$ dual to
$\P(W^*) \cap H$. 
Assume $H$ is chosen in such a way that $p_0 \in H_0$.

\begin{theo}\label{theorem_extending}
  Let $X\subset \P(W\oplus W^*) \simeq \P^{2n+1}$ be a subvariety constructed as above 
  from any irreducible subvariety $Z\subset \P(W)$.  
  On $W\oplus W^*$ consider the standard symplectic structure 
  (see \S\ref{section_symplectic_on_W_plus_W_dual}) 
  and on $\P(W\oplus W^*)$ consider the associated contact structure.
  Also assume $H$, $H_0$ and $p_0$ are chosen as above.
  Then:
  \begin{itemize}
    \item[(i)]
      $X$ is a Legendrian subvariety contained in the quadric 
      $\overline{p^{-1} \Bigl(\P\bigl(T^*\P(W)\bigr)\Bigr)}$.
    \item[(ii)]
      Let $\widetilde{X}_H$ be the Legendrian variety in $\P^{2n-1}$ 
      constructed from $X$ and $H$ as in \S\ref{section_construction}.
      Also consider the closure of $\varphi_{H_0, p_0}(Z^\sharp)$ 
      as in the construction of \cite[\S4]{landsbergmanivel04}. 
      Then the two constructions agree, 
      i.e., the closure $\varphi_{H_0, p_0}(Z^\sharp)$ is a component of $\widetilde{X}_H$.
    \item[(iii)]
      The singular locus of $X$ equal to the union of following:
        \subitem 
          on $\P(W)$ the singular points of $Z$,
        \subitem 
          on $\P(W^*)$ the singular points of $Z^*$ and 
        \subitem 
          outside $\P(W) \cup \P(W^*)$ the preimage under $p$ 
          of the singular locus of the conormal variety $Z^\sharp$.
  \end{itemize}
\end{theo}

\begin{prf}
For part (i)
consider $\widehat{Z}\subset W$, the affine cone over $Z\subset \P(W)$.
The cotangent bundle to $W$ is equal to  $W\oplus W^*$. 
Furthermore, by our definition $\widehat{X}\subset V$, 
the affine cone over $X\subset \P(W \oplus W^*)$ 
is the conormal variety of $\widehat{X}$, so a Lagrangian subvariety 
(see  Example~\ref{example_conormal_Lagrangian}).

\smallskip

For part (ii), we choose coordinates $x_0, x_1,\ldots, x_n$ on $W$ and dual coordinates $y^0, y^1,\ldots, y^n$ on $W^*$
such that in the induced coordinates on $V$ the hyperplane $H$ has the equation $x_0 - y^n = 0$. 
Now restrict to the affine piece $x_0 = y^n=1$ on both $H$ and $\P(W) \times \P(W^*)$. 
We see explicitly, that the projection map $H\rightarrow \P^{2n-1}$ 
$$
[1, x_1,\ldots, x_n, y^0,\ldots, y^{n-1}, 1] \mapsto [y^1, \ldots, y^{n-1}, y^0 - x_n, x_1,\ldots, x_{n-1},1]
$$
agrees with the map $\varphi$ from \cite[\S4]{landsbergmanivel04}.

\smallskip

To find the singularities of $X$ on $X\cap \P(W)$ as in part (iii) note that 
$X$ 
is invariant under the following action of $\C^*$ on $\P(W \oplus W^*)$: 
$$
 t \cdot [w, \alpha] := [t w, t^{-1}\alpha].
$$
In particular, points of $X\cap \P(W)$ are fixed points of the action. 
So let $[w,0] \in X$ and then $T_{[w,0]} X$ decomposes into the eigenspaces of the action:
\begin{equation}\label{decomposition_of_TX}
T_{[w,0]} X = T_{[w,0]} (X \cap \P(W)) \oplus T_{[w,0]} (X \cap F_w)
\end{equation}
where $F_w$ is the fibre of the projection $\rho: (\P(W\oplus W^*) \backslash \P(W^*)) \rightarrow \P(W)$,
$F_w:=\rho^{-1}([w])$.
Clearly the image of $X$ under the projection $\rho$ is $Z$, 
so the dimension of a general fibre of $\rho|_X: X \rightarrow Z$ is equal to 
$\dim X - \dim Z = \dim \P(W)  - \dim Z = \codim_{\P(W)} Z$. 
Therefore, since the dimension of the fibre can only grow at special points, we have:
\begin{equation}\label{equation_on_codim_Z}
\dim  T_{[w,0]} (X \cap F_w) \ge \dim(X \cap F_w) \ge \codim_{\P(W)} Z.
\end{equation}
Also $\uD_{[w,0]} (\rho|_X): T_{[w,0]} X \rightarrow T_{[w]} Z$ 
maps $T_{[w,0]} (X \cap F_w)$ to 0 and $ T_{[w,0]} (X \cap \P(W))$ onto $T_{[w]} Z$.
Therefore:
\begin{equation}\label{equation_on_dim_Z}
\dim T_{[w,0]} (X \cap \P(W)) \ge \dim T_{[w]} Z \ge \dim Z. 
\end{equation}
Now assume $[w,0]$ is a smooth point of $X$. 
Then adding Equations~\eqref{equation_on_codim_Z} and \eqref{equation_on_dim_Z} we get:
\begin{multline*}
\dim X 
= \dim T_{[w,0]} X 
\\
\stackrel{\textrm{by \eqref{decomposition_of_TX}}}{=}
\dim  T_{[w,0]} (X \cap F_w) + \dim T_{[w,0]} (X \cap \P(W)) 
\\
\ge \codim_{\P(W)} Z + \dim Z 
= \dim \P(W).
\end{multline*}
By (i) the $\dim X$ is equal to the $\dim \P(W)$, 
so in Equations~\eqref{equation_on_codim_Z} and \eqref{equation_on_dim_Z} 
all the inequalities are in fact equalities. 
In particular $\dim T_{[w]} Z = \dim Z$, 
so $[w]$ is a smooth point of $Z$. 

Conversely, assume $[w]$ is a smooth point of $Z$, then the tangent space 
$$
T_{[w,0]} X = T_{[w]}Z \oplus N^*_{[w]}(Z\subset \P(W)),
$$ 
therefore clearly $[w,0]$ is a smooth point of $X$.

Exactly the same argument shows that 
$Z^*$ is singular at $[\alpha]$
if and only if 
$X$ is singular at $[0,\alpha] \in X \cap \P(W^*)$.

For the last part of (iii) it is enough to note that 
$p$ is a locally trivial $\C^*$-bundle when restricted to 
$\P(W\oplus W^*) \backslash \left(\P(W) \cup \P(W^*)\right)$.
\end{prf}

\begin{cor}
Given a Legendrian subvariety $\widetilde{Z} \subset \P^{2n-1}$ we can take
\[
  Z^{\#}:=\varphi_{H_0,p_0}^{-1} (\widetilde{Z})
\]
to construct a Legendrian subvariety in $\P(T^* \P^n)$. 
Such a variety must be the conormal variety to some variety $Z \subset \P^n$
(see Corollary~\ref{cor_Legendrian_in_cotangent}).
Let $X\subset \P^{2n+1}$ be the Legendrian variety constructed from $Z$ as above.
By Theorem~\ref{theorem_extending}(ii), a component of a hyperplane
 section of $X$ can be
 projected onto $\widetilde{Z}$. 
\end{cor}

Unfortunately, in the setup of the theorem $X$ is almost always singular 
(see \S\ref{section_self_dual}).

\begin{ex}
\label{example_LM}
If $Z$ is a Kummer quartic surface in $\P^3$, then $X$ is a 3-fold with
 32 isolated singular points 
(it follows from Theorem~\ref{theorem_extending}(iii)
because the Kummer quartic surface has 16 singular points, 
it is isomorphic to its dual 
and it has smooth conormal variety in $\P(T^* \P^3)$). 
Therefore by Theorem~\ref{theorem_hyperplane} 
a general hyperplane section of $X$ is smooth and admits a Legendrian embedding.
By Theorem~\ref{theorem_extending} the example 
of Landsberg and Manivel is a special case of this hyperplane section. 
Even though the condition $p_0 \in H_0$ is a closed condition, 
it satisfies the generality conditions of  
Theorem~\ref{theorem_hyperplane} and therefore this hyperplane section consists
of a unique smooth component that is projected isomorphically onto
$\widetilde{Z}$.
\end{ex}

\begin{ex}
 \label{example_bryant}
Similarly, if $Z$ is a curve in $\P^2$ satisfying the generality conditions
 of Bryant \cite[Thm~G]{bryant}, then $X$ is a surface with only
 isolated singularities and its hyperplane section projects
 isomorphically onto a  Bryant's Legendrian curve. 
\end{ex}

\section{Smooth varieties with smooth dual}
\label{section_self_dual}

Furthermore we observe that a classical problem of classifying smooth varieties 
with smooth dual variety can be expressed in terms of Legendrian varieties:

\begin{cor}\label{corollary_self_dual}
  Using the notation of the previous section, 
  let $Q_W \subset \P(W\oplus W^*)$ be the quadric 
  $\overline{p^{-1} \left(\P(T^*\P(W))\right)}$
   --- see Theorem~\ref{theorem_extending}(i). 
  On $W\oplus W^*$ consider the standard symplectic structure 
  (see \S\ref{section_symplectic_on_W_plus_W_dual}) 
  and on $\P(W\oplus W^*)$ consider the associated contact structure
  (see \S\ref{section_projective_space}).
  \renewcommand{\theenumi}{(\roman{enumi})}
  \begin{enumerate}
    \item \label{item_self_dual_1}
      Let $Z\subset \P(W)$ be a smooth subvariety with $Z^* \subset \P(W^*)$ smooth.
      Let $X\subset \P(W\oplus W^*)$ be as in the above construction. 
      Then $X$ is a smooth Legendrian variety contained in $Q_W$.
    \item \label{item_self_dual_2}
      Conversely, assume $X\subset \P(W\oplus W^*)$ is irreducible, Legendrian  
      and contained in $Q_W$.
      Let $Z=X\cap \P(W)$.
      Then $Z^*= X\cap \P(W^*)$ and the variety arising from $Z$ in the above 
      construction is exactly $X$.
      Moreover, if $X$ is smooth, then $Z$ and $Z^*$ are smooth.
  \end{enumerate}
\end{cor}

We underline that although all the smooth quadrics of a given dimension 
are projectively isomorphic,
the classification of quadrics relatively to the contact structure is more complicated.
The quadric $Q_W$ can therefore be written as $x_0 y_0 +\ldots + x_n y_n=0$ 
in some \underline{symplectic} coordinates $x_0, \ldots, x_n, y_0, \ldots, y_n$ on $W \oplus W^*$.
We note (without proof), that such quadric $Q_W$ 
determines uniquely the pair of Lagrangian subspaces $W$ and $W^*$.

\begin{prf}
  Part~\ref{item_self_dual_1} follows immediately from
  Theorem~\ref{theorem_extending}(i) and (iii).
  
  \smallskip

  To prove Part~\ref{item_self_dual_2}, consider $p(X) \subset \P(T^*\P(W))$.
  By Lemma~\ref{lemma_algebraic_map_submersion} 
  and Proposition~\ref{proposition_symplectic_reduction} 
  $p(X)$ is Legendrian. 
  By Corollary~\ref{cor_Legendrian_in_cotangent},  
  $p(X)$ is a conormal variety to some subvariety $Z \subset \P(W)$.
  The next thing to prove is that $X$ coincides with the variety 
  constructed above from $Z$, i.e.~that 
  \[
    X= \overline{p^{-1}\bigl(p(X)\bigr)}.
  \]
  Equivalently, it is enough to prove that $X$ is 
  $\C^*$-invariant. 
  This is provided by Theorem~\ref{theorem_ideal_and_group} 
  since the quadric $Q_W$ produces exactly the required action.
  Finally, it follows that $Z= X \cap \P(W)$. 
  Moreover, $p(X)$ is also the conormal variety to $Z^* \subset \P(W^*)$ and hence
  $Z^*= X \cap \P(W^*)$.
  If $X$ is in addition smooth, 
  then $Z$ and $Z^*$ are smooth by Theorem~\ref{theorem_extending}(iii).
\end{prf}

\begin{table}[hbt]
\begin{tabular}{|p{0.4\textwidth}|p{0.4\textwidth}|}
\hline
smooth self-dual variety $Z\subset \P^n$ & 
the corresponding Legendrian variety $X\subset \P^{2n+1} $ \\
\hline
$Q^{m}$ & $\P^1\times Q^{m}$\\
$\P^1 \times \P^{m}$ & $\P^1 \times Q^{2m}$\\
$Gr(2,5)$ & $Gr(3,6)$\\
$\bS_5$    & $\bS_6$\\
\hline
\end{tabular}
\caption{The known self-dual varieties and their corresponding Legendrian varieties.
Note that $Q^{2m}$ and $\P^1\times \P^{m}$ lead to isomorphic Legendrian varieties. 
Yet their embeddings in the distinguished quadrics are not isomorphic.
\label{table_dual}}
\end{table}

Therefore the classification of smooth varieties with smooth dual
is equivalent to the classification of pairs $(X,Q)$, 
where $Q\subset \P^{2n+1}$ is a quadric which can be written as $x_0 y_0 +\ldots + x_n y_n=0$ 
in some symplectic coordinates 
\[
  x_0, \ldots, x_n, y_0, \ldots, y_n
\]
 on $\C^{2n+2}$ 
and $X\subset \P^{2n+1}$ is a smooth Legendrian variety, which is contained in $Q$.
So far the only known examples of smooth varieties with smooth dual 
are the smooth self-dual  varieties (see \cite{ein}). 
From these we get some of the homogeneous Legendrian varieties (see Table~\ref{table_dual}).
Therefore we cannot hope to produce new examples 
of smooth Legendrian varieties in this way. 
What we hope for is to classify the pairs $(X,Q)$ as above 
and hence finish the classification of smooth varieties with smooth dual.

\newpage
\rhead[\fancyplain{}{}]{\fancyplain{}{\footnotesize\textbf{Bibliography}}}

\bibliography{doktorat}

\bibliographystyle{alpha}

\end{document}